\magnification=\magstep1
\input amstex
\documentstyle{amsppt}

\define\defeq{\overset{\text{def}}\to=}

\def \isom {\buildrel \sim \over \rightarrow}

\def \Pri{\operatorname {Pri}}
\def \Div{\operatorname {Div}}
\def \Inn{\operatorname {Inn}}
\def \perf{\operatorname {perf}}
\def \lim{\operatorname {lim}}
\def \sep{\operatorname {sep}}
\def \rig{\operatorname {rig}}
\def \Isom{\operatorname {Isom}}
\def \Spec{\operatorname {Spec}}
\def \Gal{\operatorname {Gal}}
\def \pr{\operatorname {pr}}
\def \Ker{\operatorname {Ker}}
\def \id{\operatorname {id}}
\def \acts\ trivially\ on{\operatorname {acts\ trivially\ on}}
\def \closed\ subgroups{\operatorname {\closed\ subgroups}}
\def \Hom{\operatorname {Hom}}
\def \t{\operatorname {t}}
\def \ur{\operatorname {ur}}

\def \cd{\operatorname {cd}}

\def \ab{\operatorname {ab}}
\def \tor{\operatorname {tor}}
\def \Aut{\operatorname {Aut}}
\def \Im{\operatorname {Im}}
\def \Norm{\operatorname {Norm}}
\def \Frob{\operatorname {Frob}}
\def \Re{\operatorname {Re}}

\def \open{\operatorname {open}}
\def \pr{\operatorname {pr}}
\def \inrig{\operatorname {inrig}}

\def \et{\operatorname {et}}
\def \Pic{\operatorname {Pic}}
\def \cc{\operatorname {c\text{-}cn}}
\def \ord{\operatorname {ord}}
\def \Supp{\operatorname {Supp}}

\def \Sub{\operatorname {Sub}}

\def\Dec{\operatorname{Dec}}
\def\frakIm{\frak I \frak m}
\def\End{\operatorname{End}}
\def\adm{\operatorname{adm}}
\def\div{\operatorname{div}}
\define\w{\operatorname{w}}

\NoRunningHeads
\NoBlackBoxes
\topmatter

\title
On the Hom-form of Grothendieck's birational anabelian
conjecture 
in characteristic $p>0$
\endtitle


\author
Mohamed Sa\"\i di and Akio Tamagawa
\endauthor





\abstract
We prove that a certain class of open homomorphisms between Galois groups of 
function fields of curves over finite fields arise from embeddings between
the function fields. 
\endabstract

\toc

\subhead
\S 0. Introduction
\endsubhead

\subhead
\S 1. Generalities on Galois groups of function fields of curves
\endsubhead

\subhead
\S 2. Basic properties of homomorphisms between Galois groups
\endsubhead

\subhead
\S 3. Rigid homomorphisms between Galois groups
\endsubhead

\subhead
\S 4. Proper homomorphisms between Galois groups
\endsubhead

\subhead 
\S 5. Recovering the additive structure
\endsubhead

\endtoc

\endtopmatter

\document            

\subhead
\S 0. Introduction
\endsubhead

Let $K$ be an infinite field which is finitely generated 
over its prime field. Let $\overline K$ be an algebraic closure of $K$ 
and $K^{\sep}$ (resp. $K^{\perf}$) the separable closure (resp. perfection) of $K$ 
in $\overline K$. 
Let $G_K\defeq \Gal (K^{\sep}/K)$ be the absolute Galois group of $K$. 
(Observe $G_K=G_{K^{\perf}}$.) 
The ultimate aim of Grothendieck's birational anabelian conjectures 
is to reconstruct the field structure of $K$ from the topological group 
structure of $G_K$. More precisely, these conjectures can be formulated as follows.

\definition
{Birational anabelian conjectures} There exists a  
group-theoretic recipe in order to recover finitely 
generated infinite fields $K$ (or their perfections 
$K^{\perf}$) from their absolute Galois groups $G_K$. 
In particular, if for such fields $K$ and $L$ one has $G_K\isom G_L$, then 
$L^{\perf}\isom K^{\perf}$.
Moreover, given two such fields $K$ and $L$ one has the following.
\enddefinition

\definition
{Isom-form} Every isomorphism $\sigma:G_K\isom G_L$ is defined by a 
field isomorphism $\bar\gamma:\overline L \isom \overline K $, and $\bar\gamma$ is unique 
(resp. unique up to Frobenius twists) if the characteristic is $0$ (resp. $>0$). 
In particular, $\bar\gamma$ induces an isomorphism $L\isom K$. 
\enddefinition

\definition
{Hom-form}  Every open homomorphism $\sigma:G_K\to G_L$ is defined 
by a field embedding $\bar\gamma:\overline L \hookrightarrow \overline K$, and $\bar\gamma$ 
is unique (resp. unique up to Frobenius twists) 
if the characteristic is $0$ (resp. $>0$). 
In particular, $\bar\gamma$ induced a field embedding $L^{\perf}\hookrightarrow K^{\perf}$. 
\enddefinition

Thus, the Hom-form is stronger than the Isom-form. 
The first results concerning these 
conjectures were obtained by Neukirch and Uchida in the case of 
global fields. 

\proclaim {Theorem (Neukirch, Uchida)} Let $K$ and $L$ be global fields.
Then the natural map 
$$\Isom
(L,K)\to \Isom  (G_K,G_L)/\Inn (G_L)$$
is a bijection.
\endproclaim

More precisely, this is due to Neukirch and Uchida 
for number fields ([Neukirch1], [Neukirchi2], [Uchida1]), 
and due to Uchida for function fields of curves over 
finite fields ([Uchida2]). 
Later, Pop generalized the results of Neukirch and Uchida to the case of 
finitely generated fields of higher transcendence degree
([Pop1], [Pop3], see also [Szamuely] for a survey on Pop's results). 

In characteristic $0$, Mochizuki proved the following 
relative version of the Hom-form of the birational 
conjectures (cf. [Mochizuki1]). 

\proclaim {Theorem (Mochizuki)} Let $K$ and $L$ be two finitely 
generated, regular extensions of a field $k$. Assume that  
$k$ is a sub-$p$-adic field (i.e., $k$ can be embedded in a 
finitely generated extension of $\Bbb Q_p$) for some prime number $p$. 
Then the natural map
$$\Hom_{k}
(L,K)\to \Hom _{G_k}^{\open}  (G_K,G_L)
/\Inn (\Ker(G_L\to G_k))$$
is a bijection. Here, 
$\Hom_{k}
$ denotes the set of $k$-embeddings,
and $\Hom _{G_k}^{\open}$ denotes the set of open $G_k$-homomorphisms.
\endproclaim

\noindent
However, almost nothing is known about 
the absolute version (i.e., not relative with respect to a fixed base 
field $k$) of the Hom-form, except for Uchida's result [Uchida3] 
for $K=\Bbb Q$ and $[L:\Bbb Q]<\infty$. 

One of the major obstacles in proving the 
Hom-form of the birational anabelian conjectures 
is that one of the main common ingredients in the 
proofs of Neukirch, Uchida, and 
Pop, which is the so-called local theory (or Galois 
characterization of the decomposition subgroups) and which 
is used in order to establish a one-to-one correspondence between 
divisorial valuations, 
is not available 
in the case of open homomorphisms 
between Galois groups. More precisely, the main result of local theory 
available so far (cf. Proposition 1.5) gives very little information on the 
image of the decomposition subgroups in this case, though one can still
prove some 
partial results (cf. Proposition 2.2, Lemma 2.6, and Lemma 2.9). It seems 
quite difficult, for the moment, to establish a 
satisfactory 
local theory which is suitable to the 
Hom-form of the above conjecture. 
Also, the methods used in the proof of Mochizuki's above theorem
are quite different, and do not rely on local theory. Instead, Mochizuki 
proves his result as an application of his fundamental anabelian result 
that relative open homomorphisms between arithmetic fundamental groups of 
curves over sub-$p$-adic fields arise from  morphisms between corresponding 
curves, whose proof relies on $p$-adic Hodge theory. It is not clear how to 
use, or adapt, Mochizuki's method to the case of positive characteristics.

In this paper we investigate the Hom-form of the birational anabelian 
conjectures for function fields of curves 
over finite fields. For $i=1,2$, let $k_i$ be a finite field. 
Let $X_i$ be a proper, smooth, geometrically connected algebraic curve 
over $k_i$. 
Let 
$K_i
$ 
be 
the function field of $X_i$ 
and fix an algebraic closure $\overline K_i$ of $K_i$. 
Let $K_i^{\sep}$ (resp. $K_i^{\perf}$) 
be 
the separable closure (resp. perfection) of $K_i$
in $\overline K_i$, and  
$\bar k_i$ the algebraic closure of $k_i$ in $\overline K_i$. 
Write $G_i\defeq G_{K_i}\defeq \Gal (K_i^{\sep}/K_i)$ for the 
absolute Galois group of $K_i$, and $G_{k_i}\defeq \Gal (\bar k_i/k_i)$
for the absolute Galois group of $k_i$. We have the following natural exact 
sequence of profinite groups:
$$1\to \overline G_i\to G_i@>\pr_i>> G_{k_i}\to 1,$$
where $\overline G_i$ is the absolute Galois group $\Gal (K_i^{\sep}/
K_i\bar k_i)$ of $K_i\bar k_i$, and $\pr_i$ is the canonical projection.

Further, let $p_i$ be the characteristic of $k_i$, and let $\overline G_i^{p_i'}$ 
be the maximal prime-to-$p_i$ quotient of $\overline G_i$. The push-forward
of the above sequence with respect to the natural surjection 
$\overline G_i\to \overline G_i^{p_i'}$ gives rise to the following natural 
exact sequence
$$1\to \overline G_i^{p_i'}\to G_i^{(p_i')}@>\pr_i>> G_{k_i}\to 1.$$

Set $\frak G_i\defeq G_i$, $i=1,2$ or 
$\frak G_i\defeq G_i^{(p_i')}$, $i=1,2$, 
and refer to the first and the second cases as 
the profinite and the prime-to-characteristic cases, respectively. 
In this paper we investigate two classes of continuous, open 
homomorphisms 
--- rigid homomorphisms and proper homomorphisms --- 
between $\frak G_1$ and $\frak G_2$. 

First, we investigate a class of continuous, open homomorphisms
$\sigma: \frak G_1\rightarrow \frak G_2$, which we call rigid. 
More precisely, we say that $\sigma$ is strictly rigid if 
the image of each decomposition subgroup of $\frak G_1$ coincides 
with a decomposition subgroup of $\frak G_2$, and we say that 
$\sigma$ is rigid if 
there exist open subgroups $\frak H_1\subset \frak G_1$,  
$\frak H_2\subset \frak G_2$, such that $\sigma(\frak H_1)\subset\frak H_2$ 
and that $\frak H_1\overset{\sigma}\to{\to}\frak H_2$ is strictly rigid. 
Thus, isomorphisms between  
$\frak G_1$and $\frak G_2$ are rigid by the main 
result of local theory for the Isom-form. 
Let $\Hom (\frak G_1,\frak G_2)^{\rig}$ be the set 
of rigid homomorphisms between $\frak G_1$ and $\frak G_2$.

We say that 
a homomorphism $\gamma: K_2\to K_1$ of fields (which defines 
an extension $K_2/K_1$ of fields) is admissible, 
if the extension $K_1/K_2$ appears in the extensions of $K_2$ corresponding to 
the open subgroups of $\frak G_2$, or, equivalently,
in the profinite (resp. prime-to-characteristic) case, 
if the extension $K_1/K_2$ is separable (resp. if the extension 
$K_1/K_2$ is separable and the Galois closure of the extension 
$K_1\bar k_1/K_2\bar k_2$ is of degree prime to $p\defeq p_1=p_2$). 
We define $\Hom(K_2,K_1)^{\adm}\subset\Hom(K_2,K_1)$ 
to be the set of admissible homomorphisms $K_2\to K_1$. 

Now, our first main result is the following (cf. Theorem 3.4).

\proclaim {Theorem A} 
The natural map $\Hom(K_2,K_1)\to\Hom(\frak G_1, \frak G_2)/\Inn(\frak G_2)$ 
induces a bijection 
$$\Hom(K_2,K_1)^{\adm}\isom\Hom(\frak G_1, \frak G_2)^{\rig}/\Inn(\frak G_2).$$ 
%
%
\endproclaim

Our method to prove 
Theorem A 
is 
as follows. 
First, we prove, using 
a certain weight argument based on the Weil conjecture for curves, that
a strictly rigid homomorphism  $\sigma:\frak G_1\rightarrow \frak G_2$ 
induces a bijection $\Sigma _{X_1}\isom \Sigma _{X_2}$
between the set of closed points of $X_1$ 
and $X_2$ (cf. Lemma 3.8). 
\comment
Second, using class field theory, we reconstruct
in a functorial way an embedding $\gamma:K_2^{\times}
\hookrightarrow K_1^{\times}$ between multiplicative groups, up to 
Frobenius twist (cf. Lemma 3.9 ** deleted **). Moreover, this embedding is 
order-preserving, and value-preserving with respect to the 
above identification  $\Sigma _{X_1}\simeq \Sigma _{X_2}$ 
(cf. Definition 5.1 and Definition 5.2). Finally, we show that this embedding 
$\gamma:K_2^{\times}\hookrightarrow K_1^{\times}$ is additive. 
Recovering the additive structure is one of the main steps in the proof. 
This problem was treated by Uchida in the case of a bijective 
identification $K_2^{\times}\isom K_1^{\times}$ between multiplicative 
groups, which is order-preserving, and value-preserving. 
In fact, one needs only to restore the additivity 
between constants. For this one has to show identities of the form
$\gamma (f_2+1)=\gamma (f_2)+1$ for some specific 
non-constant function $f_2\in K_2$. 
Uchida succeeded in his case by choosing $f_2$ to be
a function with a minimal pole divisor (He called such a function a 
minimal element). His argument fails in the case of an embedding 
between multiplicative group which is not surjective, because the image 
of a minimal element is not necessarily minimal in this case.
Roughly speaking, we extend his arguments by 
using, instead, a function $f_2$ which has a unique pole. 
This one pole argument turns out to be very efficient, and leads to the
recovery of the additive structure under quite general assumptions 
(cf. Proposition 5.3).
\endcomment
By using this, we can reduce the Hom-form in this case to the Isom-form, 
which has been established in [Uchida2] (profinite case) and [ST] 
(prime-to-characteristic case). 

Next, we investigate a class of continuous, open homomorphisms
$\sigma:\frak G_1\rightarrow \frak G_2$, 
which we call proper. These are homomorphisms with the property that the image of each 
decomposition subgroup of $\frak G_1$ coincides with an open subgroup of 
a decomposition subgroup of $\frak G_2$, and such that each 
decomposition subgroup of $\frak G_2$ contains images of only finitely many 
conjugacy classes of decomposition subgroups of $\frak G_1$. We also consider
a certain rigidity condition (called ``inertia-rigidity'') 
on the various identifications between the 
modules of the roots of unity (cf. Definition 4.5). Unfortunately, we are not able 
to prove that this condition automatically holds for proper homomorphisms.
Let $\Hom(\frak G_1,\frak G_2)^{\pr,\inrig}$ be the set of 
proper and inertia-rigid homomorphisms between $\frak G_1$ and $\frak G_2$.
Our second main result is the following (cf. Theorem 4.8).

\proclaim
{Theorem B} 
%
The natural map $\Hom(K_2,K_1)\to\Hom(\frak G_1, \frak G_2)/\Inn(\frak G_2)$ 
induces a bijection 
$$\Hom(K_2,K_1)^{\sep}\isom\Hom(\frak G_1, \frak G_2)^{\pr,\inrig}/\Inn(\frak G_2).$$ 
Here, we define 
$\Hom(K_2,K_1)^{\sep}\subset\Hom(K_2,K_1)$ to be the set of 
separable homomorphisms $K_2\to K_1$. 
%
%
\endproclaim

In order to prove Theorem B, we first show, using a  
weight argument, that a homomorphism
$\sigma:\frak G_1\rightarrow \frak G_2$ as above induces a surjective map 
$\Sigma _{X_1}\to \Sigma _{X_2}$ between the sets of 
closed points of $X_1$ and $X_2$, which has 
finite fibers (cf. Lemma 2.9). Second, using Kummer theory, we reconstruct
in a functorial way an embedding $K_2^{\times}\hookrightarrow (K_1^{\perf})^{\times}$ 
between multiplicative groups (cf. Lemma 4.13). Finally, we show that this 
embedding $K_2^{\times}\hookrightarrow (K_1^{\perf})^{\times}$ is additive. 
Recovering the additive structure is one of the main steps 
in the proof. This problem was treated by Uchida in the case of a bijective 
identification $K_2^{\times}\isom K_1^{\times}$ between multiplicative 
groups, which is order-preserving and value-preserving. 
In fact, one needs only to restore the additivity 
between constants. For this one has to show identities of the form
$\gamma (f_2+1)=\gamma (f_2)+1$ for some specific 
non-constant function $f_2\in K_2$. 
Uchida succeeded in his case by choosing $f_2$ to be
a function with a minimal pole divisor (He called such a function a 
minimal element). His argument fails in the case of an embedding 
between multiplicative group which is not surjective, because the image 
of a minimal element is not necessarily minimal in this case.
Roughly speaking, we extend his arguments by 
using, instead, a function which has a unique pole. 
This one pole argument turns out to be very efficient, and leads to the
recovery of the additive structure under quite general assumptions 
(cf. Proposition 5.3).

Although rigid homomorphisms are a special case of proper homomorphisms, 
we choose to treat them 
separately for several reasons. 
First, 
the important condition of inertia-rigidity is 
automatically satisfied in the case of rigid homomorphisms (cf. Remark 4.9(i)). 
Second, 
in the case of (strictly) rigid homomorphisms we can 
reduce directly to the Isom-form, whose proof can be based on class field theory. 
This is not possible for proper homomorphisms, in general. 
In fact, in the case of proper homomorphisms, class field theory reconstructs 
only the norm map between the multiplicative groups of function fields.

This paper is organized as follows. In section 1, we review well-known facts 
concerning Galois theory of function fields of curves over finite fields, 
including the main results of local theory. In section 2, we investigate 
some basic properties of homomorphisms between absolute Galois groups of 
function fields of curves over finite fields, as well as homomorphisms 
between decomposition subgroups. In section 3, we 
investigate rigid homomorphisms between (geometrically prime-to-characteristic 
quotients of) absolute Galois groups, and prove Theorem A. 
In section 4, we investigate proper homomorphisms between 
(geometrically prime-to-characteristic 
quotients of) absolute Galois groups, and prove Theorem B.
In section 5, we investigate the problem of recovering the 
additive structure of function fields. Using the above 
``one pole argument'', we prove Proposition 5.3, which is used in the 
proof of Theorem B in section 4. 

We hope very much that this paper is a first step towards proving the Hom-form
of Grothendieck's anabelian conjecture concerning arithmetic fundamental 
groups of hyperbolic curves over finite fields, whose Isom-form was proven by Tamagawa 
for affine curves ([Tamagawa]) and Mochizuki for proper curves ([Mochizuki4]).

\definition
{Acknowledgment} The first author was holding an EPSRC advanced research fellowship
GR/R75861/02 during the preparation of this paper, and would like very much 
to thank EPSRC for its support. This work was done 
during a visit of his to the Research Institute for Mathematical Sciences (RIMS), 
Kyoto University. He would like very much to thank the staff members of RIMS 
for their hospitality. 
\enddefinition

\subhead
\S 1. Generalities on Galois groups of function fields of curves
\endsubhead

\subsubhead
1.1. Notations on profinite groups and fields
\endsubsubhead

Let $\Cal C$ be a full 
class of finite groups, i.e., $\Cal C$ is closed under taking subgroups, 
quotients, finite products, and extensions. For a profinite group $H$, denote
by $H^{\Cal C}$ the maximal pro-$\Cal C$ quotient of $H$. 
Next, given a profinite group $H$ and its closed normal
subgroup $\overline H$ we set $H^{(\Cal C)}\defeq 
H/\Ker (\overline H\twoheadrightarrow \overline H^{\Cal C})$. 
Note that $H^{(\Cal C)}$ coincides with $H^{\Cal C}$ if
and only if the quotient $A\defeq H/\overline H$ is a pro-$\Cal C$ group. By 
definition, we have the following commutative diagram:
$$
\CD
1@>>>    \overline H    @>>>     H   @>>>    A     @>>> 1\\
  @.        @VVV            @VVV              @V\id VV \\
1     @>>> \overline H^{\Cal C}  @>>>   H^{(\Cal C)}    @>>>A@>>>
1
\endCD
$$
where the rows are exact and the columns are surjective. 

When $\Cal C$
is the class of finite $l$-groups (resp. finite $l'$-groups, i.e., finite groups of order 
prime to $l$), where $l$ is a prime number, write $H^l$ and $H^{(l)}$
(resp. $H^{l'}$ and $H^{(l')}$), 
instead of $H^{\Cal C}$ and $H^{(\Cal C)}$, respectively. 

For a profinite group $H$, we write 
$H^{\ab}$ for the maximal abelian quotient of $H$; 
$\Sub (H)$ for the set of closed subgroups of $H$; 
$\Aut(H)$ for the group of (continuous) automorphisms of $H$; 
and 
$\Inn(H)$ for the group of inner automorphisms of $H$. 

For a profinite group $H$ and a prime number $l$, denote by 
$\cd (H)$ (resp. $\cd _l(H)$) the cohomological (resp. $l$-cohomological) 
dimension of $H$. 
It is well-known that if $\cd(H)<\infty$, then $H$ is torsion-free. 

Let $\kappa$ be a field, 
and 
$\kappa ^{\sep}$ a separable closure of $\kappa$. 
Denote the absolute Galois group $\Gal (\kappa ^{\sep}/\kappa)$ 
by $G_{\kappa}$. 
We shall write
$$M_{\kappa^{\sep}}
\defeq \Hom(\Bbb Q/\Bbb Z,(\kappa^{\sep})^{\times})
.$$
Thus, $M_{\kappa ^{\sep}}$ is a free ${\hat \Bbb Z^{\dag}}$-module of rank 
one, where we write $\hat \Bbb Z^{\dag}\defeq {\hat \Bbb Z}$ 
(resp. $\hat \Bbb Z^{\dag}\defeq {\hat \Bbb Z}^{p'}$ ), if the characteristic 
of $\kappa$ is $0$ (resp. $p>0$). 
Further, $M_{\kappa ^{\sep}}$ has a natural structure of $G_{\kappa}$-module, 
which is isomorphic to the ``Tate twist'' ${\hat \Bbb Z^{\dag}}(1)$, i.e., $G_{\kappa }$ 
acts on $M_{{\kappa}^{\sep}}$ via the cyclotomic character 
$\chi_{\kappa} :G_{\kappa}\to (\hat \Bbb Z^{\dag})^{\times}$.

\subsubhead
1.2. Galois groups of local fields of positive characteristic
\endsubsubhead

Let $p$ be a prime number. Let $L$ be a local field of characteristic $p$, 
i.e., a complete discrete 
valuation field of equal characteristic $p$, with finite residue 
field $\ell$. 
We denote the ring of integers  
of $L$ by $\Cal O_{L}$. 
Also, 
fix
a separable closure $L^{\sep}$ of $L$. 
We shall denote the residue field of $L^{\sep}$ by $\bar \ell$, 
since it is an algebraic closure of $\ell$. Note that 
$\ell$ (resp. $\bar \ell$) can also be regarded naturally 
as a subfield of $L$ (resp. $L^{\sep}$). 
Write $D\defeq \Gal (L^{\sep}/L)$
for the corresponding absolute Galois group of $L$, and 
define the inertia group of $L$ by 
$$I\defeq \{\gamma \in D\mid \text{$\gamma$ acts trivially on $\bar\ell$}\}.$$
We have a canonical exact sequence:
$$1\to I\to D\to G_{\ell}\defeq \Gal (\bar\ell/\ell)\to 1,$$
and, for a full class $\Cal C$ of finite groups, we get a  canonical exact sequence: 
$$1\to I^{\Cal C}\to D^{(\Cal C)}\to G_{\ell}\to 1.$$

The inertia subgroup $I$ possesses a unique $p$-Sylow 
subgroup $I^{\w}$. The quotient
$I^{\t}\defeq I/I^{\w}$ is isomorphic to $\hat \Bbb Z^{p'}$, 
and is naturally identified with the Galois group 
$\Gal (L^{\t}/L^{\ur})$, where $L^{\t}$ (resp. $L^{\ur}$) is
the maximal tamely ramified (resp. maximal unramified) extension of $L$
contained in $L^{\sep}$. We have a natural exact sequence:
$$1\to I^{\t}\to D^{\t}\to G_{\ell}\to 1,$$
where $D^{\t}\defeq \Gal (L^{\t}/L)$. 
(Observe that $I^{\t}=I^{p'}$ and $D^{\t}=D^{(p')}$.)
In particular, 
$I^{\t}$ has a natural structure of $G_{\ell}$-module. Further, there exists a natural 
identification $I^{\t}\isom M_{\bar\ell}$ of $G_{\ell}$-modules.
These follow from well-known facts in ramification theory. See 
[Serre3], Chapitre IV, for more details. 

Let $l$ be a prime number. Denote by $D_{l}$
an $l$-Sylow subgroup of $D$. 
Then the intersection $I_{l}\defeq I\cap D_{l}$ is an $l$-Sylow subgroup of $I$. 
Thus, $I_p=I^{\w}$ and, for $l\neq p$, $I_l$ is isomorphic to $\Bbb Z_l$. 
The image $G_{\ell,l}$ of $D_l$ in $G_{\ell}$ is the unique $l$-Sylow 
subgroup of $G_\ell\simeq \hat\Bbb Z$, hence $G_{\ell,l}\simeq\Bbb Z_l$. 
We have a canonical exact sequence: 
$$1\to I_l\to D_l\to G_{\ell,l}\to 1.$$
In particular, 
$I_l$ has a natural structure of $G_{\ell,l}$-module, and  there exists a natural 
identification $I_l\isom M_{\bar\ell,l}$ of $G_{\ell,l}$-modules, where 
$M_{\bar\ell,l}$ stands for the $l$-Sylow subgroup of the profinite abelian 
group $M_{\bar\ell}$. 

It is well-known that $\cd _l(D)=\cd(D_l)=2$ for any prime number $l\neq p$, and that 
$\cd _p(D)=\cd(D_p)=1$. Thus, $\cd (D)=2<\infty$. In particular, $D$ is torsion-free. 

\proclaim{Proposition 1.1} 
Let $\frak D$ be a quotient of $D$, $\frak I$ the image of $I$ in $\frak D$, 
and $\frak G_{\ell}\defeq \frak D/\frak I$. 
For each prime number $l$, let $\frak D_l$, $\frak I_l$ and 
$\frak G_{\ell,l}$ be the images of 
$D_l$, $I_l$ and $G_{\ell,l}$ 
in $\frak D$, $\frak I$ and $\frak G_{\ell}$, respectively, 
which is an $l$-Sylow subgroup of 
$\frak D$, $\frak I$ and $\frak G_{\ell}$, respectively. 
Let $l$ be a prime number $\neq p$. Then: 

\noindent
{\rm (i)} One of the following {\rm (0)}, {\rm (1)}, {\rm (2)} and {\rm ($\infty$)} occurs.

\noindent
{\rm (0)} $\cd_l(\frak D)=0$, $\frak D_l=\{1\}$, $\frak I_l=\{1\}$, and $\frak G_{\ell,l}=\{1\}.$

\noindent
{\rm (1)} $\cd_l(\frak D)=1$, $\frak D_l\simeq G_{\ell}$, $\frak I_l=\{1\}$, and $\frak G_{\ell,l}\simeq G_{\ell}$. 

\noindent
{\rm (2)} $\cd_l(\frak D)=2$, $\frak D_l\simeq D_l$, $\frak I_l\simeq I_l$, and $\frak G_{\ell,l}\simeq G_{\ell}$. 

\noindent
{\rm ($\infty$)} 
$\cd_l(\frak D)=\infty$, and $\frak I_l$ is a finite group. 

\noindent
{\rm (ii)} Assume that the above case (2) occurs. Let $\frak D'$ be an open subgroup of $\frak D$, 
$L'$ the (finite, separable) extension of $L$ corresponding to $\frak D'\subset \frak D$, and 
$D'$ the inverse image of $\frak D'$ in $D$. (Thus, $D'=G_{L'}$.) 
Then, for each finite $l$-primary $\frak D'$-module $M$ and each $k\geq 0$, 
one has $H^k(\frak D', M)\isom H^k(D',M)$. 
\endproclaim

\demo{Proof} (i) 
Since $\frak I_l$ is a quotient of $I_l\simeq \Bbb Z_l$, one of the following occurs: 
(a) $\frak I_l=\{1\}$, (b) $\frak I_l\simeq \Bbb Z/l^m\Bbb Z$ for an integer $m>0$, and 
(c) $\frak I_l\simeq \Bbb Z_l$. In case (a), $\frak D_l$ is a quotient of 
$D_l/I_l=G_\ell\simeq \Bbb Z_l$. Thus, it is easy to see that 
one of (0), (1), ($\infty$) occurs. 
In case (b), ($\infty$) occurs. In case (c), we have $\frak G_{\ell,l}\simeq G_{\ell,l}$. 
This follows from the fact that $I_l$ is isomorphic to $M_{\bar\ell,l}$ on which 
$G_{\ell,l}$ acts via the $l$-adic cyclotomic character, 
and that the $l$-adic cyclotomic charter $\chi_l: G_{\ell,l}\to\Bbb Z_l^{\times}$ 
is injective. Thus, it is easy to see that (2) occurs in this case. 

\noindent
(ii) Replacing $L$ by $L'$, we may assume that $L'=L$. (Observe that case (2) occurs also 
for the quotient $G_{L'}=D'\twoheadrightarrow \frak D'$.) 

Denote by $N$ the kernel of the surjection $D\twoheadrightarrow \frak D$. 
By the assumption that case (2) occurs, $D_l$ is injectively mapped into 
$\frak D$, hence $D_l\cap N$, which is an $l$-Sylow subgroup of $N$, 
is trivial. Namely, $N$ is of order prime to $l$, hence we have 
$H^k(D,M)=H^k(\frak D,H^0(N, M))=H^k(\frak D, M)$, as desired. 
\qed\enddemo


\subsubhead
1.3. Galois groups of function fields of curves
\endsubsubhead

Let $k$ be a finite field of characteristic $p>0$. Let $X$ be a proper, 
smooth, geometrically connected curve over $k$. 
Let $K=K_X$ be the function field of $X$ 
and fix an algebraic closure $\overline K$ of $K$. 
Write $K^{\sep}$ (resp.
$\bar k=k^{\sep}$) for the separable closure of $K$ (resp. $k$) in $\overline K$. 
Write $G=G_K\defeq \Gal (K^{\sep}/K)$ and $G_k\defeq \Gal (\bar k/k)$ 
for the absolute Galois groups of $K$ and $k$, respectively. 
We have the following exact sequence of profinite groups: 
$$1\to \overline G\to G@>\text {pr}>> G_k\to 1,\tag {1.1}$$
where $\overline G$ is the absolute Galois group $G_{K\bar k}=\Gal (K^{\sep}/K\bar k)$ 
of $K\bar k$, and $\text {pr}$ is the canonical projection.
Here, it is well-known that the right term $G_k$ is 
a profinite free group of rank $1$ which is (topologically) generated by 
the Frobenius element, while 
the left term $\overline G$ is a profinite free group of countably infinite rank 
(cf. [Pop2], [Harbater]). However, the structure of the extension $(1.1)$ itself 
is not understood well. From $(1.1)$ above, we also 
obtain the exact sequence:
$$1\to \overline G^{\Cal C}\to G^{(\Cal C)}@>\pr>> G_k\to 1
$$
for each full class $\Cal C$ of finite groups. 

In the rest of this section, let $N$ be a closed normal subgroup of $G$ and set 
$\frak G\defeq G/N$. Let $\tilde K$ denote the Galois extension of $K$ 
corresponding to $N$, i.e., $\tilde K\defeq (K^{\sep})^{N}$. 
Let $\bar \frak G$ be the image of $\bar G$ in $\frak G$, and set 
$\frak G_k\defeq \frak G/\bar \frak G$, which is a quotient of 
$G_k$. 

For a scheme $T$ denote by $\Sigma _T$ the set of closed points of $T$. 
Write $\tilde{\tilde X}$ 
for the integral closure of $X$ in $K^{\sep}$. 
The absolute Galois group $G$ acts naturally on the set 
$\Sigma _{\tilde{\tilde X}}$, and the quotient  $\Sigma _{\tilde{\tilde X}}/G$ is naturally 
identified with $\Sigma _X$. 
For a point $\tilde{\tilde x}\in \Sigma _{\tilde{\tilde X}}$, with residue field 
$k(\tilde{\tilde x})$ (which is naturally identified with $\bar k$), 
we define its decomposition group $D_{\tilde{\tilde x}}$ and inertia group  $I_{\tilde{\tilde x}}$ 
by 
$$D_{\tilde{\tilde x}}\defeq \{\gamma \in G\mid \gamma (\tilde{\tilde x})=\tilde{\tilde x}\}$$
and
$$I_{\tilde{\tilde x}}\defeq \{\gamma \in D_{\tilde{\tilde x}}\mid 
\text{$\gamma$ acts trivially on $k (\tilde{\tilde x})$}\},$$
respectively. We have a canonical exact sequence: 
$$1\to I_{\tilde{\tilde x}}\to D_{\tilde{\tilde x}}\to G_{k(x)}\to 1.$$

More generally, write $\tilde X$ for the integral closure of $X$ in $\tilde K$. 
The Galois group $\frak G$ acts naturally on the set 
$\Sigma _{\tilde X}$, and the quotient  $\Sigma _{\tilde X}/\frak G$ is naturally 
identified with $\Sigma _X$. 
For a point $\tilde x\in \Sigma _{\tilde X}$, with residue field 
$k(\tilde x)$ (which is naturally identified with a subfield of $\bar k$), 
we define its decomposition group $\frak D_{\tilde x}$ and inertia group  $\frak I_{\tilde x}$ by 
$$\frak D_{\tilde x}\defeq \{\gamma \in \frak G\  |\  \gamma (\tilde x)=\tilde x\}$$
and
$$\frak I_{\tilde x}\defeq \{\gamma \in \frak D_{\tilde x}\  |\  \gamma \ \text {acts trivially on}\  
k (\tilde x)\},$$
respectively. 
(Observe that for any $g\in \frak G$, one has 
$\frak D_{g\tilde x}=g\frak D_{\tilde x}g^{-1}$ 
and $\frak I_{g\tilde x}=g\frak I_{\tilde x}g^{-1}$.) 
Set $\frak G_{k(x)}\defeq\frak D_{\tilde x}/\frak I_{\tilde x}$. 
Thus, if we take a point $\tilde{\tilde x}\in\Sigma_{\tilde{\tilde X}}$ above 
$\tilde x\in\Sigma_{\tilde X}$, 
$\frak D_{\tilde x}$, $\frak I_{\tilde x}$ and $\frak G_{k(x)}$ 
are quotients of $D_{\tilde{\tilde x}}$, $I_{\tilde{\tilde x}}$ and $G_{k(x)}$, respectively. 
We have a canonical exact sequence:
$$1\to \frak I_{\tilde x}\to \frak D_{\tilde x}\to \frak G_{k(x)}\to 1.$$
For each closed subgroup $\frak H\subset\frak G$, 
denote by $\tilde x_{\frak H}$ the image of $\tilde x$ in $X_{\frak H}$. Define 
$\tilde K_{\tilde x}\defeq \underset {\frak H\subset \frak G} \to \bigcup (K_{\frak H})_
{\tilde x_{\frak H}}$, 
where $\frak H$ runs over all open subgroups of $\frak G$, and  
$(K_{\frak H})_{\tilde x_{\frak H}}$ means the ${\tilde x_{\frak H}}$-adic completion 
of $K_{\frak H}\defeq (\tilde K)^{\frak H}$. 
Then 
the Galois group $\Gal (\tilde K_
{\tilde x}/K_x)$ is naturally identified with $\frak D_{\tilde x}$, 
where $x\defeq 
\tilde x_{\frak G} \in\Sigma_X$. 

In the rest of this subsection, we fix a prime number $l\neq p$, and put 
the following two assumptions. First, $N^l=N$, or, equivalently, $\tilde K$ admits 
no $l$-cyclic extension. Second, $\tilde K$ contains a primitive $l$-th roots of unity. 

\definition{Remark 1.2} 
Let $\Cal C$ be a full class of finite group. 

\noindent
(i) If $\Bbb F_l \in\Cal C$, then the quotient $G^{(\Cal C)}$ of $G$ satisfies 
the above two assumptions. 

\noindent
(ii) If $\Bbb F_l \in\Cal C$ and $\Gal(K(\zeta_l)/K)\in\Cal C$, 
then the quotient $G^{\Cal C}$ of $G$ satisfies 
the above two assumptions. 
\enddefinition

\proclaim{Lemma 1.3} Let $\tilde x\in\Sigma_{\tilde X}$ and take $\tilde{\tilde x}\in
\Sigma_{\tilde{\tilde X}}$ above $\tilde x$. Let $D_{\tilde{\tilde x},l}$ be an $l$-Sylow 
subgroup of $D_{\tilde{\tilde x}}$ and $\frak D_{\tilde x,l}$ the image of $D_{\tilde{\tilde x},l}$ 
under the natural surjection 
$D_{\tilde{\tilde x}}\twoheadrightarrow \frak D_{\tilde x}$, 
which is an $l$-Sylow subgroup of $\frak D_{\tilde x}$. Then the natural surjection 
$D_{\tilde{\tilde x},l}\twoheadrightarrow \frak D_{\tilde x,l}$ is an isomorphism. 
\endproclaim

\demo{Proof}
Take $t\in K$ such that $t$ is a uniformizer at $x\defeq\tilde x_{\frak G}\in\Sigma_X$. 
Then by the two assumptions (and by Kummer theory), any $l^n$-th root $t^{1/l^n}$ 
of $t$ is contained in $\tilde K$. From this, it follows that 
$I_{\tilde{\tilde x},l}\defeq D_{\tilde{\tilde x},l}\cap I_{\tilde{\tilde x}}$ 
is injectively mapped into $\frak D_{\tilde x}$. Now, applying Proposition 1.1(i)  
to the quotient $G_{K_x}=D_{\tilde{\tilde x}}\twoheadrightarrow \frak D_{\tilde x}$, 
we conclude that case (2) of loc. cit. can only occur, as desired. 
\qed\enddemo

\proclaim{Lemma 1.4} 
Let $\frak G'$ be an open subgroup of $\frak G$, 
$K'$ the (finite, separable) extension of $K$ corresponding to $\frak G'\subset \frak G$, and 
$G'$ the inverse image of $\frak G'$ in $G$. (Thus, $G'=G_{K'}$.) 
Then, for each finite $l$-primary $\frak G'$-module $M$ and each $k\geq 0$, 
one has $H^k(\frak G', M)\isom H^k(G',M)$. 
\endproclaim

\demo{Proof} 
Replacing $K$ by $K'$, we may assume that $K'=K$. (Observe that the two assumptions 
also hold for the quotient $G'\defeq G_{K'}\twoheadrightarrow \frak G'=G'/N$.) 
By Lemma 1.3, one has $\cd_l(N)\leq 1$. (See [Serre1], Chapitre II, Proposition 9, which 
only treats the number field case but whose proof works as it is in our function field 
case.) Next, by the assumption that $N^l=N$, one has 
$H^1(N,M)=\Hom(N,M)=0$. 
Thus, we have 
$H^k(G,M)=H^k(\frak G,H^0(N, M))=H^k(\frak G, M)$, as desired. 
\qed\enddemo

\proclaim
{Proposition 1.5} (Galois Characterization of Decomposition Subgroups) 
{\rm (i)} Let $\tilde x\neq \tilde x'$ be two elements of $\Sigma _{\tilde X}$. Then 
$\frak D_{\tilde x}\cap \frak D_{\tilde x'}$ is of order prime to $l$, hence, in particular, 
is open neither in $\frak D_{\tilde x}$ nor in $\frak D_{\tilde x'}$. 

\noindent
{\rm (ii)} Let $\Dec_l(\frak G)\subset\Sub(\frak G)$ be the set of closed subgroups $\frak D$ 
of $\frak G$ satisfying the following property: 
There exists an open subgroup $\frak D_0$ of $\frak D$ such that for 
any open subgroup $\frak D'\subset\frak D_0$, 
$\dim_{\Bbb F_l}H^2(\frak D',\Bbb F_l)=1$. Define $\Dec^{\max}_l(\frak G)\subset\Dec_l(\frak G)$ 
to be the set of maximal elements of $\Dec_l(\frak G)$. 
Then, 
the map 
$\Sigma _{\tilde X}\to \Sub (\frak G)$, $\tilde x \mapsto \frak D_{\tilde x}$, 
induces a bijection $\Sigma_{\tilde X}\isom \Dec^{\max}_l(\frak G)$, and, in particular, 
it is injective. 
\endproclaim

\demo{Proof}
(i) As in [Uchida2], this follows from the approximation theorem 
(cf. [Neukirch2], Lemma 8). More precisely, let 
$\frak D_l$ be an $l$-Sylow subgroup of $\frak D_{\tilde x}\cap \frak D_{\tilde x'}$, 
and suppose that $\frak D_l\neq 1$. Since $\frak D_l\subset \frak D_{\tilde x,l}$ 
is torsion-free, $\frak D_l$ is an infinite group. Thus, one may replace 
$\frak G$ by any open subgroup, and assume that $\zeta_l\in K$, that 
the images $x$ and $x'$ 
in $\Sigma_{X}$ of $\tilde x$ and $\tilde x'$ are distinct, and that 
the image of $\frak D_l$ in $\frak G^{\ab}/(\frak G^{\ab})^l$ is nontrivial. 
In particular, this implies that the natural map 
$$\frak D_{\tilde x}^{\ab}/(\frak D_{\tilde x}^{\ab})^l \times 
\frak D_{\tilde x'}^{\ab}/(\frak D_{\tilde x'}^{\ab})^l \to 
\frak G^{\ab}/(\frak G^{\ab})^l$$ 
is not injective. By Kummer theory, this last condition is equivalent to saying that 
the natural map 
$$K^{\times}/(K^{\times})^l\to K_x^{\times}/(K_x^{\times})^l\times K_{x'}^{\times}/(K_{x'}^{\times})^l$$ 
is not surjective. This contradicts the approximation theorem. (Note that 
$(K_x^{\times})^l$ and $(K_{x'}^{\times})^l$ are open in 
$K_x^{\times}$ and $K_{x'}^{\times}$, respectively.)

\noindent
(ii) By means of Proposition 1.1(i), Lemma 1.3 and Lemma 1.4, 
the proof of [Uchida2] (which is essentially due to Neukirch, cf. 
[Neukirch1] and [Neukirch2]) 
works as it is. See Lemmas 1-3 of loc. cit. for more details. 
\qed\enddemo

\definition{Remark 1.6}
For other characterizations of decomposition groups (which 
are applicable to much more general situations), see 
[Pop1], Theorem 1.16, [Koenigsmann], Theorem 2, [EK], [EN],.... 
\enddefinition


%
%
%

\subsubhead
1.4. Fundamental groups of curves
\endsubsubhead

Write $I\defeq \langle I_{\tilde {\tilde x}}\rangle_{\tilde {\tilde x}\in \Sigma _{\tilde {\tilde X}}}$ 
for the closed subgroup of $G$ generated by the inertia subgroups 
$I_{\tilde{\tilde x}}$ for all $\tilde{\tilde x}\in \Sigma _{\tilde{\tilde X}}$, and call it the
inertia subgroup of $G$. Then $I$ is normal in $G$. The quotient 
$G/I$ is canonically identified with the fundamental group 
$\pi_1(X)
$ 
of $X$ with base point 
$\Spec(\overline K)\to X$ 
(cf. [SGA-1]). 
We have a natural exact sequence:
$$1\to \pi_1(\overline X)\to \pi_1(X) @>\pr>> G_k\to 1,$$ 
where $\pi_1(\overline X)$ is the fundamental group of 
$\overline X\defeq X\times _k\bar k$ with base point 
$\Spec(\overline K)\to \overline X$ 
and $\pr$ is 
the canonical projection. 
We have the following exact sequence:
$$1\to \pi_1(X)^{\ab,\tor}\to \pi_1(X)^{\ab} @>\pr>> G_k\to 1,$$ 
where $\pi_1(X)^{\ab,\tor}$ is the torsion subgroup of $\pi_1(X)^{\ab}$,
and $\pr$ is the canonical projection. Moreover, $\pi_1(X)^{\ab,\tor}$ is a 
finite abelian group which is canonically isomorphic to the group 
$J_X(k)$ of $k$-rational points of the Jacobian variety $J_X$ of $X$. 

More generally, write $\frak I\defeq \langle\frak I_{\tilde x}\rangle_{\tilde x\in \Sigma _{\tilde X}}$ 
for the closed subgroup of $\frak G$ generated by the inertia subgroups 
$\frak I_{\tilde x}$ for all $\tilde x\in \Sigma _{\tilde X}$, and call it the
inertia subgroup of $\frak G$. Then $\frak I$ is normal in $\frak G$. 
Set $\Pi_X\defeq\frak G/\frak I$, which is a quotient of 
$\pi_1(X)$. Define $\Pi_{\overline X}$ to be the image of $\pi_1(\overline X)$ 
in $\Pi_X$ Then we have a natural exact sequence:
$$1\to \Pi_{\overline X}\to \Pi_{X} @>\pr>> \frak G_k\to 1.$$ 

When $\frak G=G^{(\Cal C)}$ for a full class $\Cal C$ of finite groups, 
we have $\Pi_X=\pi_1(X)^{(\Cal C)}$. In this case, 
we have the following exact sequence:
$$1\to \Pi_X^{\ab,\tor}\to \Pi_X^{\ab} @>\pr>> G_k\to 1,$$ 
where $\Pi_X^{\ab,\tor}$ is the torsion subgroup of $\Pi_X^{\ab}$. 
Moreover, $\Pi_X^{\ab,\tor}$ is a 
finite abelian group which is canonically isomorphic to 
the maximal (pro-)$\Cal C$ quotient $J_X(k)^{\Cal C}$ of 
the finite group $J_X(k)$.

\subhead 
\S 2. Basic properties of homomorphisms between Galois groups
\endsubhead

In this section we investigate some basic properties of homomorphisms 
between Galois groups of function fields of curves over finite fields.
First, we shall investigate a class of homomorphisms between 
decomposition subgroups, which arise naturally from the class of 
homomorphisms between (quotients of) Galois groups that we 
will consider in $\S 3$ and $\S 4$.

\subsubhead
2.1. Homomorphisms between Galois groups of local fields of positive characteristics
\endsubsubhead

For $i\in \{1,2\}$, let $p_i>0$ be a prime number. Let $L_i$ be a complete discrete 
valuation field of equal characteristic $p_i$, with finite residue 
field $\ell _i$. We denote the ring of integers  
of $L_i$ by $\Cal O_{L_i}$. Also, 
fix 
a separable closure $L_i^{\sep}$ of $L_i$. 
We shall denote the residue field of $L_i^{\sep}$ by $\bar \ell_i$, 
since it is an algebraic closure of $\ell_i$. Note that 
$\ell_i$ (resp. $\bar \ell_i$) can also be regarded naturally 
as a subfield of 
$L_i$ (resp. $L_i^{\sep}$). 
Write $D_i\defeq \Gal (L_i^{\sep}/L_i)$
for the corresponding absolute Galois group of $L_i$, 
and $I_i\subset D_i$ the inertia subgroup. 
For each prime number $l$, let $D_{i,l}$ be an $l$-Sylow subgroup of $D_{i}$. 

By local class field theory (cf., e.g., [Serre2]), we have a natural 
isomorphism $(L_i^\times)^{\wedge}\isom D_i^{\ab}$,
where $(L_i^\times)^{\wedge}\defeq 
\underset{n}\to{\varprojlim}\  L_i^\times/(L_i^\times)^n$.  
In particular, 
$D_i^{\ab}$ fits into an exact sequence
$$0\to \Cal O_{L_i}^{\times}\to D_i^{\ab}\to \hat \Bbb Z \to 0$$
(arising from a similar exact sequence for $(L_i^{\times})^{\wedge}$),
where $\Cal O_{L_i}^{\times}$ is the group of multiplicative units in
$\Cal O_{L_i}$. Moreover, we obtain natural inclusions
$$\ell _i^{\times}\times  U_i^1 =\Cal O_{L_i}^{\times}\subset L_i^{\times}
\hookrightarrow D_i^{\ab},$$
where 
$U_i^1$ is the group of principal units in $\Cal O_{L_i}^{\times}$, 
and
$$L_i^{\times}/\Cal O_{L_i}^{\times}\isom \Bbb Z \hookrightarrow 
D_i^{\ab}/\Im (\Cal O_{L_i}^{\times})$$
(where $\isom$ is the isomorphism induced by the valuation),
by considering the Frobenius element. 


Let $\frak D_i$ be a quotient of $D_i$, $\frak I_i$ the image of $I_i$ in $\frak D_i$, 
and $\frak G_{\ell_i}\defeq \frak D_i/\frak I_i$. 
For each prime number $l$, let $\frak D_{i,l}$ 
be the 
image 
of $D_{i,l}$ 
in $\frak D_i$,  
which is an $l$-Sylow subgroup of 
$\frak D_i$. 
Write 
$$\frakIm(\ell _i^{\times}),\  \frakIm(U_i^1) \subset \frakIm(\Cal O_{L_i}^{\times})\subset \frakIm(L_i^{\times}) 
\subset\frak D_i^{\ab}$$
for the images of 
$\ell _i^{\times}$, $U_i^1$, $\Cal O_{L_i}^{\times}$ and $L_i^{\times}$ in $\frak D_i^{\ab}$, respectively. 
In the rest of this subsection, we 
assume that either $\frak D_i=D_i$, $i=1,2$ or 
$\frak D_i=D_i^{\t}=D_i^{(p_i')}$, $i=1,2$, and refer to the 
first and the second cases as 
the profinite and the tame cases, respectively. 
Thus, we have 
$\frak D_i^{\ab}=(L_i^{\times})^{\wedge}$, 
$\frakIm(L_i^{\times})=L_i^{\times}$, 
$\frakIm(\Cal O_{L_i}^{\times})=\Cal O_{L_i}^{\times}$, 
$\frakIm(\ell_i^{\times})=\ell_i^{\times}$ 
and 
$\frakIm(U_i^1)=U_i^1$ 
in the profinite case, and 
$\frak D_i^{\ab}=(L_i^{\times})^{\wedge}/U_i^1$, 
$\frakIm(L_i^{\times})=L_i^{\times}/U_i^1$, 
$\frakIm(\Cal O_{L_i}^{\times})=\Cal O_{L_i}^{\times}/U_i^1
=\frakIm(\ell_i^{\times})=\ell_i^{\times}$ 
and 
$\frakIm(U_i^1)=\{1\}$ 
in the tame case. 

Let 
$$\tau : \frak D_1\twoheadrightarrow \frak D_2$$ 
be a 
surjective homomorphism between 
profinite groups. 
Write $\tau ^{\ab}:\frak D_1^{\ab}\twoheadrightarrow \frak D_2^{\ab}$ 
for the induced 
surjective homomorphism between the 
maximal abelian quotients. 
For each prime number $l$, $\tau(\frak D_{1,l})$ is an $l$-Sylow 
subgroup of $\frak D_2$, and we shall assume that 
$\tau(\frak D_{1,l})=\frak D_{2,l}$. 

\proclaim
{Proposition 2.1} (Invariants of Arbitrary Surjective Homomorphisms between
Decomposition Groups) {\rm (i)} The equality $p_1=p_2$ holds. Set $p\defeq p_1=p_2$.

\noindent
{\rm (ii)}  Let $l\neq p$ be a prime number. 
We have $\frak D_{1,l}\cap \Ker \tau=\{1\}$. 
In particular, $\Ker\tau$ is pro-$p$. 
In the tame case, $\tau$ is an isomorphism. 

\noindent
{\rm (iii)}  The homomorphism $\tau $ induces a natural bijection 
$\ell _1^{\times} \isom \ell _2^{\times}$ between the multiplicative 
groups of residue fields. In particular, $\ell _1$ and $\ell _2$ have 
the same cardinality. 

\noindent
{\rm (iv)}  $\tau $ induces naturally an isomorphism $M_{\bar \ell _1}\isom
M_{\bar \ell _2}$, which is Galois-equivariant with respect to $\tau$. 
In particular, $\tau$ commutes with the cyclotomic 
characters $\chi _i:\frak D_i\to ({\hat \Bbb Z}^{p'})^{\times}$ of $\frak D_i$, 
i.e., the following diagram is commutative:

$$
\CD
({\hat \Bbb Z}^{p'})^{\times}  @= ({\hat \Bbb Z}^{p'})^{\times}\\
@A{\chi _1}AA    @A{\chi _2}AA   \\
\frak D_1  @>\tau>> \frak D_2 \\
\endCD
$$

\noindent
{\rm (v)} We have $\tau (\frak I_1)=\frak I_2$. 

\noindent
{\rm (vi)} The homomorphism $\tau ^{\ab}:D_1^{\ab}\to D_2^{\ab}$ 
preserves 
$\frakIm (L_i^{\times})$, 
$\frakIm (\Cal O_{L_i}^{\times})$, 
$\frakIm (\ell_i^{\times})$ 
and 
$\frakIm (U_i^1)$. 
Further, the 
isomorphism 
$\frak D_1^{\ab}/\frakIm (\Cal O_{L_1}^{\times})
\to D_2^{\ab}/\frakIm (\Cal O_{L_2}^{\times})$ induced by 
$\tau$ preserves the respective Frobenius elements.
\endproclaim

\demo
{Proof} 
Property (i) follows by considering the 
$q$-Sylow subgroups of 
$\frak D_i$ for various prime numbers $q$. 
Indeed, 
for $i\in \{1,2\}$, 
$\frak D_{i,p_i}$ is 
not (topologically) finitely generated 
(resp. is cyclic) 
in the profinite (resp. tame) case, 
while $\frak D_{i,l}$ 
for a prime number $l\neq p_i$ is (topologically) finitely generated 
and non-cyclic. 
Accordingly, the surjection $\frak D_{1,p_2}\twoheadrightarrow \frak D_{2,p_2}$ 
(resp. $\frak D_{1,p_1}\twoheadrightarrow \frak D_{2,p_1}$) cannot exist 
in the profinite (resp. tame) case, unless $p_1=p_2$. 
Thus, we must have $p_1=p_2$. 

The first assertion of (ii) follows from Proposition 1.1(i), applied to the quotient 
$D_1\twoheadrightarrow\frak D_1\overset{\tau}\to{\twoheadrightarrow} \frak D_2$. 
The second assertion follows from the first. The third assertion follows from 
the second, together with the fact (which can be checked easily) that $D_1^{\t}$ 
admits no nontrivial normal pro-$p$ subgroup. 

Next, we prove (iii). 
By local class field theory, the torsion subgroup 
$\frak D_{i}^{\ab,\tor}$ of $\frak D_{i}^{\ab}$ 
is naturally
identified with ${\ell _i}^{\times}$ 
(both in the profinite and the tame cases), hence, in particular, 
is finite of order prime to $p$. 
By (ii), the kernel of the surjective homomorphism 
$\tau^{\ab}: \frak D_{1}^{\ab}\to \frak D_{2}^{\ab}$ 
is pro-$p$. 
Thus, $\tau^{\ab}$ induces a natural isomorphism 
$\frak D_{1}^{\ab,\tor}\isom\frak D_{2}^{\ab,\tor}$, 
which is naturally identified with 
$\ell _1^{\times}\isom \ell _2^{\times}$, as desired. 

By applying the above argument to open subgroups of $\frak D_i$ (which 
correspond to each other via $\tau$); $i=1,2$, and 
passing to the projective limit with respect to the norm maps, 
we obtain a natural isomorphism
$M_{\bar \ell _1}\isom M_{\bar \ell _2}$ between the modules of 
roots of unity. Here, we use the fact that if $L_i'$ is a finite
extension of $L_i$ corresponding to the open subgroup $\frak D_i'$ of $\frak D_i$,
then the following diagram commutes:
$$
\CD
({L_i'}^{\times})^{\wedge} @>>>  {\frak D_i'}^{\ab} \\
@V{\Norm}VV     @VVV  \\
({L_i}^{\times})^{\wedge} @>>>  \frak D_i^{\ab}, 
\endCD
$$
where the horizontal maps are the natural isomorphisms from local
class field theory, and the map ${\frak D_i'}^{\ab}\to \frak D_i^{\ab}$ is induced
by the natural inclusion ${\frak D_i'}\subset {\frak D_i}$. Further, this 
identification is (by construction) Galois-compatible with respect to the 
homomorphism $\tau$. This completes the proof of (iv). 

Property (v) follows from property (iv), since $\frak I_i$ coincides with 
the kernel of $\chi_i$ for $i=1,2$. 

Next, we prove (vi). 
First, $\tau ^{\ab}$ preserves the image 
$\frakIm (\Cal O_{L_i}^{\times})$ by (v), since this image coincides with 
the image of the inertia subgroup $\frak I_i$. Since 
$\frakIm(\ell_i^{\times})$ (resp. $\frakIm(U_i^1)$) is
the maximal prime-to-$p$ (resp. pro-$p$) subgroup of 
$\frakIm (\Cal O_{L_i}^{\times})$, property (vi) for 
$\frakIm(\ell_i^{\times})$ (resp. $\frakIm(U_i^1)$) follows. 
Further, by (iii) and (iv), the homomorphism $\frak D_1^{\ab}/\frakIm (\Cal O_{L_1}^{\times})
\to \frak D_2^{\ab}/\frakIm (\Cal O_{L_2}^{\times})$ induced by $\tau$ 
preserves the respective Frobenius elements, since such an element 
is characterized as the unique element 
whose image under $\chi_i$ is $\sharp(\ell_i)$. 
Finally, since $\frakIm(L_i^\times)$ is the inverse image in $\frak D_i^{\ab}$ of 
the subgroup generated by the Frobenius element in $\frak D_i^{\ab}/\frakIm (\Cal O_{L_i}^{\times})$ for 
$i=1,2$, they are preserved by $\tau^{\ab}$. 
\qed
\enddemo

\subsubhead
2.2. Homomorphisms between Galois groups of function fields of curves over finite fields
\endsubsubhead

Next, we shall investigate some basic properties of 
homomorphisms between Galois groups of function 
fields of curves over finite fields. We follow the notations in $\S 1$, especially 
subsections 1.3 and 1.4. We also follow the following: 

\definition{Notation} 
(i) For $i\in \{1,2\}$, let $k_i$ be a finite field of characteristic 
$p_i> 0$. Let $X_i$ be a smooth, proper, geometrically connected 
curve of genus $g_i\geq 0$ over $k_i$. Let $K_i=K_{X_i}$ be the function field of $X_i$ and 
fix an algebraic closure $\overline{K}_i$ of $K_i$. Let 
$K_i^{\sep}$ be the separable closure of $K_i$ in $\overline{K}_i$,  
and $\bar k_i$ the algebraic closure of $k_i$ in $\overline{K}_i$. 
Following the notations 
in $\S 1$, we will write $G_i\defeq G_{K_i}=\Gal (K_i^{\sep}/K_i)$
for the absolute Galois group of $K_i$, and $\bar{G}_i\defeq 
G_{K_i\bar k_i}=\Gal (K_i^{\sep}/K_i\bar k_i)$
for the absolute Galois group of $K_i\bar k_i$. 

\noindent
(ii) Let $N_i$ be a normal closed subgroup of $G_i$ and set 
$\frak G_i\defeq G_i/N_i$. Let $\tilde K_i$ denote the Galois extension of $K_i$ 
corresponding to $N_i$, i.e., $\tilde K_i\defeq (K_i^{\sep})^{N_i}$. 
Let $\bar \frak G_i$ be the image of $\bar G_i$ in $\frak G_i$, and set 
$\frak G_{k_i}\defeq \frak G_i/\bar \frak G_i$, which is a quotient of 
$G_{k_i}=\Gal(\bar k_i/k_i)$. For $i=1,2$, let us denote by $\varphi_{k_i}$ the image 
in $\frak G_{k_i}$ of the $\sharp(k_i)$-th power Frobenius element of $G_{k_i}$. 

\noindent
(iii) Write $\tilde X_i$ for the integral closure of $X_i$ in $\tilde K_i$. 
The Galois group $\frak G_i$ acts naturally on the set 
$\Sigma _{\tilde X_i}$, and the quotient  $\Sigma _{\tilde X_i}/\frak G_i$ 
is naturally identified with $\Sigma _{X_i}$. Denote the natural quotient 
map $\Sigma_{\tilde X_i}\to\Sigma_{X_i}$ by $q_i$. 
For a point $\tilde x_i\in \Sigma _{\tilde X_i}$, with residue field 
$k_i(\tilde x_i)$ (which is naturally identified with a subfield of $\bar k_i$), 
we define its decomposition group $\frak D_{\tilde x_i}$ and inertia group  $\frak I_{\tilde x_i}$ by 
$$\frak D_{\tilde x_i}\defeq \{\gamma \in \frak G_i\  |\  \gamma (\tilde x_i)=\tilde x_i\}$$
and
$$\frak I_{\tilde x_i}\defeq \{\gamma \in \frak D_{\tilde x_i}\  |\  \gamma \ \text {acts trivially on}\  
k_i (\tilde x_i)\},$$
respectively. 
Set $\frak G_{k_i(x_i)}\defeq\frak D_{\tilde x_i}/\frak I_{\tilde x_i}$. 

Write $\frak I_i\defeq \langle I_{\tilde x_i}\rangle_{\tilde x_i\in \Sigma _{\tilde X_i}}$ 
for the closed subgroup of $\frak G_i$ generated by the inertia subgroups 
$\frak I_{\tilde x_i}$ for all $\tilde x_i\in \Sigma _{\tilde X_i}$, and call it the
inertia subgroup of $\frak G_i$. Then $\frak I_i$ is normal in $\frak G_i$. 

\noindent
(iv) Let 
$$\sigma:\frak G_1\to \frak G_2$$ 
be a continuous 
homomorphism between profinite groups. 
\enddefinition

\proclaim
{Proposition 2.2} (Image of a Decomposition Subgroup) 
Let $l$ be a prime number $\neq p_1,p_2$, and put 
the following two assumptions: 
{\rm (1)} $N_2^l=N_2$, or, equivalently, $\tilde K_2$ admits 
no $l$-cyclic extension; and 
{\rm (2)} $\tilde K_2$ contains a primitive $l$-th roots of unity. 
For each $\tilde x_1\in \Sigma _{\tilde X_1}$, fix 
an $l$-Sylow subgroup $\frak D_{\tilde x_1,l}$ of $\frak D_{\tilde x_1}$ and 
set $\frak I_{\tilde x_1,l}\defeq \frak I_{\tilde x_1}\cap \frak D_{\tilde x_1,l}$, 
which is an $l$-Sylow subgroup of $\frak I_{\tilde x_1}$. 
Let $\Sigma_{\tilde X_1, \sigma, l}$ be the set of 
$\tilde x_1\in\Sigma_{\tilde X_1}$ such that $\cd_l(\sigma(\frak D_{\tilde x_1}))=2$. 
Then: 

\noindent
{\rm (i)} There exists a unique map 
$\tilde\phi=\tilde\phi_{\sigma,l}: \Sigma_{\tilde X_1, \sigma, l}\to\Sigma_{\tilde X_2}$, 
such that $\sigma(\frak D_{\tilde x_1})\subset\frak D_{\tilde\phi(\tilde x_1)}$ 
for each $\tilde x_1\in \Sigma_{\tilde X_1, \sigma, l}$. 

\noindent
{\rm (ii)} For each $\tilde x_1\in\Sigma_{\tilde X_1, \sigma, l}$, 
there exists 
an 
$l$-Sylow subgroup $\frak D_{\tilde\phi(\tilde x_1),l}$
of $\frak D_{\tilde\phi(\tilde x_1)}$, such that 
$\sigma(\frak D_{\tilde x_1,l}) \subset \frak D_{\tilde\phi(\tilde x_1),l}$. 
Moreover, we have 
$\sigma(\frak I_{\tilde x_1,l}) \subset \frak I_{\tilde\phi(\tilde x_1),l}$, 
where we set 
$\frak I_{\tilde\phi(x_1),l}\defeq \frak I_{\tilde\phi(\tilde x_1)}\cap \frak D_{\tilde\phi(\tilde x_1),l}$, 
which is an $l$-Sylow subgroup of $\frak I_{\tilde\phi(\tilde x_1)}$. 

\noindent
{\rm (iii)} The subset $\Sigma_{\tilde X_1,\sigma, l}\subset \Sigma_{\tilde X_1}$ 
is $\frak G_1$-stable, or, equivalently, $\Sigma_{\tilde X_1,\sigma, l}=
q_1^{-1}(\Sigma_{X_1,\sigma,l})$, where $\Sigma_{X_1,\sigma,l}\defeq 
q_1(\Sigma_{\tilde X_1,\sigma, l})$. 
The map $\tilde\phi$ is Galois-compatible with respect to $\sigma$: 
$\tilde\phi(g_1\tilde x_1)=\sigma(g_1)\tilde\phi(\tilde x_1)$ for any 
$\tilde x_1\in \Sigma_{\tilde X_1,\sigma, l}$ and any 
$g_1\in\frak G_1$. In particular, $\tilde\phi$ induces naturally 
a map $\phi=\phi_{\sigma,l}:\Sigma_{X_1,\sigma,l}\to\Sigma_{X_2}$. 

\noindent
{\rm (iv)} For any $\tilde x_1\in
\Sigma_{\tilde X_1} \smallsetminus \Sigma_{\tilde X_1,\sigma, l}$, 
we have $\sigma(\frak I_{\tilde x_1,l})=\{1\}$. 

\noindent
{\rm (v)} For two primes $l=l_1,l_2$ satisfying the assumptions, 
$\tilde\phi_{\sigma,l_1}$ and $\tilde\phi_{\sigma,l_2}$ coincide with each other 
on the intersection $\Sigma_{\tilde X_1,\sigma, l_1}\cap\Sigma_{\tilde X_1,\sigma, l_2}$
\endproclaim

\demo{Proof} (i) Take $\tilde x_1\in \Sigma_{\tilde X_1, \sigma, l}$. 
Applying Proposition 1.1(i)(ii) to $\frak D=\sigma(\frak D_{\tilde x_1})$, we have 
$\sigma(\frak D_{\tilde x_1})\in\Dec_l(\frak G_2)$ in the notation of 
Proposition 1.5(ii). Thus, by Proposition 1.5(ii), there exists 
$\tilde x_2\in\Sigma_{\tilde X_2}$, such that $\sigma(\frak D_{\tilde x_1})\subset\frak D_{\tilde x_2}$. 
By Proposition 1.5(i), such $\tilde x_2$ is unique. So, set $\tilde\phi(\tilde x_1)=\tilde x_2$, which 
has the desired properties. 

\noindent
(ii) The existence of $\frak D_{\tilde\phi(\tilde x_1),l}$ follows from 
the fact that $\sigma(\frak D_{\tilde x_1,l})\subset\sigma(\frak D_{\tilde x_1})\subset 
\frak D_{\tilde\phi(\tilde x_1)}$ and that $\sigma(\frak D_{\tilde x_1,l})$ is pro-$l$. 
Finally, 
consider the composite map of
$$\frak D_{\tilde x_1}\overset{\sigma}\to{\to}\frak D_{\tilde\phi(\tilde x_1)}\twoheadrightarrow 
\frak D_{\tilde\phi(\tilde x_1)}/\frak I_{\tilde\phi(\tilde x_1)}=G_{k_2(q_2(\tilde\phi(\tilde x_1)))}.
$$
Then, since $\cd_l(G_{k_2(q_2(\tilde\phi(\tilde x_1)))})=1$, the image of $\frak I_{\tilde x_1,l}$ in 
$G_{k_2(q_2(\tilde\phi(\tilde x_1)))}$ must be trivial by Proposition 1.1(i), as desired. 

\noindent
(iii) Immediate from the definitions. 

\noindent
(iv) We have $\cd_l(\sigma(\frak D_{\tilde x_1}))\leq\cd_l(\frak G_2)\leq 2<\infty$, where 
the second inequality follows from Lemma 1.4. Now, the assertion follows from Proposition 1.1(i). 

\noindent
(v) This follows from the fact that 
the defining property $\sigma(\frak D_{\tilde x_1})\subset \frak D_{\tilde\phi(\tilde x_1)}$ 
of $\tilde\phi$ is independent of $l$. 
\qed\enddemo

We shall consider the following conditions: 

\definition{Condition 1}
Either $\frak G_i=G_i$, $i=1,2$ or 
$\frak G_i=G_i^{(p_i')}$, $i=1,2$. We refer to the 
first and the second cases as 
the profinite and the prime-to-characteristic cases, respectively. 
(Observe that conditions (1)(2) in Proposition 2.2 are then satisfied 
for any prime number $l\neq p_1,p_2$.)
In particular, we have $\frak G_{k_i}=G_{k_i}$ in both cases. 
\enddefinition

\definition{Condition 2}
The map $\sigma:\frak G_1\to\frak G_2$ commutes with the projections $\pr_1, \pr_2$, i.e., 
it inserts into the following commutative diagram:

$$
\CD
1@>>>    \overline {\frak G}_1   @>>>     {\frak G}_1   @>\pr_1>>    G_{k_1}     @>>> 1\\
  @.        @V{\bar \sigma} VV            @V\sigma VV              @V\sigma_0 VV \\
1     @>>> \overline {\frak G}_2  @>>>   {\frak G}_2     @>\pr _2>>   G_{k_2}     @>>>   1
\endCD 
$$
where the rows are exact. 
\enddefinition 

\definition{Condition 3}
The map $\sigma:\frak G_1\to\frak G_2$ is an open homomorphism. 
\enddefinition

In the rest of this section, we assume that Condition 1 holds. Then: 

\proclaim
{Lemma 2.3} In the prime-to-characteristic case, 
Condition 2 automatically holds. 
In the profinite case, if $\sigma (\frak I_1)\subseteq \frak I_2$, 
then 
Condition 2 holds. 
\endproclaim

\demo
{Proof} In the prime-to-characteristic case,
the quotient $\pr_i: \frak G_i\twoheadrightarrow G_{k_i}$ 
coincides with 
$\frak G_i^{\ab}$ modulo the closure of the torsion subgroup. 
Thus, $\sigma$ commutes with the projections $\pr_1,\pr_2$. 

In the profinite case, assume that $\sigma (\frak I_1)\subseteq \frak I_2$. 
Then $\sigma$ 
induces naturally, by passing to the quotients $\frak G_i/\frak I_i$, 
a 
homomorphism 
$
\pi_1(X_1)\to \pi_1(X_2)$ between fundamental groups. 
The quotient $\pr_i: \frak G_i\twoheadrightarrow \pi_1(X_i)\twoheadrightarrow G_{k_i}$ 
coincides with 
$\frak \pi_1(X_i)^{\ab}$ modulo 
the torsion subgroup. 
Thus, $\sigma$ commutes with the projections $\pr_1,\pr_2$. 
\qed
\enddemo

In the rest of this 
section, 
we assume, moreover, that 
Condition 3 holds. 
Then note that, if 
Condition 2 also holds 
and 
if $\sigma_0:G_{k_1}\to G_{k_2}$ and $\bar\sigma: \overline{\frak G}_1\to\overline{\frak G}_2$ 
are homomorphisms induced by $\sigma$, then 
automatically $\sigma_0$ is open and injective and $\bar\sigma$ is open. 

\proclaim
{Lemma 2.4} (Invariance of the Characteristics) 
The equality $p_1=p_2$ holds.
\endproclaim

\demo
{Proof} By replacing $\frak G_2$ by the open subgroup $\sigma(\frak G_1)\subset\frak G_2$, 
we may and shall assume that $\sigma$ is surjective. 

In the profinite case, the assertion 
follows by considering the 
(pro-)$q$-parts of  
$\frak G_i^{\ab}$ for various prime numbers $q$. 
Indeed, for $i\in \{1,2\}$, the $p_i$-part of $\frak G_i^{\ab}$ modulo 
the closure of the torsion subgroup 
is not finitely generated, while the $l$-part of $\frak G_i^{\ab}$ 
modulo the closure of the torsion subgroup, 
for a prime number $l\neq p_i$, is finitely generated (and even cyclic), as follows from the 
structure of $\frak G_i^{\ab}$ given by global class field theory. 
(Note, however, that the $l$-torsion subgroup of $\frak G_i^{\ab}$ is infinite.) 
Thus, $\frak G_2$ being a quotient of $\frak G_1$ (via $\sigma$) we must have $p_1=p_2$.

In the prime-to-characteristic case, 
the assertion follows by considering the 
$q$-Sylow subgroups $\frak G_{i,q}$ of 
$\frak G_i$ for various prime numbers $q$. As $\sigma$ is assumed to be surjective, 
we may and shall take $\frak G_{2,q}=\sigma(\frak G_{1,q})$. 
Indeed, 
for $i\in \{1,2\}$, 
$\frak G_{i,p_i}$ is 
cyclic, while $\frak G_{i,l}$ 
for a prime number $l\neq p_i$ is non-cyclic. 
Accordingly, the surjection $\frak G_{1,p_1}\twoheadrightarrow \frak G_{2,p_1}$ 
cannot exist, unless $p_1=p_2$. Thus, we must have $p_1=p_2$. 
\qed
\enddemo

So, from now on, set $p\defeq p_1=p_2$. 

\definition
{Remark 2.5}\ The same argument used in the proof of (the prime-to-characteristic 
case of) Lemma 2.3 shows
that an open homomorphism $\sigma:\frak G_1\to \frak G_2$ between profinite 
groups automatically commute with the natural projections 
$\pr_i':\frak G_i\to G_{k_i}^{p'}$, induced by $\pr_i$, for $i=1,2$. Thus, we 
have a commutative diagram:
$$
\CD
\frak G_1 @>{\pr_1'}>> G_{k_1}^{p'} \\
   @V{\sigma} VV           @V \sigma_0' VV  \\
\frak G_2 @>{\pr_2'}>> G_{k_2}^{p'}, 
\endCD
$$
where the right column is automatically open and injective. 
The authors do not know, 
at least at the time of writing, whether or not 
Condition 2 follows from Conditions 1 and 3 
in general (i.e., even in the profinite case). 

\enddefinition

In the rest of this subsection we assume, moreover, that Condition 2 holds. 

\proclaim{Lemma 2.6} 
The map $\sigma$ induces a natural open homomorphism 
$\sigma': G_1^{(p')}\to G_2^{(p')}$, which commutes with the 
canonical projections $G_i^{(p')}\to G_{k_i}$; $i=1,2$. 
For $i=1,2$, let $\frak I'_i$ be the 
image of $\frak I_i\subset \frak G_i$ in $G_i^{(p')}$. 
Then $\sigma'(\frak I_1')\subset\frak I_2'$. Thus, $\sigma$ induces a 
natural open homomorphism 
$\tau': \pi_1(X_1)^{(p')}\to\pi_1(X_2)^{(p')}$, 
which commutes with the canonical projections 
$\pi_1(X_i)^{(p')}\to G_{k_i}$; $i=1,2$. 
In particular, we have $g_1\geq g_2$. 
\endproclaim

\demo{Proof}
The first assertion is clear. 
The second assertion follows from Proposition 2.2(ii)(iv). 
The third assertion follows from the second. 
Now, $\tau':\pi_1(X_1)^{(p')}\to\pi_1(X_2)^{(p')}$ induces an open 
homomorphism 
$\pi_1(\overline X_1)^{p'}\to\pi_1(\overline X_2)^{p'}$, 
hence an open homomorphism 
$\pi_1(\overline X_1)^{p',\ab}\to\pi_1(\overline X_2)^{p',\ab}$. 
Since the 
$\pi_1(\overline X_i)^{p',\ab}$ is a free 
${\hat\Bbb Z}^{p'}$-module of rank $2g_i$ for $i=1,2$, 
this implies the last assertion. 
\qed\enddemo

%

\proclaim
{Lemma 2.7} 
For a prime number $l\neq p$, the map 
$\phi=\phi_{\sigma,l}: \Sigma_{X_1,\sigma,l}\to\Sigma_{X_2}$ 
is almost surjective, i.e.,
$\Sigma _{X_2}\smallsetminus \phi (\Sigma  _{X_1,\sigma,l})$ 
is finite. In particular, $\Sigma _{X_1,\sigma,l}$ is infinite (hence, a fortiori, nonempty). 
\endproclaim

\demo
{Proof} Assume that the set $S\defeq \Sigma _{X_2}\smallsetminus \phi (\Sigma  _{X_1,\sigma,l})$ 
is infinite. Set $U_2\defeq X_2\smallsetminus S$. 

As in (the third assertion of) Lemma 2.6, then $\sigma$ induces an open homomorphism 
$\tau^{(l)}_1: \pi_1(X_1)^{(l)}\to\pi_1(U_2)^{(l)}$, 
which is a lifting of the homomorphism 
$\tau^{(l)}: \pi_1(X_1)^{(l)}\to\pi_1(X_2)^{(l)}$ 
induced by 
$\tau': \pi_1(X_1)^{(p')}\to\pi_1(X_2)^{(p')}$. 
We have a commutative diagram:
$$
\CD
1 @>>> \pi_1(\overline X_1)^l    @>>>  \pi_1(X_1)^{(l)} @>\pr_1>> G_{k_1}@>>> 1  \\
@.   @V \bar\tau^{l}_1VV           @V \tau^{(l)}_1 VV               @V\sigma_0 VV  \\
 1 @>>>\pi_1(\overline U_2)^l    @>>>     \pi_1(U_2)^{(l)}  @>\pr_2>> G_{k_2} @>>>
1 
\endCD 
$$
where ${\overline U_2}\defeq U_2\times _{k_2}{\bar k_2}$.
Since $\tau^{(l)}_1: \pi_1(X_1)^{(l)}\to\pi_1(U_2)^{(l)}$ is open and 
$\sigma_0: G_{k_1}\to G_{k_2}$ is (open and) injective, 
we see that $\bar\tau^{l}_1: \pi_1(\overline X_1)^l\to \pi_1(\overline U_2)^l$ is open. 
This is a contradiction, since $\pi_1(\overline X_1)^l$ is (topologically) finitely generated, 
while $\pi_1(\overline U_2)^l$ (hence $\bar\tau^{l}_1( \pi_1(\overline X_1)^l)$) 
is not (topologically) finitely generated, as $S$ is infinite. 
\qed
\enddemo

\proclaim{Lemma 2.8}
Let $\sigma_0: G_{k_1}\to G_{k_2}$ be the (open, injective) homomorphism 
induced by $\sigma$. Set $d_0\defeq [G_{k_2}: \sigma_0(G_{k_1})]$. 
Then: 

\noindent
{\rm (i)} The following diagram is commutative:
$$
\CD 
(\hat \Bbb Z^{p'})^{\times}    @=   (\hat \Bbb Z^{p'})^{\times} \\
@A\chi_{k_1}AA                                  @A\chi_{k_2}AA \\
G_{k_1}  @>\sigma_0>>   G_{k_2} \\
@A\pr_1AA                    @A\pr_2AA   \\
\frak G_1   @>\sigma>>  \frak G_2  \\
\endCD 
$$
where $\chi_{k_i}$ is the cyclotomic character of $G_{k_i}$ for $i=1,2$.

\noindent
{\rm (ii)} We have 
$\sharp (k_1)=\sharp (k_2)^{d_0}$ 
and 
$\sigma_0(\varphi_{k_1})=\varphi_{k_2}^{d_0}$. 
\endproclaim

\demo{Proof} 
(i) Since the bottom square is commutative by the definition of $\sigma_0$, 
we only have to prove that the top square is commutative. 
As $G_{k_2}$ is (topologically) generated by $\varphi_{k_2}$, 
we may write $\sigma_0(\varphi_{k_1})=\varphi_{k_2}^{\alpha}$, 
where $\alpha\in{\hat\Bbb Z}$. Now, the desired commutativity 
$\chi_{k_2}\circ\sigma_0=\chi_{k_1}$ is equivalent to 
saying that 
$\chi_{k_2}(\sigma_0(\varphi_{k_1}))=\chi_{k_1}(\varphi_{k_1})$ 
(as $G_{k_1}$ is (topologically) generated by $\varphi_{k_1}$). 
Since $\chi_{k_1}(\varphi_{k_1}) =\sharp(k_1)=p^{[k_1:\Bbb F_p]}$ and 
$$\chi_{k_2}(\sigma_0(\varphi_{k_1}))=
\chi_{k_2}(\varphi_{k_2}^{\alpha})=
\chi_{k_2}(\varphi_{k_2})^{\alpha}=\sharp(k_2)^\alpha
=p^{\alpha[k_2:\Bbb F_p]},$$
the desired commutativity is thus equivalent to the 
equality $\alpha[k_2:\Bbb F_p]=[k_1:\Bbb F_p]$ in $\hat \Bbb Z$. 
(Here, note that the homomorphism $\hat\Bbb Z\to (\hat\Bbb Z^{p'})^{\times}$, 
$\beta\mapsto p^\beta$ is injective by 
[Chevalley], Th\'eor\`eme 1.) 
In particular, 
it suffices to prove the desired commutativity on an open subgroup 
$H\subset G_{k_1}$. 
Indeed, set $m\defeq[G_{k_1}:H]$. Then, 
since $\varphi_{k_1}^m$ is the Frobenius element for $H$, 
the commutativity on $H$ is equivalent to the equality 
$m\alpha[k_2:\Bbb F_p]=m[k_1:\Bbb F_p]$ in $\hat\Bbb Z$, 
which implies 
$\alpha[k_2:\Bbb F_p]=[k_1:\Bbb F_p]$, as desired. 
Thus, by replacing $\frak G_1$ and $\frak G_2$ by suitable open subgroups, 
we may and shall assume 
that $g_2>0$. 

Next, for each prime number $l\neq p$ and $i\in\{1,2\}$, let 
$\chi_{k_i,l}: G_{k_i}\to\Bbb Z_l^{\times}$ denote the 
$l$-adic cyclotomic character. Thus, corresponding to the 
decomposition $(\hat\Bbb Z^{p'})^{\times}=
\prod_{l\neq p}\Bbb Z_l^{\times}$, we have 
$\chi_{k_i}=(\chi_{k_i,l})_{l\neq p}$. 
We have to prove that $\chi_{k_2}\circ\sigma_0=\chi_{k_1}$, 
which is equivalent to saying that $\chi_{k_2,l}\circ\sigma_0=\chi_{k_1,l}$ 
for all $l\neq p$. 

We shall first prove that the last equality holds up to torsion. More precisely, 
denote by $\bar \chi_{k_i,l}$ 
the composite of $G_{k_i}\overset{\chi_{k_i,l}}\to{\to}\Bbb Z_l^{\times}\twoheadrightarrow 
\Bbb Z_l^{\times}/(\Bbb Z_l^{\times})^{\tor}$. By Lemma 2.7, we can take 
$\tilde x_1\in \Sigma_{\tilde X_1,\sigma,l}\neq\emptyset$. Set $\tilde x_2\defeq\tilde \phi(\tilde x_1)$. 
Let $x_i$ denote the image of $\tilde x_i$ in $\Sigma_{X_i}$ for $i=1,2$. By Proposition 2.2(ii), 
we have $\sigma: \frak D_{\tilde x_1,l} \to \frak D_{\tilde x_2,l} $ and 
$\sigma: \frak I_{\tilde x_1,l} \to \frak I_{\tilde x_2,l} $, which are injective by Proposition 1.1(i). 
This implies that $\chi_{k_2,l}\circ\sigma_0=\chi_{k_1,l}$ holds on the image of 
$\frak D_{\tilde x_1,l}$ in $G_{k_1}$, which is an open subgroup of 
the $l$-Sylow subgroup $G_{k_1,l}$ of $G_{k_1}$. 
As $\Bbb Z_l^{\times}/(\Bbb Z_l^{\times})^{\tor}\simeq \Bbb Z_l$ is torsion-free and pro-$l$, 
this implies that 
$\bar \chi_{k_2,l}\circ\sigma_0=\bar\chi_{k_1,l}$. 

In particular, we have 
$\bar\chi_{k_2,l}(\sigma_0(\varphi_{k_1}))=\bar\chi_{k_1,l}(\varphi_{k_1})$. 
This implies the equality $\sharp(k_2)^{\alpha}=\sharp(k_1)$ in 
$\Bbb Z_l^{\times}/(\Bbb Z_l^{\times})^{\tor}\simeq \Bbb Z_l$. Since $p\in \Bbb Z_l^{\times}$ is 
not a torsion, this last equality shows that 
$\alpha_l[k_2:\Bbb F_p]=[k_1:\Bbb F_p]$ in $\Bbb Z_l$. Here, 
corresponding to the decomposition $\hat\Bbb Z=\prod_{\text{$l$: prime}}\Bbb Z_l$, 
we write $\alpha=(\alpha_l)_{\text{$l$: prime}}$. 
Or, equivalently, we have 
$$\alpha[k_2:\Bbb F_p]=[k_1:\Bbb F_p]+\iota_p(\epsilon)$$
in $\hat\Bbb Z$, where $\iota_p:\Bbb Z_p\hookrightarrow\hat\Bbb Z$ is a natural injection 
and 
$\epsilon
\defeq \alpha_p[k_2:\Bbb F_p]-[k_1:\Bbb F_p]
\in \Bbb Z_p$. 

On the other hand, by Lemma 2.6, we get an open homomorphism 
$\pi_1 (\overline X_1)^{p'}\to\pi_1(\overline X_2)^{p'}$, 
hence a surjection 
$\pi_1 (\overline X_1)^{p',\ab}\otimes_{\Bbb Z}\Bbb Q\twoheadrightarrow 
\pi_1 (\overline X_2)^{p',\ab}\otimes_{\Bbb Z}\Bbb Q$, 
which are Galois-compatible with respect to 
$\sigma_0: G_{k_1}\to G_{k_2}$. 
For each $i=1,2$, let $P_i(T)$ be the characteristic polynomial of 
$\varphi_{k_i}^{[k_{i'}: \Bbb F_p]}$ on the free $\hat\Bbb Z^{{p'}}$-module 
$\pi_1 (\overline X_i)^{p',\ab}$ 
(of rank $2g_i$), where $i'$ is defined by $\{i,i'\}=\{1,2\}$. 
Then, it is known that $P_i(T) \in \Bbb Z[T]$. 

Write $\rho_i$ for the natural representation 
$G_{k_i}\to
\Aut_{\hat\Bbb Z^{p'}}(\pi_1 (\overline X_i)^{p',\ab})$. 
Let $R_{\Bbb Q}$ be the (commutative) $\Bbb Q$-subalgebra of 
$\End_{\hat\Bbb Z^{p'}\otimes_{\Bbb Z}\Bbb Q}(\pi_1 (\overline X_2)^{p',\ab}\otimes_{\Bbb Z}\Bbb Q)$ 
generated by $\rho_2(G_{k_2})$. We have 
$P_2(\rho_2(\varphi_{k_2}^{[k_1:\Bbb F_p]}))=0$ 
in $R_{\Bbb Q}$. By the Galois-compatibility, we also have 
$P_1(\rho_2(\sigma_0(\varphi_{k_1}^{[k_2:\Bbb F_p]})))=0$ 
in $R_{\Bbb Q}$. These identities imply that 
both of $\rho_2(\varphi_{k_2}^{[k_1:\Bbb F_p]})$, 
$\rho_2(\sigma_0(\varphi_{k_2}^{[k_2:\Bbb F_p]})) \in R_{\Bbb Q}$ 
are algebraic over $\Bbb Q$, hence so is the ratio 
$$
\rho_2(\sigma_0(\varphi_{k_1}^{[k_2:\Bbb F_p]})
(\varphi_{k_2}^{[k_1:\Bbb F_p]})^{-1})
=\rho_2(\varphi_{k_2}^{\alpha[k_2:\Bbb F_p]-[k_1:\Bbb F_p]})
=\rho_2(\varphi_{k_2}^{\iota_p(\epsilon)})\defeq \eta$$
in $R_{\Bbb Q}$. So, take a monic polynomial $Q(T)\in\Bbb Q[T]$ 
satisfying $Q(\eta)=0$ in $R_{\Bbb Q}$. 
Set $b\defeq\deg(Q)$. 

Let $l$ be a prime number $\neq p$, and 
$R_{l,\Bbb Q}$ the image of 
$R_{\Bbb Q}
$ 
in 
$\End_{\Bbb Q_l}(\pi_1 (\overline X_2)^{l,\ab} 
\otimes_{\Bbb Z}\Bbb Q)$. 
Then observe that 
the image $\eta_l$ of $\eta$ 
in $R_{l,\Bbb Q}\subset \End_{\Bbb Q_l}(\pi_1 (\overline X_2)^{l,\ab}\otimes_{\Bbb Z}\Bbb Q)$ 
is a pro-$p$ element of 
$\End_{\Bbb Z_l}(\pi_1 (\overline X_2)^{l,\ab})$, 
hence a torsion element of $p$-power order. So, let $p^{a_l}$ be 
the order of $\eta_l$. As $Q(\eta_l)=0$ in the commutative 
$\Bbb Q$-algebra $R_{l,\Bbb Q}$, we conclude: 
$\frac{p-1}{p}p^{a_l}\leq \varphi(p^{a_l})\leq b$, 
where $\varphi$ stands for Euler's function. 
(Use the fact $\Bbb Q\hookrightarrow R_{l,\Bbb Q}$, which follows from 
$g_2>0$.) 
Thus, $a_l$ is bounded: there exists $a\geq 0$ such that $a_l\leq a$ for 
all $l\neq p$. Namely, $(\eta_l)^{p^a}=1$ for all $l\neq p$. 

Set $\zeta_l\defeq\det(\eta_l)$, where the determinant is taken as an element of 
$\End_{\Bbb Q_l}(\pi_1 (\overline X_2)^{l,\ab}\otimes_{\Bbb Z}\Bbb Q)$. 
Since $\det$ is a homomorphism, we have 
$(\zeta_l)^{p^a}=1$ for all $l\neq p$. Set $\zeta\defeq(\zeta_l)_{l\neq p}$
in $(\hat\Bbb Z^{p'})^{\times}=\prod_{l\neq p}\Bbb Z_l^{\times}$
Now, by construction, we have 
$$\zeta=\chi_{k_2,l}^{g_2}(\varphi_{k_2}^{\iota_p(\epsilon)})
=\sharp(k_2)^{g_2 \iota_p(\epsilon)},$$
hence 
$\sharp(k_2)^{p^ag_2 \iota_p(\epsilon)}=1$ 
in $(\hat \Bbb Z^{p'})^{\times}$. 
Since the homomorphism $\hat\Bbb Z\to(\hat \Bbb Z^{p'})^{\times}$, 
$\beta\mapsto p^{\beta}$ is injective, this last equality 
forces 
$[k_2:\Bbb F_p]p^ag_2 \iota_p(\epsilon)=0$
in $\hat\Bbb Z$. 
As $[k_2:\Bbb F_p]p^ag_2>0$, this implies $\iota_p(\epsilon)=0$. Namely, we have 
$\alpha[k_2:\Bbb F_p]=[k_1:\Bbb F_p]$
in $\hat \Bbb Z$, 
as desired. 

\noindent
(ii) As in the proof of (i), set $\sigma_0(\varphi_{k_1})=\varphi_{k_2}^{\alpha}$. 
Since $G_{k_2}\simeq\hat\Bbb Z$ and $[G_{k_2}:\sigma_0(G_{k_1})]=d_0$, we must have 
$\alpha=d_0u$, where $u\in\hat\Bbb Z^{\times}$. Now, since 
$\alpha[k_2:\Bbb F_p]=[k_1:\Bbb F_p]$ by (i), we obtain 
$d_0u[k_2:\Bbb F_p]=[k_1:\Bbb F_p]$, hence 
$u=[k_1:\Bbb F_p]/(d_0[k_2:\Bbb F_p])\in\Bbb Q_{>0}\ (\subset \hat\Bbb Z\otimes_{\Bbb Z}\Bbb Q)$. 
Since $\hat\Bbb Z^{\times}\cap\Bbb Q_{>0}=\{1\}$, we conclude $u=1$. 
Thus, $d_0[k_2:\Bbb F_p]=[k_1:\Bbb F_p]$ and $\sigma_0(\varphi_{k_1})=\varphi_{k_2}^{d_0}$, 
as desired. 
\qed\enddemo

\proclaim{Lemma 2.9} 
For 
each prime number $l\neq p$, the map 
$\tilde\phi_{\sigma,l}: \Sigma_{\tilde X_1,\sigma,l}\to\Sigma_{\tilde X_2}$ 
is surjective. In particular, the map 
$\phi_{\sigma,l}: \Sigma_{X_1,\sigma,l}\to\Sigma_{X_2}$ 
is surjective. 
\endproclaim

\demo{Proof} 
As in the proof of Lemma 2.7, set $S\defeq \Sigma _{X_2}\smallsetminus \phi (\Sigma  _{X_1,\sigma,l})$ 
and $U_2\defeq X_2\smallsetminus S$. Thus, by Lemma 2.7, $S$ is a finite set. 
Let $r<\infty$ be the cardinality of $S(\bar k_2)$. 
Then $\sigma$ induces an open homomorphism 
$\tau^{(l)}_1: \pi_1(X_1)^{(l)}\to\pi_1(U_2)^{(l)}$, 
which is a lifting of the homomorphism 
$\tau^{(l)}: \pi_1(X_1)^{(l)}\to\pi_1(X_2)^{(l)}$ 
induced by 
$\tau': \pi_1(X_1)^{(p')}\to\pi_1(X_2)^{(p')}$ 
in Lemma 2.6. We have a commutative diagram:
$$
\CD
1 @>>> \pi_1(\overline X_1)^l    @>>>  \pi_1(X_1)^{(l)} @>\pr_1>> G_{k_1}@>>> 1  \\
@.   @V \bar\tau^{l}_1VV           @V \tau^{(l)}_1 VV               @V\sigma_0 VV  \\
 1 @>>>\pi_1(\overline U_2)^l    @>>>     \pi_1(U_2)^{(l)}  @>\pr_2>> G_{k_2} @>>>
1 
\endCD 
$$
where ${\overline U_2}\defeq U_2\times _{k_2}{\bar k_2}$.
Since $\tau^{(l)}_1: \pi_1(X_1)^{(l)}\to\pi_1(U_2)^{(l)}$ is open and 
$\sigma_0: G_{k_1}\to G_{k_2}$ is (open and) injective, 
we see that $\bar\tau^{l}_1: \pi_1(\overline X_1)^l\to \pi_1(\overline U_2)^l$ is open. 
The open homomorphism $\bar\tau^{l}_1: \pi_1(\overline X_1)^l\to \pi_1(\overline U_2)^l$ 
induces an open homomorphism 
$\bar\tau^{l,\ab}_1: \pi_1(\overline X_1)^{l,ab}\to \pi_1(\overline U_2)^{l,ab}$. 
This last homomorphism is, by construction, Galois-compatible with respect to 
$\sigma_0: G_{k_1}\to G_{k_2}$. In other words, if we regard 
$\pi_1(\overline U_2)^{l,ab}$ as a $G_{k_1}$-module via $\sigma_0$, 
$\bar\tau^{l,\ab}_1$ is a homomorphism as $G_{k_1}$-modules. 

The absolute values of eigenvalues of $\varphi_{k_1}\in G_{k_1}$ 
in $\pi_1(\overline X_1)^{l,ab}$ 
are 
all $\sharp(k_1)^{1/2}$, with multiplicity $2g_1$. On the other hand, 
by the second assertion of 
Lemma 2.8(ii), the absolute values of eigenvalues of $\varphi_{k_1}$ in 
$\pi_1(\overline U_2)^{l,ab}$ are the same as those of $\varphi_{k_2}^{d_0}$, 
which are $\sharp(k_2)^{d_0/2}$ with multiplicity $2g_2$ and 
$\sharp(k_2)^{d_0}$ with multiplicity $\max(r-1,0)$. By the first assertion of 
Lemma 2.8(ii), they coincide with $\sharp(k_1)^{1/2}$ and $\sharp(k_1)$, 
respectively. Thus, we conclude $r\leq 1$. 
However, if $r\neq 0$, by replacing $\frak G_1, \frak G_2$ with suitable open subgroups, 
we may assume that $r>1$, a contradiction. So, we have established $r=0$. 

To prove the surjectivity of $\tilde\phi_{\sigma,l}$, we may replace freely 
$\frak G_1, \frak G_2$ by open subgroups $\frak H_1,\frak H_2$, respectively, such that 
$\sigma(\frak H_1)\subset \frak H_2$. 
(Indeed, the map $\tilde\phi_{\sigma,l}: \Sigma_{\tilde X_1,\sigma,l}\to\Sigma_{\tilde X_2}$ 
remains unchanged.) In particular, we may assume that 
$\sigma:\frak G_1\to\frak G_2$ is surjective. Then, the surjectivity of 
$\tilde\phi_{\sigma,l}: \Sigma_{\tilde X_1,\sigma,l}\to\Sigma_{\tilde X_2}$ 
is equivalent to the surjectivity of 
$\phi_{\sigma,l}: \Sigma_{X_1,\sigma,l}\to\Sigma_{X_2}$, 
which is then equivalent to $r=0$. This completes the proof. 
\qed\enddemo

\subhead
\S 3. Rigid homomorphisms between Galois groups
\endsubhead

In this section we shall investigate a class of homomorphisms 
between (geometrically prime-to-characteristic quotients of) absolute 
Galois groups of function fields of curves over finite 
fields, which we call rigid. We follow the notations in $\S1$ and $\S2$. 
In particular, we follow the Notation at the beginning of subsection 2.2. 
We assume that Condition 3 holds. 

\comment
For $i\in \{1,2\}$, let $k_i$ be a finite field of characteristic 
$p_i> 0$. Let $X_i$ be a smooth, proper, geometrically connected 
curve of genus $g_i$ over $k_i$. Let $K_i=K_{X_i}$ be the function field of $X_i$ 
and fix 
an algebraic closure $\overline K_i$ of $K_i$. Let 
$K_i^{\sep}$ be a separable closure of $K_i$ in $\overline K_i$, and 
$\bar k_i$ the algebraic closure of $k_i$ in $\overline K_i$. 
Following the notations 
in $\S 1$, we will write $G_i\defeq G_{K_i}=\Gal (K_i^{\sep}/K_i)$
for the absolute Galois group of $K_i$, and $\bar{G}_i\defeq 
G_{K_i\bar k_i}=\Gal (K_i^{\sep}/K_i\bar k_i)$
for the absolute Galois group of $K_i\bar k_i$. 

Let $N_i$ be a normal closed subgroup of $G_i$ and set 
$\frak G_i\defeq G_i/N_i$. Let $\tilde K_i$ denote the Galois extension of $K_i$ 
corresponding to $N_i$, i.e., $\tilde K_i\defeq (K_i^{\sep})^{N_i}$. 
Let $\bar \frak G_i$ be the image of $\bar G_i$ in $\frak G_i$, and set 
$\frak G_{k_i}\defeq \frak G_i/\bar \frak G_i$, which is a quotient of 
$G_{k_i}$. For $i=1,2$, let us denote by $\varphi_{k_i}$ the image 
in $\frak G_{k_i}$ of the $\sharp(k_i)$-th power Frobenius element of $G_{k_i}$. 
Write $\tilde X_i$ for the integral closure of $X_i$ in $\tilde K_i$. 
The Galois group $\frak G_i$ acts naturally on the set 
$\Sigma _{\tilde X_i}$, and the quotient  $\Sigma _{\tilde X_i}/\frak G_i$ 
is naturally identified with $\Sigma _{X_i}$. 
For a point $\tilde x_i\in \Sigma _{\tilde X_i}$, 
let $\frak G_i\supset\frak D_{\tilde x_i}\supset\frak I_{\tilde x_i}$ be the 
decomposition and the inertia subgroups, respectively, and set 
$\frak G_{k_i(x_i)}\defeq\frak D_{\tilde x_i}/\frak I_{\tilde x_i}$. 


Now, let 
$$\sigma:\frak G_1\to \frak G_2$$ 
be 
an open 
homomorphism between profinite groups. 
\endcomment

\definition
{Definition 3.1} (Rigid Homomorphisms) 
(i) We say that $\sigma :\frak G_1\to \frak G_2$ is strictly rigid,  
if there exists a map
$$\tilde \phi :\Sigma_{\tilde X_1} \to \Sigma_{\tilde X_2},$$
such that 
$$\sigma (\frak D_{\tilde x_1})=\frak D_{\tilde \phi(\tilde x_1)}$$
for each $\tilde x_1\in \Sigma_ {\tilde X_1}$. 

\noindent
(ii) We say that $\sigma :\frak G_1\to \frak G_2$ is rigid, 
if there exist open subgroups $\frak H_1\subset \frak G_1$,  
$\frak H_2\subset \frak G_2$, such that $\sigma(\frak H_1)\subset\frak H_2$ 
and that $\frak H_1\overset{\sigma}\to{\to}\frak H_2$ is strictly rigid. 
(Here, $\frak H_i$ is considered as a quotient of the absolute Galois group 
that is the inverse image in $G_i$ of $\frak H_i\subset\frak G_i$.) 

\noindent
(iii) Define $\Hom(\frak G_1,\frak G_2)^{\rig}\subset\Hom(\frak G_1,\frak G_2)$ 
to be the set of rigid (hence continuous open) homomorphisms $\frak G_1\to\frak G_2$. 
\enddefinition

\definition
{Remark 3.2} 
{\rm (i)} Consider a commutative diagram of maps between profinite groups:
$$
\CD
\frak G_1 @>\sigma>> \frak G_2 \\
@VVV     @VVV \\
\frak G_1' @>\sigma'>> \frak G_2' \\
\endCD
$$
where the vertical arrows are surjective. 
Then, if $\sigma: \frak G_1\to \frak G_2$ is strictly rigid (resp. rigid),
then $\sigma': \frak G_1'\to \frak G_2'$ is strictly rigid (resp. rigid).

\noindent
{\rm (ii)} Let $\frak H_2$ be an open subgroup of $\frak G_2$ and 
$\frak H_1\defeq\sigma^{-1}(\frak H_2)$. Then, 
if $\sigma: \frak G_1\to\frak G_2$ is strictly rigid (resp. rigid), 
then the natural homomorphism $\frak H_1\to \frak H_2$ induced by $\sigma$ 
is strictly rigid (resp. rigid). 

\noindent
{\rm (iii)} Assume that $\sigma: \frak G_1\to\frak G_2$ 
is strictly rigid with respect to 
$\tilde \phi :\Sigma _{\tilde X_1}
\to \Sigma _{\tilde X_2}$. Then, 
if $\tilde\phi$ is surjective, then 
$\sigma$ is surjective. Indeed, this follows 
immediately from the fact, by Chebotarev's density theorem, that 
$\frak G_2$ 
is (topologically) generated by its decomposition subgroups. 


\noindent
(iv) As in Proposition 2.2, let $l$ be a prime number $\neq p_1,p_2$, and put 
the following two assumptions: 
{\rm (1)} $N_2^l=N_2$, or, equivalently, $\tilde K_2$ admits 
no $l$-cyclic extension; and 
{\rm (2)} $\tilde K_2$ contains a primitive $l$-th roots of unity. 

If $\sigma$ is strictly rigid, with respect to 
$\tilde \phi: \Sigma_{\tilde X_1} \to \Sigma_{\tilde X_2}$, 
then we must have $\Sigma_{\tilde X_1,\sigma,l}=\Sigma_{\tilde X_1}$ 
and $\tilde\phi=\tilde\phi_{\sigma,l}$. In particular, 
then $\tilde\phi$ is unique and Galois-equivariant with respect to $\sigma$, 
hence induces naturally a map
$\phi\  (=\phi_{\sigma,l}):\Sigma _{X_1} \to \Sigma_ {X_2}$. 

If $\sigma$ is rigid, then we must have $\Sigma_{\tilde X_1,\sigma,l}=\Sigma_{\tilde X_1}$, 
and, if we set $\tilde\phi\defeq\tilde\phi_{\sigma,l}$, then we have 
$\sigma (\frak D_{\tilde x_1})\underset \open \to \subset \frak D_{\tilde \phi(\tilde x_1)}$
for each $\tilde x_1\in \Sigma_ {\tilde X_1}$. 
The map $\tilde\phi$ is uniquely characterized by this property, 
and Galois-equivariant with respect to $\sigma$, 
hence induces naturally a map
$\phi\  (=\phi_{\sigma,l}):\Sigma _{X_1} \to \Sigma_ {X_2}$. 
\enddefinition

In the rest of this section, we 
assume that 
Condition 1 holds. 

\definition{Definition 3.3} 
(i) Let $\gamma: K_2\to K_1$ be a homomorphism of fields, which defines 
an extension $K_1/K_2$ of fields. Set $p\defeq p_1=p_2$. Then 
we say that $\gamma$ is admissible, if 
the extension $K_1/K_2$ appears in the extensions of $K_2$ corresponding to 
the open subgroups of $\frak G_2$. More precisely, 
in the profinite (resp. prime-to-characteristic) case, 
we say that $\gamma$ is 
admissible, if the extension $K_1/K_2$ is separable (resp. if the extension 
$K_1/K_2$ is separable and the Galois closure of the extension 
$K_1\bar k_1/K_2\bar k_2$ is of degree prime to $p$). 

Equivalently, $\gamma: K_2\to K_1$ is admissible if and only if it extends 
to an isomorphism $\tilde\gamma: \tilde K_2\isom \tilde K_1$. 

\noindent
(ii) Define $\Hom(K_2,K_1)^{\adm}\subset\Hom(K_2,K_1)$ 
to be the set of admissible homomorphisms $K_2\to K_1$. 
\enddefinition

Our aim in this section is to prove the following. 

\proclaim
{Theorem 3.4} 
The natural map $\Hom(K_2,K_1)\to\Hom(\frak G_1, \frak G_2)/\Inn(\frak G_2)$ 
induces a bijection 
$$\Hom(K_2,K_1)^{\adm}\isom\Hom(\frak G_1, \frak G_2)^{\rig}/\Inn(\frak G_2).$$ 
More precisely, 

\noindent
{\rm (i)} If $\gamma:K_2\to K_1$ is an admissible homomorphism between fields, then 
the homomorphism $\frak G_1\to\frak G_2$ induced by $\gamma$ (up to inner automorphisms) 
is rigid. 

\noindent
{\rm (ii)} If $\sigma:\frak G_1\to \frak G_2$ is a rigid 
homomorphism between profinite groups, then there 
exists a unique isomorphism $\tilde \gamma: \tilde K_2\to
\tilde K_1$ of fields, such that $\tilde\gamma\circ \sigma (g_1)=g_1
\circ \tilde \gamma$, for all $g_1\in \frak G_1$, which induces an 
admissible homomorphism $K_2\to K_1$. 
\endproclaim

\definition{Remark 3.5}
(i) By local theory for the Isom-form, 
any isomorphism $\frak G_1\isom\frak G_2$ is strictly rigid. In particular, 
we have $\Isom(\frak G_1,\frak G_2)
\subset\Hom(\frak G_1,\frak G_2)^{\rig}$. Thus, Theorem 3.4 can be viewed as 
a generalization of the Isom-form: 
$$\Isom(K_2,K_1)\isom \Isom(\frak G_1,\frak G_2)/\Inn(\frak G_2),$$
which is the main theorem of [Uchida2] (resp. [ST]) in the profinite 
(resp. prime-to-characteristic) case. 

\noindent
(ii) Let $\gamma: K_2^{\perf}\to K_1^{\perf}$ be a homomorphism of fields, which defines 
an extension $K_1^{\perf}/K_2^{\perf}$ of fields. Set $p\defeq p_1=p_2$. Then 
we say that $\gamma$ is admissible, if 
the extension $K_1^{\perf}/K_2^{\perf}$ appears in the extensions of $K_2^{\perf}$ 
corresponding to the open subgroups of $\frak G_2$, which is regarded as a quotient 
of the absolute Galois group $G_{K_2^{\perf}}=G_{K_2}$. More precisely, 
in the profinite (resp. prime-to-characteristic) case, 
$\gamma$ is always admissible (resp. admissible if and only if the extension 
the Galois closure of the extension 
$K_1^{\perf}\bar k_1/K_2^{\perf}\bar k_2$ is of degree prime to $p$). 
Define $\Hom(K_2^{\perf},K_1^{\perf})^{\adm}\subset\Hom(K_2^{\perf},K_1^{\perf})$ 
to be the set of admissible homomorphisms $K_2^{\perf}\to K_1^{\perf}$. Then 
the natural map $\Hom(K_2^{\perf},K_1^{\perf})\to\Hom(\frak G_1, \frak G_2)/\Inn(\frak G_2)$ 
induces a bijection 
$$\Hom(K_2^{\perf},K_1^{\perf})^{\adm}
/\Frob^{\Bbb Z}\ \isom\ \Hom(\frak G_1, \frak G_2)^{\rig}/\Inn(\frak G_2).$$ 
Indeed, this follows from Theorem 3.4, since the natural map 
$\Hom(K_2,K_1)\to\Hom(K_2^{\perf},K_1^{\perf})
$ 
induces 
$$\Hom(K_2,K_1)^{\adm}\isom\Hom(K_2^{\perf},K_1^{\perf})^{\adm}/\Frob^{\Bbb Z}.$$
\enddefinition

The rest of this section will be devoted to the proof of Theorem 3.4. 

First, to prove (i), let $\gamma:K_2\to K_1$ be an admissible homomorphism. 
Then, by the definition of admissibility, the extension $K_1/K_2$ is 
isomorphic to some extension $L/K_2$ corresponds to an open subgroup 
$\frak H_2$ of $\frak G_2$. Set $\frak H_1\defeq\frak G_1$. 
Now, let $\sigma:\frak G_1\to\frak G_2$ be the homomorphism induced by 
$\gamma$ (up to conjugacy). Then it is easy to see that $\sigma$ restricts 
to an isomorphism $\frak H_1\isom\frak H_2$ (corresponding to 
the isomorphism $L\isom K_1$), which is strictly rigid. Thus, $\sigma$ 
is rigid, as desired. 

Next, to prove (ii), let $\sigma:\frak G_1\to\frak G_2$ be a rigid homomorphism. 
By definition, 
there exist open subgroups $\frak H_1\subset \frak G_1$,  
$\frak H_2\subset \frak G_2$, such that $\sigma(\frak H_1)\subset\frak H_2$ 
and that $\frak H_1\overset{\sigma}\to{\to}\frak H_2$ is strictly rigid 
with respect to, say, $\tilde\phi: \Sigma_{\tilde X_1}\to\Sigma_{\tilde X_2}$. 
Then, by Remark 3.2(iv), $\tilde\phi$ is Galois-equivariant with respect 
to $\sigma: \frak G_1\to\frak G_2$ (i.e., not only with respect to 
$\sigma: \frak H_1\to\frak H_2$), and, for each $\tilde x_1\in
\Sigma_{\tilde X_1}$, we have 
$\sigma (\frak D_{\tilde x_1})\underset \open \to \subset \frak D_{\tilde \phi(\tilde x_1)}$ 
and $\sigma (\frak D_{\tilde x_1}\cap\frak H_1)= \frak D_{\tilde \phi(\tilde x_1)}\cap\frak H_2$. 

\proclaim
{Lemma 3.6} Condition 2 holds for $\sigma: \frak G_1\to\frak G_2$. 
\endproclaim

\demo
{Proof} By Proposition 2.1(v), we have $\sigma(\frak I_{\tilde x_1})\subset\frak I_{
\tilde\phi(\tilde x_1)}$ for each $\tilde x_1\in\Sigma_{\tilde X_1}$. In particular, 
we have $\sigma(\frak I_1)\subset\frak I_2$. Now, the assertion follows from 
Lemma 2.3. 
\qed\enddemo

Thus, we may apply Lemmas 2.6--2.9 to $\sigma$. Further, we have the 
following: 

\proclaim{Lemma 3.7}
We have 
$\sigma(\frak H_1)=\frak H_2$ and 
$\frak H_1=\sigma^{-1}(\frak H_2)$. 
\endproclaim

\demo{Proof}
By Lemma 2.9, $\tilde\phi$ is surjective, hence, by Remark 3.2(iii), 
$\sigma: \frak H_1\to\frak H_2$ is surjective, that is, 
$\sigma(\frak H_1)=\frak H_2$. 

Next, let 
$X_{1,\frak H_1}\to X_{1,\frak \sigma^{-1}(\frak H_2)}\to X_1$ and 
$X_{2,\frak H_2}\to X_2$ 
be (finite, generically \'etale) covers corresponding to 
open subgroups $\frak H_1\subset\sigma^{-1}(\frak H_2)\subset\frak G_1$ 
and $\frak H_2\subset\frak G_2$, respectively. 
Suppose that $\frak H_1\subsetneq\sigma^{-1}(\frak H_2)$. 
Then, by Chebotarev's density theorem, there exists $\tilde x_1\in \Sigma_{\tilde X_1}$ 
such that $k(x_{1,\frak H_1})\supsetneq k(x_{1,\sigma^{-1}(\frak H_2)})$, where 
$x_{1,\frak H_1}$ and $x_{1,\sigma^{-1}(\frak H_2)}$ denote the 
images of $\tilde x_1$ in $\Sigma_{1,\frak H_1}$ and $\Sigma_{1,\sigma^{-1}(\frak H_2)}$, 
respectively. 
Set $\tilde x_2\defeq\tilde\phi(\tilde x_1)\in\Sigma_{\tilde X_2}$. 
We have 
$\sigma(\frak D_{\tilde x_1})\subset \frak D_{\tilde x_2}$, 
hence 
$$\sigma(\frak D_{\tilde x_1}\cap\frak H_1)\subset 
\sigma(\frak D_{\tilde x_1}\cap\frak \sigma^{-1}(\frak H_2))\subset 
\frak D_{\tilde x_2}\cap\frak H_2.$$
Now, since $\frak H_1\overset{\sigma}\to{\to}\frak H_2$ is strictly rigid, 
we must have 
$$\sigma(\frak D_{\tilde x_1}\cap\frak H_1)=
\sigma(\frak D_{\tilde x_1}\cap\frak \sigma^{-1}(\frak H_2))=
\frak D_{\tilde x_2}\cap\frak H_2.$$
By Proposition 2.1(iii), this implies that 
$\sharp(k(x_{1,\frak H_1}))
=\sharp(k(x_{2,\frak H_2}))
=\sharp(k(x_{1,\sigma^{-1}(\frak H_2)})$, 
where $x_{2,\frak H_2}$ denotes the image of $\tilde x_2$ in $\Sigma_{X_{2,\frak H_2}}$. 
This contradicts $k(x_{1,\frak H_1})\supsetneq k(x_{1,\sigma^{-1}(\frak H_2)})$. 
\qed\enddemo

First, we shall treat the special case that $\sigma: \frak G_1\to\frak G_2$ 
is strictly rigid (hence, in particular, surjective). 

\proclaim
{Lemma 3.8} Assume that $\sigma: \frak G_1\to\frak G_2$ 
is strictly rigid. Then: 

\noindent
{\rm (i)} We have $g_1=g_2$. 

\noindent
{\rm (ii)} The map $\phi :\Sigma _{X_1}\to \Sigma _{X_2}$ is bijective. 
\endproclaim

\demo
{Proof} By Lemma 2.6, the homomorphism $\sigma$ induces naturally a commutative
diagram:

$$
\CD
1 @>>> \pi_1 (\overline X_1)^{p',\ab} @>>> \Pi_1 @>>> G_{k_1} @>>> 1 \\
@.            @VVV           @VVV                @VVV    \\
1 @>>> \pi_1 (\overline X_2)^{p',\ab} @>>> \Pi_2 @>>> G_{k_2} @>>> 1
\endCD
$$
where 
$\Pi _i$ is the quotient 
$\pi_1 (X_i)^{(p')}/\Ker (\pi_1 (\overline X_i)^{p'}\twoheadrightarrow \pi_1 
(\overline X_i)^{p',\ab})$, and the maps $\Pi_i \to G_{k_i}$ 
are the natural projections; $i=1,2$. 
Further, the vertical maps are surjective. In particular, 
the representation $G_{k_1}\to G_{k_2} \to \Aut (\pi_1 (\overline X_2)^{p',\ab})$,
where $ G_{k_2} \to \Aut (\pi_1 (\overline X_2)^{p',\ab})$ is the natural 
representation, and $G_{k_1}\to G_{k_2}$ is the right vertical 
map in the above diagram,
is a quotient representation of the natural representation 
$G_{k_1}\to \Aut (\pi_1 (\overline X_1)^{p',\ab})$.
For $i\in \{1,2\}$, let $E_i$ be the set of eigenvalues,
counted with multiplicities, of the Frobenius element 
$\varphi _{k_i}$ acting on  $\pi_1 (\overline X_i)
^{p',\ab}$. Then $E_2\subset E_1$, since the map 
$G_{k_1}\to G_{k_2}$ maps $\varphi _{k_1}$ to  $\varphi _{k_2}$ (cf. Lemma 2.8(ii)). 
We will show that $E_1=E_2$. 

For an integer $n\ge 1$, let $k_{i,n}$ be the unique extension 
of $k_i$ of degree $n$; $i=1,2$. Then, by the Lefschetz trace formula, 
$\sharp X_i(k_{i,n})=1-\sum _{\alpha_i\in E_i} \alpha _{i}^n+q^n$, 
where $q\defeq \sharp (k_i)$ (cf. Lemma 2.8(ii) for the equality 
$\sharp (k_1)=\sharp (k_2)$). 
Recall that the map 
$\phi :\Sigma _{X_1}\to \Sigma _{X_2}$ is surjective (cf. Lemma 2.9), 
and if $x_2=\phi (x_1)$, then the residue fields $k(x_1)$ and $k(x_2)$ have the same 
cardinality (cf. Proposition 2.1(iii)). In particular, $\sharp (X_1(k_{1,n}))\ge 
\sharp (X_2(k_{2,n}))$ for all $n$. Thus, $\sum _{j=1}^r \beta _j^n\le 0$ 
for any $n\ge 1$, where $E\defeq E_1\smallsetminus E_2\defeq \{\beta_1,...,\beta_r\}$ 
($r=2g_1-2g_2\geq 0$). 
Write $\beta _j=\rho _je^{i\theta_j}$ 
($\rho_j\in\Bbb R_{>0}$, $\theta_j\in[0,2\pi)$), for $j\in \{1,...,r\}$ (note that
$\rho _j=q^{1/2}$ by the Riemann hypothesis for curves). 
Let $\Cal T$ be the set consisting of the $4$ quadrants of 
$\Bbb C=\Bbb R^2$. More precisely, $\Cal T=\{T_k\mid k\in\{1,2,3,4\}\}$, where 
$T_k\defeq\{\rho e^{i\theta}\mid \rho\in\Bbb R_{>0}, \theta\in[\frac{(k-1)\pi}{2},\frac{k\pi}{2})\}$. 
Thus, each $\alpha\in\Bbb C^{\times}$ belongs to a unique element of $\Cal T$, which 
we shall denote by $T(\alpha)$. 
Consider the map $\mu :\Bbb N\to \Cal T^r$ that maps 
an integer $n$ to $\{T(\beta_j^n)\}_{j=1}^r$.
Then there must exist integers $m_1<m_2$ such that 
$\mu (m_1)=\mu (m_2)$, since $\sharp (\Cal T^r)=4^r$ is finite. This implies 
that $e^{im_1\theta_j}$, and $e^{im_2\theta_j}$, belong to the same quadrant of the
$\Bbb C=\Bbb R^2$ 
for all $j\in \{1,...,r\}$. In particular,
$\Re (\beta _j^n)=\rho_j^n\cos n\theta_j>0$, where $n\defeq m_2-m_1\ge 1$. 
Suppose that $r>0$, then this implies that
$\Re (\sum _{j=1}^r\beta _j^n)=\sum _{j=1}^r \Re \beta _j^n >0$, which 
contradicts the above fact that $\sum _{j=1}^r \beta _j^n\le 0$, 
for all $n$. Thus, $r=0$, i.e., $E= E_1\smallsetminus E_2$ must be empty, and $E_1=E_2$. 

In particular, the $\hat \Bbb Z^{p'}$-ranks of  
$\pi_1 (\overline X_i)^{p',\ab}$, which equal to 
$2g_i$, are equal; $i=1,2$. This completes the proof of (i). 

Finally, we can conclude that 
$\phi$ is injective. For otherwise, there would exist an integer $n\ge 1$ 
such that $\sharp (X_1(k_n))>\sharp(X_2(k_n))$. Hence, $E\neq \emptyset $, 
which is a contradiction. This completes the proof of (ii). 
\qed
\enddemo

\proclaim
{Lemma 3.9} Assume that $\sigma: \frak G_1\to\frak G_2$ 
is strictly rigid. Then $\sigma^{(p')}:\frak G_1^{(p')}\to\frak G_2^{(p')}$ is 
an isomorphism. 
(In particular, in the prime-to-characteristic case, 
$\sigma$ is an isomorphism.) 
\endproclaim

\demo
{Proof} By Lemma 3.8(ii), the map $\phi :\Sigma _{X_1}\to \Sigma _{X_2}$ 
induced by $\sigma$ is bijective. For a finite subset $S_2$ of $\Sigma_{X_2}$ let
$S_1\defeq\phi^{-1} (S_2)$. Then $\sigma$ induces naturally a continuous, and 
surjective homomorphism $\tau' _{S_1,S_2}:\pi_1(U_1)^{(p')}
\twoheadrightarrow \pi_1(U_2)^{(p')}$, where  $\pi_1(U_i)^{(p')}$ is the 
maximal geometrically prime-to-$p$ quotient of the fundamental 
group $\pi_1(U_i)$ of $U_i\defeq X_i-S_i$; $i=1,2$. Further, we have the 
following commutative diagram:
$$
\CD
1    @>>>  \pi_1 (\overline U_1)^{p'}  @>>>  \pi_1 (U_1)^{(p')} @>>> G_{k_1} @>>> 1 \\
@.        @VVV         @V{\tau' _{S_1,S_2}}VV       @VVV    \\
1    @>>>  \pi_1 (\overline U_2)^{p'}  @>>>  \pi_1 (U_2)^{(p')} @>>> G_{k_2} @>>> 1
\endCD  
$$
The surjective homomorphism $\pi_1 (\overline U_1)^{p'}\twoheadrightarrow \pi_1 
(\overline U_2)^{p'}$ must be an isomorphism by [FJ], Proposition 15.4, 
since $X_i-S_i$ have the same topological 
type $(g_i,\sharp (\overline S_i))$, where $\overline S_i$ denotes the 
inverse image of $S_i$ in $\Sigma_{\overline X_i}$ ; $i=1,2$, by Lemma 3.8. 
(For the bijectivity $\bar S_1\isom\bar S_2$, apply Lemma 3.8(ii) to various open 
subgroups of $\frak G_1,\frak G_2$ corresponding to constant field extensions.) 
Thus, the map $\tau' _{S_1,S_2}$ is an isomorphism
(note that the surjective map $G_{k_1}\to G_{k_2}$ is an isomorphism).
Also, $\frak G_i^{(p')}=\underset {S_i}\to {\varprojlim}\  
\pi_1(X_i-S_i)^{(p')}$, where the projective 
limit is taken over all finite subsets $S_i$ of $\Sigma_{X_i}$; $i=1,2$. 
Further, $\sigma^{(p')}=\underset {\{S_1,S_2\}}\to {\varprojlim}
\tau' _{S_1,S_2}$, where the projective limit is taken over all finite 
subsets $S_1$ and $S_2$, corresponding to each other via $\phi$. Thus, 
$\sigma^{(p')}$ must be an isomorphism. 
\qed
\enddemo

Now, return to the general case. As above, 
let $\frak H_1\subset \frak G_1$,  $\frak H_2\subset \frak G_2$ be 
open subgroups such that $\sigma(\frak H_1)\subset\frak H_2$ 
and that the map $\sigma_{\frak H_1,\frak H_2}: 
\frak H_1\to\frak H_2$ obtained by restricting $\sigma: \frak G_1\to\frak G_2$, 
is strictly rigid 
with respect to $\tilde\phi: \Sigma_{\tilde X_1}\to\Sigma_{\tilde X_2}$. 
By Remark 3.2(iv), $\tilde\phi$ is Galois-equivariant with respect 
to $\sigma: \frak G_1\to\frak G_2$ (i.e., not only with respect to 
$\sigma_{\frak H_1,\frak H_2}: \frak H_1\to\frak H_2$), and, for each $\tilde x_1\in
\Sigma_{\tilde X_1}$, we have 
$\sigma (\frak D_{\tilde x_1})\underset \open \to \subset \frak D_{\tilde \phi(\tilde x_1)}$ 
and $\sigma (\frak D_{\tilde x_1}\cap\frak H_1)= \frak D_{\tilde \phi(\tilde x_1)}\cap\frak H_2$. 
Moreover, by Lemma 3.9, the map $(\sigma_{\frak H_1,\frak H_2})^{(p')}: 
\frak H_1^{(p')}\to\frak H_2^{(p')}$ induced by 
$\sigma_{\frak H_1,\frak H_2}$ is an isomorphism. 

Now, let us denote the finite separable extension of $K_i$ corresponding 
to $\frak H_i\subset\frak G_i$ by $K_{i,\frak H_i}$, and 
the (infinite) Galois extension of $K_{i,\frak H_i}$ corresponding to 
$\frak H_i\twoheadrightarrow \frak H_i^{(p')}$ by 
$\tilde K_{i,\frak H_i}^{(p')}$. 
Then, by applying the Isom-form proved in [ST], 
$(\sigma_{\frak H_1,\frak H_2})^{(p')}: \frak H_1^{(p')}\isom\frak H_2^{(p')}$ 
arises from a unique field isomorphism 
$\gamma_{\frak H_1^{(p')}, \frak H_2^{(p')}}: 
\tilde K_{2,\frak H_2}^{(p')}\isom \tilde K_{1,\frak H_1}^{(p')}$ 
that induces an isomorphism 
$K_{2, \frak H_2}\isom K_{1,\frak H_1}$. 

\proclaim{Lemma 3.10} Let $\frak H_i'\subset \frak H_i$, $i=1,2$ be open subgroups, 
such that $\sigma(\frak H_1')\subset \frak H_2'$ and that 
$\sigma_{\frak H_1',\frak H_2'}: \frak H_1'\to\frak H_2'$ is strictly rigid. 
Then, the field isomorphism 
$\gamma_{(\frak H_1')^{(p')}, (\frak H_2')^{(p')}}: 
\tilde K_{2,\frak H_2'}^{(p')}\isom \tilde K_{1,\frak H_1'}^{(p')}$ 
restricts to 
$\gamma_{\frak H_1^{(p')}, \frak H_2^{(p')}}: 
\tilde K_{2,\frak H_2}^{(p')}\isom \tilde K_{1,\frak H_1}^{(p')}$. 
\endproclaim

\demo{Proof} 
This follows just formally from the statement of the Isom-form proved in [ST], as follows, 
without recalling any construction in [ST]. 

Take an open subgroup $\frak H_2''$ of $\frak H_2'$ that is normal in $\frak H_2$. 
Then, by Lemma 3.7, we have 
$$\frak H_1''\defeq\sigma^{-1}(\frak H_2'')\subset
\sigma^{-1}(\frak H_2')=\frak H_1',$$
hence, by Remark 3.2(ii), $\frak H_1''\overset{\sigma}\to{\to} \frak H_2''$ is 
strictly rigid. 
Assume that $\gamma_{(\frak H_1'')^{(p')}, (\frak H_2'')^{(p')}}: 
\tilde K_{2,\frak H_2''}^{(p')}\isom \tilde K_{1,\frak H_1''}^{(p')}$ 
restricts to 
$\gamma_{(\frak H_1')^{(p')}, (\frak H_2')^{(p')}}: 
\tilde K_{2,\frak H_2'}^{(p')}\isom \tilde K_{1,\frak H_1'}^{(p')}$ 
and to
$\gamma_{\frak H_1^{(p')}, \frak H_2^{(p')}}: 
\tilde K_{2,\frak H_2}^{(p')}\isom \tilde K_{1,\frak H_1}^{(p')}$. 
Then $\gamma_{(\frak H_1')^{(p')}, (\frak H_2')^{(p')}}: 
\tilde K_{2,\frak H_2'}^{(p')}\isom \tilde K_{1,\frak H_1'}^{(p')}$ 
restricts to 
$\gamma_{\frak H_1^{(p')}, \frak H_2^{(p')}}: 
\tilde K_{2,\frak H_2}^{(p')}\isom \tilde K_{1,\frak H_1}^{(p')}$, as desired. 
So, it suffices to prove the desired property in the case where 
$\frak H_i'\subset\frak H_i$ is normal for $i=1,2$, and $\sigma$ 
induces naturally an isomorphism $\frak H_1/\frak H_1'\isom \frak H_2/\frak H_2'$ 
between finite groups. 

For $i=1,2$, let $\frak J_i'$ be the image of $\frak H_i'$ in $\frak H_i^{(p')}$, 
which is an open normal subgroup of $\frak H_i^{(p')}$. 
Let 
$\frak J_i\subset \frak H_i$ 
be the inverse image of $\frak J_i'$ in $\frak H_i$. Thus, 
we have the natural identification 
$\frak J_i^{(p')}=\frak J_i'$ and the following commutative diagram: 
$$
\matrix
\frak H_i' &\subset& \frak J_i &\subset& \frak H_i &\subset& \frak G_i \\
&&&&&&\\
\downarrow &&\downarrow &&\downarrow &&\downarrow \\
&&&&&&\\
(\frak H_i')^{(p')} &\twoheadrightarrow & \frak J_i^{(p')} 
&\hookrightarrow& \frak H_i^{(p')} &\to& \phantom{,}\frak G_i^{(p')}, 
\endmatrix
$$
in which the vertical arrows are natural surjective maps. 

Now, since the isomorphism $(\sigma_{\frak H_1',\frak H_2'})^{(p')}:
(\frak H_1')^{(p')}\isom(\frak H_2')^{(p')}$ 
is compatible with the natural (conjugate) actions of $\frak H_1$ and $\frak H_2$ 
with respect to $\frak H_1\overset{\sigma}\to{\twoheadrightarrow} \frak H_2$, 
the corresponding field isomorphism 
$\gamma_{(\frak H_1')^{(p')}, (\frak H_2')^{(p')}}: 
\tilde K_{2,\frak H_2'}^{(p')}\isom \tilde K_{1,\frak H_1'}^{(p')}$ 
is also compatible with the natural actions of $\frak H_1$ and $\frak H_2$ 
with respect to $\frak H_1\overset{\sigma}\to{\twoheadrightarrow} \frak H_2$. 
In particular, 
$\gamma_{(\frak H_1')^{(p')}, (\frak H_2')^{(p')}}$ 
restricts to $K_{2,\frak H_2}\isom K_{1,\frak H_1}$, hence induces 
an isomorphism $\alpha: K_{2,\frak H_2}^{(p')}\isom K_{1,\frak H_1}^{(p')}$ 
which is compatible with $\sigma: \frak H_1\twoheadrightarrow \frak H_2$, 
hence with $\sigma_{\frak H_1^{(p')},\frak H_2^{(p')}}: \frak H_1^{(p')}
\isom \frak H_2^{(p')}$. 
On the other hand, 
the isomorphism $\gamma_{\frak H_1^{(p')}, \frak H_2^{(p')}}: 
\tilde K_{2,\frak H_2}^{(p')}\isom \tilde K_{1,\frak H_1}^{(p')}$ 
is also compatible with $\sigma_{\frak H_1^{(p')},\frak H_2^{(p')}}: \frak H_1^{(p')}
\isom \frak H_2^{(p')}$. Thus, we conclude the desired identity 
$\alpha=\sigma_{\frak H_1^{(p')},\frak H_2^{(p')}}$ by the uniqueness of 
such a Galois-compatible isomorphism. (Indeed, this is included in the statement 
of the Isom-form proved in [ST].)  
\qed\enddemo

Now, consider the set $\Cal S$ ($\subset\Sub(\frak G_1)\times\Sub(\frak G_2)$) of 
all pairs of open subgroups $\frak H_1\subset\frak G_1$, 
$\frak H_2\subset\frak G_2$ such that $\sigma(\frak H_1)\subset \frak H_2$, 
that $\frak H_1\overset{\sigma}\to{\to}\frak H_2$ is strictly rigid, and 
that $\frak H_2$ is normal in $\frak G_2$. Then, as in the proof of Lemma 3.10, 
it follows from Lemma 3.7 and Remark 3.2(ii) that 
$(\frak H_1,\frak H_2)\in\Cal S$ implies 
that $\sigma(\frak H_1)=\frak H_2$, that 
$\frak H_1=\sigma^{-1}(\frak H_2)$
and that the image of $\Cal S$ in $\Sub(\frak G_2)$ is cofinal in the set of 
open subgroups of $\frak G_2$. 

For each pair $(\frak H_1,\frak H_2)\in\Cal S$, we get 
an isomorphism $\sigma_{\frak H_1^{(p')}, \frak H_2^{(p')}}: 
\frak H_1^{(p')}\isom\frak H_2^{(p')}$ 
by Lemma 3.9, which is Galois-compatible with respect to 
$\sigma:\frak G_1\to\frak G_2$. By the Isom-form proved in [ST], 
$\sigma_{\frak H_1^{(p')}, \frak H_2^{(p')}}$ 
induces an isomorphism $\gamma_{\frak H_1^{(p')}, \frak H_2^{(p')}}: 
\tilde K_{2,\frak H_2}^{(p')}\isom\tilde K_{1,\frak H_1}^{(p')}$, 
which is Galois-compatible with respect to 
$\sigma:\frak G_1\to\frak G_2$. By Lemma 3.10, 
$\gamma_{\frak H_1^{(p')}, \frak H_2^{(p')}}$ can be patched together 
and define an isomorphism $\tilde \gamma: \tilde K_{2}\isom(\tilde K_1)^{\frak N}$, 
where $\frak N\defeq\Ker(\sigma: \frak G_1\to\frak G_2)$, which is Galois-compatible 
with respect to $\sigma:\frak G_1\to\frak G_2$. 

In the profinite (resp. prime-to-characteristic) case, 
$\tilde K_{2}$ admits no nontrivial 
separable (resp. geometrically prime-to-$p$) extension, hence 
so is $(\tilde K_1)^{\frak N}\ (\simeq \tilde K_2)$. 
This implies that $(\tilde K_1)^{\frak N}=\tilde K_1$, i.e., $\frak N=\{1\}$. 
Thus, we obtain $\tilde\gamma: \tilde K_{2}\isom \tilde K_1$, 
which is Galois-compatible with respect to $\sigma: \frak G_1\to\frak G_2$, 
as desired. Finally, the uniqueness of such $\tilde\gamma$ follows 
formally from the uniqueness in the statement of the Isom-form (proved in [Uchida2][ST]). 
This finishes the proof of Theorem 3.4. \qed

\definition{Remark 3.11} 
We have proved Theorem 3.4 by reducing it to the {\it statement} of the Isom-form, by means of 
Lemma 3.9. Instead, we could mimic the {\it proof} of the Isom-form. 
\enddefinition

\subhead 
\S 4. Proper homomorphisms between Galois groups
\endsubhead

In this section we shall investigate a class of homomorphisms 
between (geometrically prime-to-characteristic quotients of) 
absolute Galois groups of function fields of curves over finite 
fields, which we call proper. We follow the notations in \S\S1--3. 
In particular, we follow the Notation at the beginning of subsection 2.2. 
We assume that Condition 3 holds. 

\comment
For $i\in \{1,2\}$, let $k_i$ be a finite field 
of characteristic 
$p_i> 0$. Let $X_i$ be a smooth, proper, geometrically connected 
curve of genus $g_i$ over $k_i$. Let $K_i=K_{X_i}$ be the function field of $X_i$ 
and fix 
an algebraic closure $\overline K_i$ of $K_i$. Let 
$K_i^{\sep}$ be the separable closure of $K_i$ in $\overline K_i$, 
and $\bar k_i$ the algebraic closure of $k_i$ in $\overline K_i$. 
Following the notations 
in $\S 1$, we will write $G_i\defeq G_{K_i}= \Gal (K_i^{\sep}/K_i)$
for the absolute Galois group of $K_i$, and 
$\bar{G}_i\defeq G_{K_i\bar k_i}=\Gal (K_i^{\sep}/K_i\bar k_i)$
for the absolute Galois group of $K_i\bar k_i$. 

Let $N_i$ be a normal closed subgroup of $G_i$ and set 
$\frak G_i\defeq G_i/N_i$. Let $\tilde K_i$ denote the Galois extension of $K_i$ 
corresponding to $N_i$, i.e., $\tilde K_i\defeq (K_i^{\sep})^{N_i}$. 
Let $\bar \frak G_i$ be the image of $\bar G_i$ in $\frak G_i$, and set 
$\frak G_{k_i}\defeq \frak G_i/\bar \frak G_i$, which is a quotient of 
$G_{k_i}$. For $i=1,2$, let us denote by $\varphi_{k_i}$ the image 
in $\frak G_{k_i}$ of the $\sharp(k_i)$-th power Frobenius element of $G_{k_i}$. 
Write $\tilde X_i$ for the integral closure of $X_i$ in $\tilde K_i$. 
The Galois group $\frak G_i$ acts naturally on the set 
$\Sigma _{\tilde X_i}$, and the quotient  $\Sigma _{\tilde X_i}/\frak G_i$ 
is naturally identified with $\Sigma _{X_i}$. 
For a point $\tilde x_i\in \Sigma _{\tilde X_i}$, 
let $\frak G_i\supset\frak D_{\tilde x_i}\supset\frak I_{\tilde x_i}$ be the 
decomposition and the inertia subgroups, respectively, and set 
$\frak G_{k_i(x_i)}\defeq\frak D_{\tilde x_i}/\frak I_{\tilde x_i}$. 


Now, let 
$$\sigma:\frak G_1\to \frak G_2$$ 
be 
an open 
homomorphism between profinite groups. 
\endcomment

\definition 
{Definition 4.1} (Well-Behaved Homomorphisms) We say that 
$\sigma :\frak G_1\to \frak G_2$ is well-behaved, 
if there exists a map
$\tilde \phi :\Sigma_ {\tilde X_1} \to \Sigma _{\tilde X_2}$, 
such that 
$\sigma (\frak D_{\tilde x_1})\underset \open \to \subset \frak D_{\tilde \phi(\tilde x_1)}$
for each $\tilde x_1\in \Sigma_ {\tilde X_1}$. 
\enddefinition

\definition
{Remark 4.2} {\rm (i)} Given a commutative diagram of maps between 
profinite groups:
$$
\CD
\frak G_1 @>>> \frak G_2 \\
@VVV     @VVV \\
\frak G_1' @>>> \frak G_2' \\
\endCD
$$
where the vertical arrows are surjective, and the map $\frak G_1\to \frak G_2$ is
well-behaved, then the map $\frak G_1'\to \frak G_2'$ is well-behaved.

\noindent
{\rm (ii)} 
Let $\frak H_1\subset\frak G_1$, $\frak H_2\subset\frak G_2$ be open subgroups 
such that $\sigma(\frak H_1)\subset \frak H_2$. 
Then if $\sigma:\frak G_1\to \frak G_2$ is well-behaved, 
then the natural homomorphism $\frak H_1\to \frak H_2$ induced by $\sigma$ is 
well-behaved. (Here, $\frak H_i$ is considered as a quotient of the absolute 
Galois group that is the inverse image in $G_i$ of $\frak H_i\subset\frak G_i$.) 

\noindent
{\rm (iii)} If $\sigma:\frak G_1\to\frak G_2$ is strictly rigid (cf. Definition 3.1), 
then it is well-behaved. 

\noindent
(iv) As in Proposition 2.2, 
let $l$ be a prime number $\neq p_1,p_2$, and put 
the following two assumptions: 
{\rm (1)} $N_2^l=N_2$, or, equivalently, $\tilde K_2$ admits 
no $l$-cyclic extension; and 
{\rm (2)} $\tilde K_2$ contains a primitive $l$-th roots of unity. 
Then, first, if $\sigma: \frak G_1\to\frak G_2$ is rigid, then it is well-behaved by 
Remark 3.2(iv). 
Second, if $\sigma: \frak G_1\to\frak G_2$ is well-behaved with respect to 
$\tilde \phi: \Sigma_{\tilde X_1} \to \Sigma_{\tilde X_2}$, 
then we must have $\Sigma_{\tilde X_1,\sigma,l}=\Sigma_{\tilde X_1}$ 
and $\tilde\phi=\tilde\phi_{\sigma,l}$. In particular, 
then $\tilde\phi$ is unique and Galois-equivariant with respect to $\sigma$, 
hence induces naturally a map
$\phi\  (=\phi_{\sigma,l}):\Sigma _{X_1} \to \Sigma_ {X_2}$. 
\enddefinition

\definition
{Definition 4.3} (Proper Homomorphisms) 
We say that 
$\sigma :\frak G_1\to \frak G_2$ is proper, 
if $\sigma$ is well-behaved with respect to 
$\tilde\phi:\Sigma_ {\tilde X_1} \to \Sigma _{\tilde X_2}$, 
such that $\tilde\phi$ is Galois-equivariant with respect to $\sigma$, 
and the map 
$\phi :\Sigma _{X_1} \to \Sigma _ {X_2}$ induced by $\tilde\phi$ 
has finite fibers, 
i.e., for each $x_2\in \Sigma _ {X_2}$, the fiber 
$\phi ^{-1}(x_2)$ is a (possibly empty) finite set. 
\enddefinition

\definition
{Remark 4.4} {\rm (i)} Given a commutative diagram of maps between 
profinite groups:
$$
\CD
\frak G_1 @>>> \frak G_2 \\
@VVV     @VVV \\
\frak G_1' @>>> \frak G_2' \\
\endCD
$$
where the vertical arrows are surjective, and the map $\frak G_1\to \frak G_2$ is
proper, then the map $\frak G_1'\to \frak G_2'$ is proper.

\noindent
{\rm (ii)} 
Let $\frak H_1\subset\frak G_1$, $\frak H_2\subset\frak G_2$ be open subgroups 
such that $\sigma(\frak H_1)\subset \frak H_2$. 
Then if $\sigma:\frak G_1\to \frak G_2$ is proper, 
then the natural homomorphism $\frak H_1\to \frak H_2$ induced by $\sigma$ is 
proper. (Here, $\frak H_i$ is considered as a quotient of the absolute 
Galois group that is the inverse image in $G_i$ of $\frak H_i\subset\frak G_i$.)  
\enddefinition

In the rest of this section, we 
assume that 
Condition 1 holds. 
Assume, moreover, that the continuous open homomorphism 
$\sigma:\frak G_1\to\frak G_2$ is well-behaved with respect to 
$\tilde\phi: \Sigma_ {\tilde X_1} \to \Sigma _{\tilde X_2}$. 
By Lemma 2.4, we have $p\defeq p_1=p_2$. 
Let $\tilde x_1\in\Sigma_{\tilde X_1}$ and set 
$\tilde x_2\defeq\tilde\phi(\tilde x_1)$. 
Denote by $x_1$ (resp. $x_2$) the image of $\tilde x_1$ (resp. $\tilde x_2$) 
in $\Sigma_{X_1}$ (resp. $\Sigma_{X_2}$). 
Then we have 
$$\frak D_{\tilde x_1}\overset\sigma\to\twoheadrightarrow \sigma(\frak D_{\tilde x_1})
\underset \open \to \subset \frak D_{\tilde x_2}.$$ 
By this and Proposition 2.1(v), we have 
$$\frak I_{\tilde x_1}\overset\sigma\to\twoheadrightarrow \sigma(\frak I_{\tilde x_1})
\underset \open \to \subset \frak I_{\tilde x_2}.$$ 

In particular, $\sigma$ induces an open injective homomorphism 
$\tau^{\t}_{\tilde x_1}: \frak I_{\tilde x_1}^{\t}\hookrightarrow \frak I_{\tilde x_2}^{\t}$, 
where $\frak I_{\tilde x_1}^{\t}$ (resp. $\frak I_{\tilde x_2}^{\t}$) 
denotes the inertia subgroup of $\frak D_{\tilde x_1}^{\t}$ 
(resp. of $\frak D_{\tilde x_2}^{\t}$). Note that we have natural 
identifications 
$M_1 \isom M_{k (x_1)^{\sep}} \isom  \frak I_{\tilde x_1}^{\t}$ 
and 
$M_2 \isom M_{k (x_2)^{\sep}} \isom \frak I_{\tilde x_2}^{\t}$, 
where $M_i\defeq M_{K_i^{\sep}}$ is the (global) module of roots of unity for $i=1,2$. 

Now, we introduce the following important concept of 
rigidity of inertia. 

\definition
{Definition 4.5} (Inertia-Rigid Homomorphisms) 
We say that the well-behaved homomorphism $\sigma :\frak G_1\to \frak G_2$ is inertia-rigid, 
if there exists an isomorphism $\tau :M_1\isom M_2$ 
of $\hat \Bbb Z^{p'}$-modules, such that for each $\tilde x_1\in \Sigma _{\tilde X_1}$, 
there exists a positive integer $e_{\tilde x_1}$ such that 
the following diagram commutes:
$$
\matrix
M_1 &\isom &M_{k (x_1)^{\sep}} &\isom  &\frak I_{\tilde x_1}^{\t}\\
&&&&\\
e_{\tilde x_1}\cdot\tau\downarrow\phantom{e_{\tilde x_1}\cdot\tau}&&&& 
\phantom{\tau_{\tilde x_1}^{\t}}\downarrow\tau_{\tilde x_1}^{\t}   \\
&&&&\\
M_2 &\isom &M_{k (x_2)^{\sep}} &\isom &\frak I_{\tilde x_2}^{\t}\\
\endmatrix\tag {$4.1$}
$$
where $\tilde x_2\defeq\tilde\phi(\tilde x_1)$; 
$x_1$ (resp. $x_2$) is the image of $\tilde x_1$ (resp. $\tilde x_2$) 
in $\Sigma_{X_1}$ (resp. $\Sigma_{X_2}$); and 
the isomorphisms are the canonical identifications. 
\enddefinition

\definition
{Remark 4.6} {\rm (i)} Given a commutative diagram of maps between 
profinite groups:
$$
\CD
G_1 @>>> G_2 \\
@VVV     @VVV \\
G_1^{(p')} @>>> G_2^{(p')} \\
\endCD
$$
where the vertical arrows are natural surjective maps, 
and the map $G_1\to G_2$ is
inertia-rigid, then the map $G_1^{(p')}\to G_2^{(p')}$ is inertia-rigid.

\noindent
{\rm (ii)} 
Let $\frak H_1\subset\frak G_1$, $\frak H_2\subset\frak G_2$ be open subgroups 
such that $\sigma(\frak H_1)\subset \frak H_2$. 
Then if $\sigma:\frak G_1\to \frak G_2$ is inertia-rigid, 
then the natural homomorphism $\frak H_1\to \frak H_2$ induced by $\sigma$ is 
inertia-rigid. (Here, $\frak H_i$ is considered as a quotient of the absolute 
Galois group that is the inverse image in $G_i$ of $\frak H_i\subset\frak G_i$.)  
\enddefinition

\definition{Remark 4.7} (i) 
Set $\frak e_{\tilde x_1}\defeq[\frak I_{\tilde x_2}: \sigma(\frak I_{\tilde x_1})]$ 
and $\frak e_{\tilde x_1}^{\t}\defeq[\frak I_{\tilde x_2}^{\t}: \tau_{\tilde x_1}^{\t}(\frak I_{\tilde x_1}^{\t})]$. 
Note that $p\nmid \frak e_{\tilde x_1}^{\t}$ 
and that there exists an integer $b_{\tilde x_1}\geq 0$, such that 
$\frak e_{\tilde x_1}=p^{b_{\tilde x_1}}\frak e_{\tilde x_1}^{\t}$. 
(In the prime-to-characteristic case, we have $\frak e_{\tilde x_1}=\frak e_{\tilde x_1}^{\t}$ 
and $b_{\tilde x_1}= 0$.)
Now, in Definition 4.5, 
there must exist an integer $a_{\tilde x_1}\geq 0$, such that 
$e_{\tilde x_1}=p^{a_{\tilde x_1}}\frak e_{\tilde x_1}^{\t}$, or, equivalently,
$e_{\tilde x_1}=p^{c_{\tilde x_1}}\frak e_{\tilde x_1}$, where 
$c_{\tilde x_1}\defeq a_{\tilde x_1}-b_{\tilde x_1}\in\Bbb Z$. Moreover, set 
$a\defeq \min\{a_{\tilde x_1}\mid \tilde x_1\in\Sigma_{\tilde X_1}\}$. Then, 
replacing $\tau$ by $p^a\tau$ and $e_{\tilde x_1}$ by $p^{-a}e_{\tilde x_1}=
p^{a_{\tilde x_1}-a}\frak e_{\tilde x_1}^{\t}$, we may assume that $a=0$. 

Assume, moreover, that $\sigma$ is proper and that we are in the 
profinite case. Then, in fact, we have $c_{\tilde x_1}=0$ for every $\tilde x_1\in\tilde X_1$ 
eventually, if we choose $\tau$ with $a=0$. (This follows from Theorem 4.8 below and its proof.) 
Thus, in the profinite case, we may assume $e_{\tilde x_1}=\frak e_{\tilde x_1}$ 
in Definition 4.5 from the beginning. In the prime-to-characteristic case, however, it seems 
difficult to specify the value of $e_{\tilde x_1}$ a priori. (If we 
assumed $e_{\tilde x_1}=\frak e_{\tilde x_1}$ in the prime-to-characteristic 
case, then inertia-rigid homomorphisms would cover only tame homomorphisms 
$K_2\to K_1$.) 

\noindent
(ii) In the situation of Definition 4.5, we have 
$$\frak D_{\tilde x_1}\overset\sigma\to\twoheadrightarrow \frak E_{\tilde x_1}
\defeq\sigma(\frak D_{\tilde x_1})\subset\frak D_{\tilde x_2}.$$
The subgroup $\frak E_{\tilde x_1}\subset \frak D_{\tilde x_2}$ 
corresponds to a finite extension 
$L_{x_1}/(K_2)_{x_2}$ of the $x_2$-adic completion $(K_2)_{x_2}$ of $K_2$
Thus, the residue field $\ell _{x_1}$ of $L_{x_1}$ is a finite extension 
of the residue field $k (x_2)$ at $x_2$.
We have the following commutative diagram:
$$
\CD
\frak D_{\tilde x_1} @>>>  \frak E_{\tilde x_1} \\
@VVV   @VVV   \\
\frak D_{\tilde x_1}^{\t} @>>>  \frak E_{\tilde x_1}^{\t} \\
\endCD
$$
where the vertical maps are the canonical surjections onto the maximal 
tame quotients, and the horizontal maps are naturally induced by $\sigma$.
Further, the lower horizontal map, which is surjective, induces naturally
an isomorphism $\frak I_{\tilde x_1}^{\t}\isom \frak J_{\tilde x_1}^{\t}$
by Proposition 2.1(v). Here, $\frak I_{\tilde x_1}^{\t}$ (resp. $\frak J_{\tilde x_1}^{\t}$) 
denotes the inertia subgroup of $\frak D_{\tilde x_1}^{\t}$ 
(resp. of $\frak E_{\tilde x_1}^{\t}$). 
We have a natural identification 
$\frak J_{\tilde x_1}^{\t}\isom\frak I_{\tilde x_2}^{\t}$, 
where $\frak I_{\tilde x_2}^{\t}$ is the inertia subgroup of $\frak D_{\tilde x_2}^{\t}$ 
(obtained via the natural identifications 
$M_{(K_2)_{x_2}^{\sep}}\isom \frak I_{\tilde x_2}^{\t}$, 
$M_{L_{x_1}^{\sep}}\isom  \frak J_{\tilde x_1}^{\t}$, 
and $(K_2)_{x_2}^{\sep}=L_{x_1}^{\sep}$), 
which, composed with the natural map 
$\frak J_{\tilde x_1}^{\t}\to \frak I_{\tilde x_2}^{\t}$ induced by the inclusion 
$\frak E_{\tilde x_1}\to \frak D_{\tilde x_2}$, is the $\frak e_{\tilde x_1}$-th 
power map $\frak I_{\tilde x_2}^{\t}@>[\frak e_{\tilde x_1}]>>\frak I_{\tilde x_2}^{\t}$, 
as is easily verified. We define 
$$\tau _{\tilde x_1,\tilde x_2}^{\t}:\frak I^{\t}_{\tilde x_1}\isom \frak I_{\tilde x_2}^{\t}$$
to be the natural isomorphism obtained by composing the natural 
isomorphism $\frak I^{\t}_{x_1}\isom  \frak J_{\tilde x_1}^{\t}$ 
induced by $\sigma$ (cf. Proposition 2.1(v)), with the  
above natural identification $\frak J_{\tilde x_1}^{\t} \isom \frak I_{\tilde x_2}^{\t}$. 

Now, the inertia-rigidity is equivalent to requiring the commutativity of the 
following diagram: 
$$
\matrix
M_1 &\isom &M_{(K_1)_{x_1}^{\sep}} &\isom  &\frak I_{\tilde x_1}^{\t}\\
&&&&\\
p^{c_{\tilde x_1}}\cdot\tau\downarrow\phantom{p^{c_{\tilde x_1}}\cdot\tau}&&&& 
\phantom{\tau_{\tilde x_1}^{\t}}\downarrow\tau_{\tilde x_1,\tilde x_2}^{\t}\\
&&&&\\
M_2 &\isom &M_{(K_2)_{x_2}^{\sep}} &\isom &\frak I_{\tilde x_2}^{\t}
\endmatrix
$$
in which both vertical arrows are isomorphisms. 
\enddefinition

Define $\Hom(K_2,K_1)^{\sep}\subset\Hom(K_2,K_1)$ to be the set of 
separable homomorphisms $K_2\to K_1$. 
Define 
$\Hom(\frak G_1,\frak G_2)^{\pr,\inrig}\subset\Hom(\frak G_1,\frak G_2)$ to be 
the set of proper (hence continuous open), inertia-rigid homomorphisms 
$\frak G_1\to\frak G_2$. 
Our aim in this section is to prove the following. 

\proclaim
{Theorem 4.8} 
The natural map $\Hom(K_2,K_1)\to\Hom(\frak G_1, \frak G_2)/\Inn(\frak G_2)$ 
induces a bijection 
$$\Hom(K_2,K_1)^{\sep}\isom\Hom(\frak G_1, \frak G_2)^{\pr,\inrig}/\Inn(\frak G_2).$$ 
More precisely, 

\noindent
{\rm (i)} If $\gamma:K_2\to K_1$ is a separable homomorphism between fields, then 
the homomorphism $\frak G_1\to\frak G_2$ induced by $\gamma$ (up to inner automorphisms) 
is proper, inertia-rigid. 

\noindent
{\rm (ii)} If $\sigma:\frak G_1\to \frak G_2$ is a proper, inertia-rigid 
homomorphism between profinite groups, then there 
exists a unique homomorphism $\tilde \gamma: \tilde K_2\to
\tilde K_1$ of fields, such that $\tilde\gamma\circ \sigma (g_1)=g_1
\circ \tilde \gamma$, for all $g_1\in \frak G_1$, which induces a 
separable homomorphism $K_2\to K_1$. 
\endproclaim

\definition
{Remark 4.9} 
(i) Assume that $\sigma:\frak G_1\to\frak G_2$ is a rigid homomorphism. 
Then it follows from Lemma 3.8(ii) that $\sigma$ is proper. Further, 
$\sigma$ is inertia-rigid. This can be reduced to the 
case where $\sigma$ is strictly rigid, and then deduced from class field theory 
as in the arguments preceding Lemma 4.12. (Note that then $\phi$ is bijective 
by Lemma 3.8(ii).) Thus, Theorem 4.8 can be viewed as 
a generalization of Theorem 3.4. 

\noindent
(ii) The natural map $\Hom(K_2^{\perf},K_1^{\perf})\to\Hom(\frak G_1, \frak G_2)/\Inn(\frak G_2)$ 
induces a bijection 
$$\Hom(K_2^{\perf},K_1^{\perf})
/\Frob^{\Bbb Z}\ \isom\ \Hom(\frak G_1, \frak G_2)^{\pr,\inrig}/\Inn(\frak G_2).$$ 
Indeed, this follows from Theorem 4.8, since the natural map 
$\Hom(K_2,K_1)\to\Hom(K_2^{\perf},K_1^{\perf})
$ 
induces 
$$\Hom(K_2,K_1)^{\sep}\isom\Hom(K_2^{\perf},K_1^{\perf})/\Frob^{\Bbb Z}.$$
\enddefinition

The rest of this section is devoted to the proof of Theorem 4.8. 

First, to prove (i), let $\gamma:K_2\to K_1$ be a separable homomorphism. 
Then $\gamma$ induces naturally an open injective homomorphism 
$G_1\hookrightarrow G_2$ (up to $\Inn(G_2)$) and then an open homomorphism 
$\sigma: \frak G_1\to\frak G_2$ (up to $\Inn(\frak G_2)$). 
The map $\sigma$ is well-behaved with respect to 
the map $\phi: \Sigma_{X_1}\to \Sigma_{X_2}$ that arises from a finite 
morphism $X_1\to X_2$ of schemes corresponding to $\gamma:K_2\to K_1$. 
Thus, each fiber of $\phi$ is finite, hence $\sigma$ is proper. Next, if 
we define $\tau: M_1\isom M_2$ to be the identification 
$M_{K_1^{\sep}}\isom M_{K_2^{\sep}}$ (with respect to a suitable extension 
$K_2^{\sep}\isom K_1^{\sep}$ of $\gamma: K_2\to K_1$), then diagram (4.1) 
commutes with $e_{\tilde x_1}$ defined to be the ramification index of 
$K_1/K_2$ at $\tilde x_1$. Thus, $\sigma$ is inertia-rigid. 

Next, to prove (ii), let $\sigma:\frak G_1\to\frak G_2$ be a proper, inertia-rigid 
homomorphism. 

\proclaim
{Lemma 4.10} 
Condition 2 holds for $\sigma:\frak G_1\to\frak G_2$. 
\endproclaim

\demo
{Proof} 
Same as that of Lemma 3.6. 
\qed
\enddemo

Thus, we may apply Lemmas 2.6--2.9 to $\sigma$. 

Next, let $\tau:M_1\isom M_2$ be the isomorphism appearing in the 
definition of inertia-rigid homomorphism, so that diagram (4.1) 
commutes for each $\tilde x_1\in\Sigma_{\tilde X_1}$ and for some 
$e_{\tilde x_1}\in\Bbb Z_{>0}$. 

\proclaim
{Lemma 4.11} 
{\rm (i)} 
The isomorphism $\tau :M_1\isom M_2$ is 
Galois-equivariant with respect to $\sigma$. 

\noindent
{\rm (ii)} The positive integers $e_{\tilde x_1}$, $\frak e_{\tilde x_1}$ and 
$\frak e_{\tilde x_1}^{\t}$ depend only on the image $x_1\in\Sigma_{X_1}$ 
of $\tilde x_1$. 
\endproclaim

\demo
{Proof} (i) For each $\tilde x_1\in\Sigma_{\tilde X_1}$, the commutativity 
of diagram (4.1), together with Proposition 2.1(iv), implies 
that $\tau$ is Galois-equivariant with respect to 
$\frak D_{\tilde x_1}\overset\sigma\to{\to}\frak D_{\tilde\phi(\tilde x_1)}$. 
Our assertion then 
follows, since $\frak G_1$ is generated 
by the decomposition subgroups $\frak D_{\tilde x_1}$ for all 
$\tilde x_1\in \Sigma _{\tilde X_1}$, as follows from Chebotarev's density theorem.

\noindent
(ii) Take another $\tilde x_1'\in\Sigma_{\tilde X_1}$ above $x_1\in\Sigma_{X_1}$ 
and set $\tilde x_2'\defeq\tilde\phi(\tilde x_1')$. Fix $\gamma\in\frak G_1$ 
such that $\tilde x_1'=\gamma\tilde x_1$. By the Galois-equivariance property 
of $\tilde\phi$ 
(cf. Remark 4.2(iv)), we have then $\tilde x_2'=\sigma(\gamma)\tilde x_2$. 
Denote by $[\gamma]$ (resp. $[\sigma(\gamma)]$) the inner automorphism of 
$\frak G_1$ (resp. $\frak G_2$) induced by $\gamma$ (resp. $\sigma(\gamma)$). 
Then we have the following commutative diagram
$$
\CD
\frak I_{\tilde x_1}@>[\gamma]>>  \frak I_{\tilde x_1'} \\
@V{\sigma}VV   @V{\sigma}VV   \\
\frak I_{\tilde x_2}@>[\sigma(\gamma)]>>  \frak I_{\tilde x_2'} 
\endCD
$$
in which both rows are isomorphisms. From this, 
it follows that $\frak e_{\tilde x_1'}=\frak e_{\tilde x_1}$. 
Next, this commutative diagram induces the following commutative diagram 
$$
\CD
\frak I_{\tilde x_1}^{\t}@>[\gamma]>>  \frak I_{\tilde x_1'}^{\t} \\
@V{\tau_{\tilde x_1}^{\t}}VV   @V{\tau_{\tilde x_1'}^{\t}}VV   \\
\frak I_{\tilde x_2}^{\t}@>[\sigma(\gamma)]>>  \frak I_{\tilde x_2'}^{\t} 
\endCD
$$
in which both rows are isomorphisms. From this, 
it follows that $\frak e_{\tilde x_1'}^{\t}=\frak e_{\tilde x_1}^{\t}$. 
Finally, combined with (i), this last commutative diagram also implies that 
$ e_{\tilde x_1'}=e_{\tilde x_1}$. 
\qed
\enddemo

{}From now on, we shall write $e_{x_1}$, $\frak e_{x_1}$ and 
$\frak e_{x_1}^{\t}$ for $e_{\tilde x_1}$, $\frak e_{\tilde x_1}$ and 
$\frak e_{\tilde x_1}^{\t}$, respectively. Further, according to this, 
we shall write $a_{x_1}$, $b_{x_1}$ and $c_{x_1}$ 
for the invariants $a_{\tilde x_1}$, $b_{\tilde x_1}$ and $c_{\tilde x_1}$ in 
Remark 4.7(i), respectively. 
We may and shall also assume $a\ (=\min\{a_{x_1}\mid x_1\in\Sigma_{X_1}\})=0$ (cf. Remark 4.7(i)). 

We have the following commutative diagram of exact sequences:

$$
\CD
1@>>>k_1^{\times}@>>> 
{\underset {x_2\in \Sigma _{X_2}}\to \prod}
(\underset{x_1\in\phi^{-1}(x_2)}\to{\prod}k (x_1)^{\times}) 
@>>> {\frak G_1}^{(p'),\ab}\\
@.  @VVV            @VVV           @VVV      \\
1@>>>k_2^{\times}@>>> {\underset {x_2\in \Sigma _{X_2}}\to \prod} 
k(x_2)^{\times} @>>> {\frak G_2}^{(p'),\ab}
\endCD 
$$
{}from global class field theory. 
Here, the map ${\frak G_1}^{(p'),\ab}\twoheadrightarrow  
{\frak G_2}^{(p'),\ab}$ is naturally induced by $\sigma$. The right horizontal 
maps are induced by Artin's reciprocity map, and the map 
${\underset {x_2\in \Sigma _{X_2}}\to \prod}(\underset{x_1\in\phi^{-1}(x_2)}\to{\prod}
k (x_1)^{\times})\to {\underset {x_2\in \Sigma _{X_2}}\to \prod} k 
(x_2)^{\times}$ maps each component 
$
k (x_1)^{\times}$ to $k (x_2)^{\times}$ as follows. First, 
$k (x_1)^{\times}$ maps isomorphically onto 
$\ell_{x_1}^{\times}$ via the natural 
identification induced by $\sigma$ (cf. Remark 4.7(ii) and Proposition 2.1(iii)). 
Then 
$\ell_{x_1}^{\times}$ maps to $k(x_2)^{\times}$  
by the $\frak e_{x_1}$-th power of the norm map. 

The above diagram induces, for each $x_2\in \Sigma _{X_2}$, the following 
commutative diagram: 
$$
\CD
k_1^{\times}@>>> \underset{x_1\in\phi^{-1}(x_2)}\to{\prod} k (x_1)^{\times} \isom 
\underset{x_1\in\phi^{-1}(x_2)}\to{\prod} \ell_{x_1} ^{\times}\\
@VVV            @VVV \\
k_2^{\times}@>>> k (x_2)^{\times}\\
\endCD 
$$
where the map $k_2^{\times}\to k(x_2)^{\times}$ is the natural 
embedding, the map $k_1^{\times}\to 
\underset{x_1\in\phi^{-1}(x_2)}\to{\prod}
k (x_1)^{\times}$ is the natural diagonal embedding, 
the isomorphism $\underset{x_1\in\phi^{-1}(x_2)}\to{\prod} k (x_1)^{\times}\isom \
\underset{x_1\in\phi^{-1}(x_2)}\to{\prod} \ell_{x_1} ^{\times}$, and the map 
$\underset{x_1\in\phi^{-1}(x_2)}\to{\prod}\ell_{x_1} ^{\times}\to k (x_2)^{\times}$ 
are as above. By passing to various open subgroups corresponding to 
extensions of the constant fields, and to the projective limit 
via the norm maps, we obtain the following commutative diagram:
$$
\CD
M_{k_1^{\sep}} @>>> \oplus_{\bar x_1\in\bar\phi^{-1}(\bar x_2)} M_{k (x_1)^{\sep}}
\overset {\oplus_{\bar x_1\in\bar\phi^{-1}(\bar x_2)}\rho_{x_1}}\to \isom  
\oplus_{\bar x_1\in\bar\phi^{-1}(\bar x_2)} M_{\ell _{x_1}^{\sep}} \\ 
@VVV   @VVV \\
M_{k_2^{\sep}} @>>> M_{k (x_2)^{\sep}}\\
\endCD
$$
where $\rho _{x_1}: M_{k (x_1)^{\sep}}\isom  M_{\ell _{x_1}^{\sep}}$
is the natural isomorphism induced by $\sigma$ (cf. 
Remark 4.7(ii) and Proposition 2.1(v)). 
Here, $\bar x_2\in \Sigma_{\overline X_2}$ is any point above $x_2$ and 
$\bar\phi: \Sigma_{\overline X_1}\to\Sigma_{\overline X_2}$ is 
obtained as the inductive limit 
of $\phi$'s for various open subgroups corresponding to extensions 
of the constant fields. Observe that 
$\bar\phi: \Sigma_{\overline X_1}\to\Sigma_{\overline X_2}$ 
has finite fibers, since 
(i) $\phi:\Sigma_{X_1}\to\Sigma_{X_2}$ 
has finite fibers, 
(ii) the projection $\Sigma_{\overline X_1} \to\Sigma_{X_1}$ 
has finite fibers, and 
(iii) $\bar\phi$ is compatible with $\phi$. 

This can be rewritten as:
$$
\CD
M_1 @>>> 
\oplus_{\bar x_1\in\bar\phi^{-1}(\bar x_2)} M_1
\overset{\oplus_{\bar x_1\in\bar\phi^{-1}(\bar x_2)}\rho_{x_1}}\to
\isom 
\oplus_{\bar x_1\in\bar\phi^{-1}(\bar x_2)} M_2\\ 
@VVV   @VVV \\
M_2 @=  M_2 \\
\endCD \tag{$4.2$}
$$
via the natural identifications $M_{k (x_1)^{\sep}}\isom M_1$
(resp. $M_{\ell _{x_1}^{\sep}}\isom M_2$) for 
$x_1\in \phi^{-1}(x_2)$; $M_{k (x_2)^{\sep}}\isom M_2$; and 
$M_{k_i^{\sep}}\isom M_i$, $i=1,2$. Thus, in diagram (4.2)
the map $M_1\to \oplus _{\bar x_1\in\bar\phi^{-1}(\bar x_2)}M_1$ 
is the natural diagonal embedding, and the map 
$\oplus _{\bar x_1\in\bar\phi^{-1}(\bar x_2)} M_2\to M_2$
is the map $\oplus _{\bar x_1\in\bar\phi^{-1}(\bar x_2)}[\frak e_{x_1}].$ 
We shall denote by $\tau':M_1\to M_2$
the homomorphism which is the left vertical arrow in diagram (4.2) 
(note that $\tau '$ is independent of the choice of $x_2\in \Sigma _{X_2}$).

\proclaim
{Lemma 4.12} (The Product Formula) 
The sum $\sum _{\bar x_1\in\bar\phi^{-1}(\bar x_2)} e_{x_1}$
is independent of the choice of $x_2\in \Sigma  _{X_2}$. Set
$n \defeq \sum _{\bar x_1\in\bar\phi^{-1}(\bar x_2)} e_{x_1}
>0$. Then we have
$\tau'=[n] \circ \tau,$ 
where $[n]:M_2\to M_2$ denotes the map of elevation 
to the power $n$ in $M_2$.
\endproclaim

\demo 
{Proof} This follows from the commutativity of diagram (4.2), 
by observing that the homomorphism $\sigma$ being inertia-rigid means that 
the isomorphism $\rho _{x_1}$ in diagram (4.2) equals $p^{c_{x_1}}\tau$ for all 
$\bar x_1\in\bar\phi^{-1}(\bar x_2)$.
\qed
\enddemo

For the rest of this section 
all cohomology groups will be 
continuous Galois cohomology groups, unless otherwise specified.

The Galois-equivariant identification 
$\tau ^{-1}:M_2\isom M_1$ 
induces naturally an injective homomorphism
$H^1(\frak G_2,M_2)\to H^1(\frak G_1,M_1)$
between Galois cohomology groups. 
Indeed, this homomorphism fits into the following commutative diagram: 
$$
\matrix
0&\to &H^1(G_{k_2},M_2)&\to& H^1(\frak G_2,M_2)&\to &H^1(\bar \frak G_2,M_2)\\
&&&&&&\\
&&\downarrow&&\downarrow&&\downarrow\\
&&&&&&\\
0&\to &H^1(G_{k_1},M_1)&\to& H^1(\frak G_1,M_1)&\to &H^1(\bar \frak G_1,M_1),
\endmatrix
$$
in which both rows are exact and vertical maps are natural maps induced by 
$(\sigma, \tau^{-1})$. Here, the left vertical arrow is injective by 
the fact $H^0(H_{k_1}, M_2)=0$, where $H_{k_1}$ stands for the (isomorphic) 
image of $G_{k_1}$ in $G_{k_2}$, and the right vertical arrow is injective since 
$M_2$ is torsion-free and $[\bar\frak G_2:\sigma(\bar \frak G_1)]<\infty$. 
Therefore, the middle vertical arrow is also injective. 

Further, for each $x_2\in \Sigma _{X_2}$, the following diagram is 
commutative:
$$
\CD
H^1(\frak G_1,M_1) @>>> \oplus _{\bar x_1\in\bar\phi^{-1}(\bar x_2)} H^1(\frak I_{\tilde x_1},M_1)\isom
\oplus _{\bar x_1\in\bar\phi^{-1}(\bar x_2)} H^1(\frak J_{\tilde x_1},M_2)\\
@AAA    @AAA  \\
H^1(\frak G_2,M_2) @>>> H^1(\frak I_{\tilde x_2},M_2)\\
\endCD
$$
where the horizontal maps are the natural restriction maps, 
the left vertical map is the above map, the map $H^1(\frak I_{\tilde x_2},M_2)\to
\oplus_{\bar x_1\in\bar\phi^{-1}(\bar x_2)} H^1(\frak J_{\tilde x_1},M_2)$ is the natural map induced by the 
inclusion $\frak J_{\tilde x_1}\subset \frak I_{\tilde x_2}$ for $\bar x_1\in\bar\phi^{-1}(\bar x_2)$, 
and the isomorphism
$H^1(\frak I_{\tilde x_1},M_1)\isom  H^1(\frak J_{\tilde x_1},M_2)$ is naturally induced by the 
natural surjective map $\frak I_{\tilde x_1}\twoheadrightarrow \frak J_{\tilde x_1}$, 
which is induced by $(\sigma,\tau^{-1})$. 

Moreover, we have natural identifications:
$$H^1(\frak I_{\tilde x_1},M_1)\isom \Hom (\frak I_{\tilde x_1},M_1)\isom  \Hom (\frak I_{\tilde x_1}^{\t},M_1)     
\isom \Hom (M_1,M_1)\isom \hat \Bbb Z^{p'},$$
$$H^1(\frak J_{\tilde x_1},M_2)\isom \Hom (\frak J_{\tilde x_1},M_2)\isom \Hom (\frak J_{\tilde x_1}^{\t},M_2)
\isom \Hom (M_2,M_2)\isom \hat \Bbb Z^{p'},$$
and: 
$$ H^1(\frak I_{\tilde x_2},M_2)\isom \Hom (\frak I_{\tilde x_2},M_2)\isom \Hom (\frak I_{\tilde x_2}^{\t},M_2)
\isom \Hom (M_2,M_2)\isom \hat \Bbb Z^{p'}.$$
In light of these identifications, the above diagram can be 
rewritten as:
$$
\CD
H^1(\frak G_1,M_1) @>>>  \oplus _{\bar x_1\in\bar\phi^{-1}(\bar x_2)} \hat \Bbb Z^{p'}\\
@AAA    @A{\oplus _{\bar x_1\in\bar\phi^{-1}(\bar x_2)}[e_{x_1}]}AA  \\
H^1(\frak G_2,M_2) @>>> \hat \Bbb Z^{p'}\\
\endCD
$$
where the vertical map $\hat \Bbb Z^{p'}\to 
\oplus _{\bar x_1\in\bar\phi^{-1}(\bar x_2)} \hat \Bbb Z^{p'}$ is the map 
$\oplus _{\bar x_1\in\phi^{-1}(\bar x_2)}[e_{x_1}]$, and $[e_{x_1}]$ denotes the map of multiplication 
by $e_{x_1}$ in $\hat \Bbb Z^{p'}$. By considering 
all $x_2\in \Sigma _{X_2}$, we obtain the following 
commutative diagram:
$$
\CD
H^1(\frak G_1,M_1) @>>>  \widehat {\Div}_{\overline X_1}\defeq 
\prod' _{\bar x_1\in \Sigma _{\overline X_1}}\hat 
\Bbb Z ^{p'} \isom \prod' _{\bar x_2\in \Sigma _{\overline X_2}}
(\oplus _{\bar x_1\in\phi^{-1}(\bar x_2)} \hat \Bbb Z^{p'}) \\
@AAA    @AAA  \\
H^1(\frak G_2,M_2) @>>> \widehat {\Div}_{\overline X_2}\defeq 
\prod' _{\bar x_2\in \Sigma _{\overline X_2}}\hat 
\Bbb Z^{p'}\\
\endCD 
$$
Here, given an index set $\Lambda$, we define
$\prod'_{\lambda\in \Lambda}\hat\Bbb Z^{p'}\defeq 
\underset{p\nmid n}\to{\varprojlim}(\oplus_{\lambda\in\Lambda} \Bbb Z/n\Bbb Z)$. 
(Accordingly, one has 
$\oplus_{\lambda\in \Lambda}\hat\Bbb Z^{p'}\subset
\prod'_{\lambda\in \Lambda}\hat\Bbb Z^{p'}\subset
\prod_{\lambda\in \Lambda}\hat\Bbb Z^{p'}$, and the equalities hold if and only if 
$\sharp(\Lambda)<\infty$.) 
Thus, the map $\widehat {\Div}_{\overline X_2}\to \widehat {\Div}_{\overline X_1}$ maps
$\bar x_2$ to $\sum _{\bar x_1\in\phi^{-1}(\bar x_2)}e_{x_1}\bar x_1$. 
In particular, the subgroup 
$\widehat {\Div}_{X_2}$ of 
$\widehat {\Div}_{\overline X_2}$ 
maps into the subgroup 
$\widehat {\Div}_{X_1}$ of 
$\widehat {\Div}_{\overline X_1}$. 
Here, for $i=1,2$, $\widehat {\Div}_{X_i}\defeq 
\prod' _{x_1\in \Sigma _{X_i}}\hat \Bbb Z ^{p'}$ 
is naturally embedded into
$\widehat {\Div}_{\overline X_i}$ and is 
regarded as a subgroup of $\widehat {\Div}_{\overline X_i}$. 
Moreover, it follows from various constructions 
that, for $i=1,2$, the image of the map 
$H^1(\frak G_i,M_i) \to  \widehat {\Div}_{\overline X_i}$ 
is contained in $\widehat {\Div}_{X_i}$. Thus, we obtain 
the following 
commutative diagram:
$$
\CD
H^1(\frak G_1,M_1) @>>>  \widehat {\Div}_{X_1}\\
@AAA    @AAA  \\
H^1(\frak G_2,M_2) @>>> \widehat {\Div}_{X_2}\\
\endCD \tag {$4.3$}
$$
Further, for $i=1,2$, set 
$\Div_{X_i}\defeq \oplus _{x_i\in \Sigma _{X_i}}\Bbb Z$, 
which is the group of divisors on $X_i$. 
Then the subgroup
$\Div_{X_2}$ of $\widehat {\Div}_{X_2}$ maps into the subgroup 
$\Div_{X_1}= \oplus _{x_2\in \Sigma _{X_2}}(\oplus _{x_1\in\phi^{-1}(x_2)} \Bbb Z)$ 
of $\widehat {\Div}_{X_1}.$
Thus, we have a natural map
$$\Div _{X_2}\to \Div _{X_1}.$$
We will denote by $\Pri_{X_i}$
the subgroup of $\Div_{X_i}$ which consists of principal divisors.
Note that we have a natural map 
$K_i^{\times}\to \Div _{X_i},$
which maps a function $f_i$ to its divisor $\div (f_i)$ of zeros and poles.
Further, Let $J_{X_i}$ be the Jacobian variety of $X_i$. 
Let $\Div _{X_i}^0\subset\Div_{X_i}$ be the group of degree zero divisors 
on $X_i$. Thus, there exists a natural isomorphism 
$\Div_{X_i}^0/\Pri_{X_i}=J_{X_i}(k_i)$. 
Write $D_{X_i}$ for the kernel of the natural homomorphism 
$\Div _{X_i}^0\to J_{X_i}(k_i)^{p'}$. 
Here, $J_{X_i}(k_i)^{p'}$ stands for the maximal prime-to-$p$ quotient 
$J_{X_i}(k_i)/(J_{X_i}(k_i)\{p\})$ of $J_{X_i}(k_i)$, where, 
for an abelian group $M$, $M\{p\}$ stands 
for the subgroup of torsion elements $a$ of $M$ of 
$p$-power order. 
Then $D_{X_i}$ sits naturally in the following exact sequence:
$$0\to \Pri  _{X_i}\to D_{X_i}\to J_{X_i}(k_i)\{p\} \to 0.$$

For $i\in \{1,2\}$, and a positive integer $n$ prime to $p$, the Kummer exact 
sequence
$$1\to \mu_n\to\Bbb G_m @>[n]>> \Bbb G_m \to 1$$
induces a natural isomorphism
$$ K_i^{\times}/(K_i^{\times})^n\isom  H^1(\frak G_i,\mu _n(K_i^{\sep}))$$
(cf. Lemma 1.4). 
By passing to the projective limit over all integers $n$ prime to $p$, 
we obtain a natural isomorphism
$$(K_i^{\times})^{\wedge p'}\isom H^1(\frak G_i,M_i),$$
where 
$(K_i^{\times})^{\wedge p'}\defeq \underset {p\nmid n}\to 
{\varprojlim}\  K_i^{\times}/(K_i^{\times})^n$. 
As we have a natural embedding
$K_i^{\times}\hookrightarrow (K_i^{\times})^{\wedge p'}$, 
we obtain a natural embedding
$$K_i^{\times} \hookrightarrow  H^1(\frak G_i,M_i).$$
In what follows we will identify $K_i^{\times}$ with its image in $
H^1(\frak G_i,M_i)$; $i=1,2$.
Observe that the natural maps 
$K_i^{\times} \to\Div_{X_i}$ 
and 
$H^1(\frak G_i,M_i)\to\widehat\Div_{X_i}$ 
are compatible with each other, hence that 
the image of $K_i^{\times}$ in $\widehat {\Div}_{X_i}$,
via the map $H^1(\frak G_i,M_i)\to \widehat {\Div}_{X_i}$ in diagram (4.3), 
coincides with the subgroup $\Pri_{X_i}$ of principal divisors. 

\proclaim 
{Lemma 4.13} (Recovering the Multiplicative Group
)  
{\rm (i)} The homomorphism 
$\widehat {\Div}_{X_2}\to \widehat {\Div}_{X_1}$
in diagram (4.3) maps $D_{X_2}$ into $D_{X_1}$. 

\noindent
{\rm (ii)} The above map 
$H^1(\frak G_2,M_2)\to H^1(\frak G_1,M_1)$ induces 
a natural injective (multiplicative) homomorphism 
$$\gamma : K_2^{\times}\hookrightarrow 
(K_1^{\times})^{p^{-n}}= (K_1^{p^{-n}})^{\times},$$ 
where $p^n$ is the exponent 
of the $p$-primary finite abelian group $J_{X_1}(k_1)\{p\}$. 
We have 
$[\gamma(K_2^{\times}): \gamma(K_2^{\times})\cap K_1^{\times}]<\infty$ 
and 
$[\gamma(K_2^{\times}): \gamma(K_2^{\times})\cap (K_1^{\times})^p]>1$. 

Moreover, this injective homomorphism
is functorial in the following sense: 
Let $\frak H_1\subset\frak G_1$, $\frak H_2\subset\frak G_2$ 
be open subgroups such that $\sigma(\frak H_1)\subset\frak H_2$, 
and, for $i=1,2$, let $L_i/K_i$ be 
the separable extension corresponding to $\frak H_i\subset\frak G_i$, 
$Y_i$ the integral closure of $X_i$ in $L_i$, 
and 
$\ell_i$ the constant field of $L_i$ (i.e., the algebraic closure of 
$k_i$ in $L_i$). Then we have a commutative diagram:
$$
\CD
L_2^{\times} @>>> (L_1^{\times})^{p^{-m}}\\
@AAA                        @AAA  \\
K_2^{\times} @>>>  (K_1^{\times})^{p^{-n}}\\
\endCD
$$
where $p^{m}\ge p^n$ is the exponent 
of the $p$-primary finite abelian group
$J_{Y_1}(\ell_1)\{p\}$, and 
the vertical arrows are the natural embeddings.
\endproclaim

\demo
{Proof} 
(i) 
We have the following diagram of maps:
$$
\CD
\Div_{X_1} @>>>   H^2(\pi_1(X_1)^{(p')},M_1)\\
@AAA         @AAA  \\
\Div_{X_2} @>>>   H^2(\pi_1(X_2)^{(p')},M_2)\\
\endCD
$$
where the map $\Div_{X_2}\to \Div_{X_1}$ is the one induced by the map
$\widehat {\Div}_{X_2}\to \widehat {\Div}_{X_1}$ in diagram (4.3). 
For $i\in \{1,2\}$, the group $H^2(\pi_1(X_i)^{(p')},M_i)$ denotes the second cohomology 
group of the profinite group $\pi_1(X_i)^{(p')}$ with coefficients in the 
(continuous) $\pi_1(X_i)^{(p')}$-module $M_i$. 

First, we shall treat the special case that $(g_1\geq)g_2>0$. In this case, 
we have a natural isomorphism
$H^2(\pi_1(X_i)^{(p')},M_i)\isom H^2_{\et}(X_i,M_i)$ (cf. [Mochizuki4], 
Proposition 1.1), where $H^2_{\et}(X_i,M_i)$ denotes the second 
\'etale cohomology group of $X_i$ with coefficients in $M_i$. 
In what follows we will identify the groups $H^2(\pi_1(X_i)^{(p')},M_i)$ 
and $H^2_{\et}(X_i,M_i)$ via the above 
identifications. Further, the map
$H^2(\pi_1(X_2)^{(p')},M_2)\to H^2(\pi_1(X_1)^{(p')},M_1)$ is the map induced by 
the natural map $\pi_1(X_1)^{(p')}\to \pi_1(X_2)^{(p')}$ between  
fundamental groups, which is induced by $\sigma$ (cf. Lemma 2.6),
and the Galois-equivariant 
identification $\tau ^{-1}:M_2\isom  M_1$. The map $\Div_{X_i}\to   
H^2(\pi_1(X_i)^{(p')},M_i)$ maps a divisor $D$ to its 
first arithmetic (\'etale) Chern class $c_1(D)$, and is naturally induced by the 
Kummer exact sequence $1\to\mu_n\to\Bbb G_m @>[n]>> \Bbb G_m \to 1$
in \'etale topology (cf. [Mochizuki2], 4.1). 
In particular, the map $\Div_{X_i}\to  H^2(\pi_1(X_i)^{(p')},M_i)$ factors as
$\Div_{X_i}\to  \Pic (X_i)/(J_{X_i}(k_i)\{p\})\hookrightarrow H^2(\pi_1(X_i)^{(p')},M_i)$, where
$\Pic (X_i)\defeq H^1_{\et}(X_i,\Bbb G_m)$ is the Picard group of $X_i$. 
Note that the kernel of 
the above map $\Div_{X_i}\to   H^2(\pi_1(X_i)^{(p')},M_i)$ coincides with  
$D_{X_i}$. We claim that the above diagram is commutative. Thus, 
it induces a natural map $D_{X_2}\to D_{X_1}$, as desired (in the case that $g_2>0$). 

To prove the above claim, 
let $x_2\in \Sigma _{X_2}$. 
We shall investigate the images of 
$x_2\in \Div_{X_2}$ in  $H^2(\pi_1(X_1)^{(p')},M_1)$ under the two (composite) maps 
in the above diagram. First, consider the special case where $x_2\in\Sigma_{X_2}$ 
is $k_2$-rational and each point of $\phi ^{-1}(x_2)\subset\Sigma_{X_1}$ 
is $k_1$-rational. 
Then the image $c_1(x_2)$ of the divisor 
$x_2\in \Div_{X_2}$ in  $H^2(\pi_1(X_2)^{(p')},M_2)$ coincides with the class of the 
extension $1\to M_2\to \pi_1(\Bbb L_{x_2}^{\times})^{(p')} \to \pi_1(X_2)^{(p')}\to 1$, 
where $\pi_1(\Bbb L_{x_2}^{\times})^{(p')}$ is the geometrically prime-to-$p$ 
fundamental group of the geometric line bundle $\Bbb L_{x_2}$ corresponding 
to the invertible sheaf $\Cal O_{X_2}(x_2)$ with the zero section removed 
(cf. [Mochizuki3], Lemma 4.2 and [Mochizuki2], 4.1). Further, $\pi_1(\Bbb L_{x_2}^{\times})^{(p')}$ is 
naturally identified 
with the maximal tame cuspidally central quotient 
$\pi_1(X_2\smallsetminus\{x_2\})^{(p'),\cc}$ of $\pi_1(X_2\smallsetminus\{x_2\})^{(p')}$ (cf. [Mochizuki3], Lemma 4.2(iii)). 
Similarly, the maximal tame
cuspidally central quotient $\pi_1(X_1\smallsetminus\phi^{-1}(x_2))^{(p'),\cc}$ 
of $\pi_1(X_1\smallsetminus\phi^{-1}(x_2))^{(p')}$ gives the extension of 
$\pi_1(X_1)^{(p')}$ by $\oplus _{x_1\in\phi^{-1}(x_2)}M_1$ that corresponds 
to $(c_1(x_1))_{x_1\in\phi^{-1}(x_2)}\in 
\oplus _{x_1\in\phi^{-1}(x_2)}H^2(\pi_1(X_1)^{(p')},M_1)=
H^2(\pi_1(X_1)^{(p')},\oplus _{x_1\in\phi^{-1}(x_2)}M_1)
$. 
Being well-behaved (with respect to $\tilde\phi$), $\sigma$ 
induces naturally a homomorphism 
$\pi_1(X_1\smallsetminus\phi^{-1}(x_2))^{(p')}\to\pi_1(X_2\smallsetminus\{x_2\})^{(p')}$, 
which is a lifting of 
$\pi_1(X_1)^{(p')}\to\pi_1(X_2)^{(p')}$ and 
which further induces a homomorphism 
$\pi_1(X_1\smallsetminus\phi^{-1}(x_2))^{(p'),\cc}\to\pi_1(X_2\smallsetminus\{x_2\})^{(p'),\cc}$. 
These homomorphisms fit into the following commutative diagram: 
$$
\matrix
1&\to& \oplus _{x_1\in\phi^{-1}(x_2)}M_1 & \to & \pi_1(X_1\smallsetminus\phi^{-1}(x_2))^{(p'),\cc} 
&\to&  \pi_1(X_1)^{(p')} &\to& 1 \\
&&&&&&&&\\
&&\downarrow&&\downarrow&&\downarrow&&\\
&&&&&&&&\\
1&\to& M_2 & \to & \pi_1(X_2\smallsetminus\{x_2\})^{(p'),\cc} 
&\to&  \pi_1(X_2)^{(p')} &\to& 1 \\
\endmatrix
$$
in which both rows are exact and the left vertical arrow is 
$\oplus _{x_1\in\phi^{-1}(x_2)}e_{x_1}\tau$, by the 
inertia-rigidity of $\sigma$. 
The commutativity of this last diagram implies that the image of 
the extension class 
of the top row (i.e., $(c_1(x_1))_{x_1\in\phi^{-1}(x_2)}$) 
under the map 
$H^2(\pi_1(X_1)^{(p')},\oplus _{x_1\in\phi^{-1}(x_2)}M_1)
\to H^2(\pi_1(X_1)^{(p')},M_1)$ 
induced by $\oplus _{x_1\in\phi^{-1}(x_2)}[e_{x_1}]$ coincides 
with the image of the extension class 
of the bottom row (i.e., $c_1(x_2)$) under 
the map 
$H^2(\pi_1(X_2)^{(p')},M_2) \to H^2(\pi_1(X_1)^{(p')},M_1)$ 
induced by $\sigma$ and $\tau^{-1}$. In other words, 
the image of $c_1(x_2)$ in $H^2(\pi_1(X_1)^{(p')},M_1)$ 
coincides with 
$\sum_{x_1\in\phi^{-1}(x_2)}e_{x_1}c_1(x_1)$. From this follows the claim 
(in the special case), since the 
divisor $x_2$ maps to $\sum _{x_1\in\phi^{-1}(x_2)}e_{x_1}x_1$ via the above 
map $\Div _{X_1}\to \Div _{X_2}$.  
Finally, consider the general case where $x_2$ may not be $k_2$-rational and 
each point of $\phi^{-1}(x_2)$ may not be $k_1$-rational. But this is 
reduced to the special case by considering suitable open subgroups of $\frak G_i$, 
$i=1,2$, corresponding to constant field extensions $k_i'$ of $k_i$. 
(Here, use the fact that the natural map $H^2(\pi_1(X_i)^{(p')},M_i)\to 
H^2(\pi_1(X_i\times_{k_i}k_i')^{(p')}, M_i)$ is injective, which follows from 
the injectivity of the natural map 
$J_{X_i}(k_i)\to J_{X_i}(k_i')=J_{X_i\times_{k_i}k_i'}(k_i')$.) 
Thus, the claim follows. 

Next, to treat the general case that we may possibly have $g_2=0$, consider 
any open subgroup $\frak H_2$ of $\frak G_2$ and set 
$\frak H_1\defeq \sigma^{-1}(\frak H_2)$, which is an open subgroup of $\frak G_1$. 
For each $i=1,2$, let $Y_i$ be the cover of $X_i$ corresponding to 
the open subgroup $\frak H_i\subset \frak G_i$, and $\ell_i$ the constant 
field of $Y_i$ (i.e., the algebraic closure of $k_i$ in the function field 
of $Y_i$). Now, assume that the genus of $Y_2$ is $>0$. 
Then it follows from the preceding argument that 
the homomorphism $\Div_{Y_1}\to \Div_{Y_2}$ maps $D_{Y_1}$ into $D_{Y_2}$. 
In particular, by functoriality, the image of $D_{X_1}$ in $\Div_{X_2}$ is 
mapped into $D_{Y_2}\subset \Div_{Y_2}$ under the natural map 
$\Div_{X_2}\to\Div_{Y_2}$. Or, equivalently, the image of 
$D_{X_1}$ in $\Div_{X_2}/D_{X_2}$ lies in the kernel of 
$\Div_{X_2}/D_{X_2}\to \Div_{Y_2}/D_{Y_2}$. This last map is identified with 
the natural map 
$\Pic_{X_2}/(J_{X_2}(k_2)\{p\})\to 
\Pic_{Y_2}/(J_{Y_2}(\ell_2)\{p\})$ 
induced by the pull-back of line bundles. Thus, by considering the 
norm map, we see that the kernel in question is killed by the 
degree $[\frak G_2: \frak H_2]$ of the cover $Y_2\to X_2$, 
hence so is the image of $D_{X_1}$ in $\Div_{X_2}/D_{X_2}$. 

Observe that the greatest common divisor of $[\frak G_2: \frak H_2]$, where 
$\frak H_2$ runs over all open subgroups of $\frak G_2$ such that 
the corresponding cover 
has genus $>0$, is $1$. (Indeed, if 
$g_2>0$, this is trivial, and, if $g_2=0$, this follows, 
for example, from Kummer theory.) Thus,  
the image of $D_{X_1}$ in $\Div_{X_2}/D_{X_2}$ must be trivial, as desired. 

\noindent
(ii) For $i=1,2$, let $\tilde D_{X_i}$ denote the inverse image of 
$D_{X_i}\subset\Div_{X_i}\ (\subset\widehat\Div_{X_i})$ in 
$H^1(\frak G_i,M_i)$. 
It follows from (i) and the commutativity of diagram (4.3) that 
the natural injective homomorphism 
$H^1(\frak G_2, M_2)\hookrightarrow H^1(\frak G_1, M_1)$ induces 
a natural injective homomorphism 
$\tilde D_{X_2}\hookrightarrow \tilde D_{X_1}$. 
Since 
$K_i^{\times}$ is the inverse image of $\Pri_{X_i}\subset\Div_{X_i}$ in 
$H^1(\frak G_i,M_i)$ (cf. [Mochizuki4], Proposition 2.1(ii)), we have 
$\tilde D_{X_i}/K_i^{\times}\isom D_{X_i}/\Pri_{X_i}\isom J_{X_i}(k_i)\{p\}$. 
Thus, the above injective homomorphism $\tilde D_{X_2}\hookrightarrow \tilde D_{X_1}$ 
induces $(K_2^{\times})^{p^n}\hookrightarrow K_1^{\times}$, or, equivalently, 
$K_2^{\times}\hookrightarrow (K_1^{\times})^{p^{-n}}$. 

Since $\gamma(K_2^{\times})/(\gamma(K_2^{\times})\cap K_1^{\times})$ is 
injectively mapped into 
$\tilde D_{X_1}/K_1^{\times}\isom J_{X_1}(k_1)\{p\}$, which is 
finite, $\gamma(K_2^{\times})\cap K_1^{\times}$ is of finite index 
in $\gamma(K_2^{\times})$. 
Next, suppose that $\gamma(K_2^{\times})=\gamma(K_2^{\times})\cap (K_1^{\times})^p$, 
or, equivalently, $\gamma(K_2^{\times})\subset (K_1^{\times})^p$. By the assumption 
that $a=0$, there exists an $x_1\in\Sigma_{X_1}$ such that $e_{x_1}=\frak e_{x_1}^{\t}$. 
In particular, then $e_{x_1}$ is prime to $p$. Set $x_2\defeq \phi(x_1)\in\Sigma_{X_2}$ and take any 
$g\in K_2^{\times}$ such that $\ord_{x_2}(g)=1$. Then, by the commutativity of 
diagram (4.3), we have $\ord_{x_1}(\gamma(g))=e_{x_1}\ord_{x_2}(g)=e_{x_1}$, which is prime to $p$. 
On the other hand, since $\gamma(g)\in (K_1^{\times})^p$, $\ord_{x_1}(\gamma(g))$ must 
be divisible by $p$, which is absurd. 

Finally, the desired commutativity of diagram follows easily from the functoriality 
of Kummer theory. 
\qed
\enddemo

Next, let $x_1\in \Sigma_{X_1}$ and set 
$x_2\defeq \phi(x_1)\in\Sigma_{X_2}$. Then (by choosing $\tilde x_1\in\Sigma_{\tilde X_1}$ above $x_1$ 
and $\tilde x_2\in\Sigma_{\tilde X_2}$ above $x_2$ such that $\tilde \phi(\tilde x_1)=\tilde x_2$) 
we have the following natural commutative diagram: 
$$
\matrix
H^1(\frak G_1, M_1)&\to&H^1(\frak D_{\tilde x_1}, M_1)&\to&H^1(\frak I_{\tilde x_1},M_1)\\
&&&&\\
\uparrow&&\uparrow&&\uparrow\\
&&&&\\
H^1(\frak G_2, M_2)&\to&H^1(\frak D_{\tilde x_2},M_2)&\to&H^1(\frak I_{\tilde x_2},M_2),
\endmatrix
$$
where the horizontal arrows are natural restriction maps and 
the vertical arrows are induced by $(\sigma,\tau^{-1})$. By Kummer theory, this diagram 
can be identified with the following natural commutative diagram: 
$$
\matrix
(K_1^{\times})^{\wedge p'}&\to&((K_1)_{x_1}^{\times})^{\wedge p'}&
\overset{\ord_{x_1}}\to{\to}&\hat\Bbb Z^{p'}\\
&&&&\\
\uparrow&&\uparrow&&\uparrow\\
&&&&\\
(K_2^{\times})^{\wedge p'}&\to&((K_2)_{x_2}^{\times})^{\wedge p'}&
\overset{\ord_{x_2}}\to{\to}&\hat\Bbb Z^{p'},\\
\endmatrix\tag{$4.4$}
$$
where the left horizontal arrows in the two rows arise from natural field 
homomorphisms $K_1\to (K_1)_{x_1}$ and $K_2\to (K_2)_{x_2}$ and the vertical arrows are 
induced by $(\sigma,\tau^{-1})$. Further, 
the kernels of $((K_1)_{x_1}^{\times})^{\wedge p'}\overset{\ord_{x_1}}\to{\to}\hat\Bbb Z^{p'}$ and 
$((K_2)_{x_2}^{\times})^{\wedge p'}\overset{\ord_{x_2}}\to{\to}\hat\Bbb Z^{p'}$ 
are naturally identified with 
$H^1(G_{k(x_1)},M_1)=(k(x_1)^{\times})^{\wedge p'}=k(x_1)^{\times}$ and 
$H^1(G_{k(x_2)},M_2)=(k(x_2)^{\times})^{\wedge p'}=k(x_2)^{\times}$, respectively. Thus, in particular, 
the above homomorphism $((K_2)_{x_2}^{\times})^{\wedge p'}\to((K_1)_{x_1}^{\times})^{\wedge p'}$ 
induces naturally a homomorphism $\iota_{x_1}: k(x_2)^{\times}\to k(x_1)^{\times}$ that is identified 
with the homomorphism 
$H^1(G_{k(x_2)},M_2)\to H^1(G_{k(x_1)},M_1)$ induced by $(\sigma,\tau^{-1})$. 
Here, the last homomorphism is injective by 
the fact $H^0(H_{k(x_1)}, M_2)=0$, where $H_{k(x_1)}$ stands for the (isomorphic) 
image of $G_{k(x_1)}$ in $G_{k(x_2)}$, which is open in $G_{k(x_2)}$. 

We have two natural field homomorphisms $K_1\to K_1^{p^{-n}}$: 
the first one is a natural embedding $i: K_1\hookrightarrow K_1^{p^{-n}}$ 
of degree $p^n$ and the second one is the isomorphism 
$j: K_1\isom K_1^{p^{-n}}$ induced by the $p^{-n}$-th power map. 
According to these we obtain two scheme morphisms 
$X_1^{p^{-n}}\to X_1$, where 
$X_1^{p^{-n}}$ stands for the integral closure of $X_1$ in $K_1^{p^{-n}}$. 
First, for closed points, these two morphisms give the same 
bijection $\pi:\Sigma_{X_1^{p^{-n}}}\isom \Sigma_{X_1}$. So, let $x_1\in\Sigma_{X_1}$ 
and set $x_1^{p^{-n}}\defeq \pi^{-1}(x_1)$. The two field homomorphism $i$ and $j$ induce 
two isomorphisms $k(x_1)\to k(x_1^{p^{-n}})$ of residue fields, which we shall denote by 
$\bar i(x_1)$ and $\bar j(x_1)$, respectively. Then we have 
$\bar i(x_1)=F^n\circ \bar j(x_1)$, where $F$ stands for the $p$-th power Frobenius map. 
Now, for valuations of functions, we have 
$\ord_{x_1^{p^{-n}}}\circ i=p^n\ord_{x_1}$ and 
$\ord_{x_1^{p^{-n}}}\circ j=\ord_{x_1}$. Finally, for values of functions, 
we have 
$i(f)(x_1^{p^{-n}})=\bar i(x_1)(f(x_1))$ 
and 
$j(f)(x_1^{p^{-n}})=\bar j(x_1)(f(x_1))$ 
for each $f\in K_1^{\times}$ with $\ord_{x_1}(f)\geq 0$. Thus, in particular, 
$i(f)(x_1^{p^{-n}})=j(f)(x_1^{p^{-n}})^{p^n}$. 

\proclaim
{Lemma 4.14} Let $\gamma : K_2^{\times}\hookrightarrow (K_1^{\times})^{p^{-n}}$ 
be the injective homomorphism in Lemma 4.13. Let $x_1\in \Sigma_{X_1}$ and set 
$x_2\defeq \phi(x_1)\in\Sigma_{X_2}$. 
Then: 

\noindent
(i) For each $g\in K_2^{\times}$, we have $\ord_{x_1^{p^{-n}}}(\gamma(g))=p^ne_{x_1}\ord_{x_2}(g)$. 
(Namely, $\gamma$ is order-preserving with respect $\pi^{-1}\circ\phi$. See Definition 5.1.) 

\noindent
(ii) For each $g\in K_2^{\times}$ with $\ord_{x_2}(g)=0$, we have 
$(\gamma(g))(x_1^{p^{-n}})=i(x_1)(\iota_{x_1}(g(x_2)))$. 
(Namely, $\gamma$ is value-preserving with respect $\pi^{-1}\circ\phi$ 
and $\{i(x_1)\circ \iota_{x_1}\}_{x_1^{p^{-n}}\in \Sigma_{X_1^{p^{-n}}}}$. 
See Definition 5.2.)
\endproclaim

\demo{Proof}
(i) and (ii) follow immediately from the commutativity of diagrams (4.3) and (4.4). 
\qed\enddemo

Fix a prime number $l\neq p$. For each $i=1,2$, let $k_i^l$ 
be the (unique) $\Bbb Z_l$-extension of $k_i$, set 
$K_i^l\defeq K_ik_i^l$, 
and write $X_i^l$ for the normalization of $X_i$ in $K_i^l$. 
(Thus, $X_i^l=X_i\times_{k_i}k_i^l$.) 
Then the $p$-primary 
abelian subgroup $J_{X_i}(k_i^l)\{p\}$ 
of $J_{X_i}(k_i^l)$ is finite for $i=1,2$. 
(See, e.g., [Rosen], Theorem 11.6, or, 
alternatively, [ST], Proof of Theorem 3.7.) 
So, write $p^{n_0}$ for the exponent of $J_{X_1}(k_1^l)\{p\}$. 
By passing to the limit over the finite extensions of $k_i$ 
contained in $k_i^l$ for $i=1,2$ (cf. Lemma 4.13(ii)), 
we get a natural embedding 
$(K_2^l)^{\times} \hookrightarrow 
((K_1^l)^{\times})^{p^{-{n_0}}}$. 
Now, we shall apply a result from \S 5. (Observe that there are no vicious 
circles since the discussion of \S 5 does not depend on the contents of 
earlier sections.) More specifically, by Lemma 4.14 and Proposition 5.3, 
the above embedding $(K_2^l)^{\times} \hookrightarrow 
((K_1^l)^{\times})^{p^{-{n_0}}}$ 
arises from an (a uniquely determined) embedding $K_2^l
\hookrightarrow (K_1^l)^{p^{-{n_0}}}$ of fields. 
This embedding of fields 
restricts to the original embedding of multiplicative groups
$K_2^{\times} \hookrightarrow (K_1^{\times})^{p^{-{n}}}$. 
Thus, we conclude that this original embedding also arises from an 
(a uniquely determined)  
embedding $K_2\hookrightarrow K_1^{p^{-{n}}}$ of fields. 

Define the subfields $K_2\supset K_2'\supset K_2''$ to be the inverse images of
the subfields $K_1^{p^{-{n}}}\supset K_1\supset K_1^p$ in $K_2$. By Lemma 4.13(ii), 
there exists a finite subset $S\subset K_2$ such that 
$K_2=\cup_{\alpha\in S}K_2'\alpha$. As $K_2$ is an infinite field, this 
implies that $K_2'$ is also an infinite field and that $K_2$ must be 
of dimension $1$ as a $K_2'$-vector space. Namely, $K_2=K_2'$, or, equivalently, 
the above field homomorphism $K_2\hookrightarrow K_1^{p^{-{n}}}$ 
induces a field homomorphism $\gamma: K_2\hookrightarrow K_1$. Next, 
again by Lemma 4.13(ii), we have $[K_2^{\times}: (K_2'')^{\times}]>1$, 
i.~e., $K_2\supsetneq K_2''$. Equivalently, 
the field homomorphism $K_2\hookrightarrow K_1$ is separable. 

Passing to the open subgroups $\frak H_1\subset \frak G_1$, $\frak H_2\subset\frak G_2$ 
with $\sigma(\frak H_1)\subset\frak H_2$ and applying the above arguments to 
$\frak H_1\overset\sigma\to\to\frak H_2$, we obtain naturally a (separable) field homomorphism 
$\tilde\gamma: \tilde K_2\to\tilde K_1$ which restricts to the above (separable) field homomorphism 
$K_2\to K_1$. 

\proclaim
{Lemma 4.15} (Compatibility with the Galois Action) Let 
$g_1\in \frak G_1$, and 
let $g_2\defeq \sigma (g_1)\in \frak G_2$. Then the following diagram is commutative: 
$$
\CD
\tilde K_2 @>\tilde \gamma>> \tilde K_1 \\
@Ag_2AA         @Ag_1AA  \\
\tilde K_2 @>\tilde \gamma>> \phantom{.}\tilde K_1. 
\endCD
$$
\endproclaim

\demo
{Proof} Let $\frak H_2\subset\frak G_2$ be an open normal subgroup 
and set $\frak H_1\defeq\sigma^{-1}(\frak H_2)$, which is an 
open normal subgroup of $\frak G_1$. For $i=1,2$, let 
$F_i/K_i$ be the finite Galois sub-extension of 
$\tilde K_i/K_i$ corresponding to $\frak H_i\subset\frak G_i$, 
and denote by $Y_i$ the integral closure of $X_i$ in $F_i$. 
We have commutative diagrams:
$$
\CD
H^1(\frak H_2,M_2)@>>> H^1(\frak H_1,M_1)\\
@Ag_2AA            @Ag_1AA \\
H^1(\frak H_2,M_2)@>>> H^1(\frak H_1,M_1) \\
\endCD
$$
where $g_i:H^1(\frak H_i,M_i)\to H^1(\frak H_i,M_i)$ denotes the automorphism 
induced by the action of $g_i$ on $\frak H_i$, and the horizontal maps are 
naturally induced by $(\sigma,\tau^{-1})$ 
(cf. Lemma 4.11(i)), and:
$$
\CD
\widehat {\Div}_{Y_2} @>>> \widehat {\Div}_{Y_1}\\
@Ag_2AA            @Ag_1AA\\
\widehat {\Div}_{Y_2} @>>> \widehat {\Div}_{Y_1}\\
\endCD
$$
where the map $g_i:\widehat {\Div}_{Y_i}\to \widehat {\Div}_{Y_i}$
is the automorphism naturally induced by the action of $g_i$ on
$Y_i$ (cf. Remark 4.2(iv)). Further, the above diagrams commute with each other, via the 
maps $H^1(\frak H_i,M_i)\to \widehat {\Div}_{Y_i}$ in diagram (4.3) for $i=1,2$. 
Note that in the above diagrams the map $g_i:H^1(\frak H_i,M_i)\to H^1(\frak H_i,M_i)$ 
restricted to $F_i^{\times}$ coincides with the automorphism
$g_i:F_i^{\times}\to F_i^{\times}$. Therefore, we deduce the following 
commutative diagram:
$$
\CD
F_2^{\times} @>\tilde\gamma>> F_1^{\times} \\
@Ag_2AA         @Ag_1AA  \\
F_2^{\times} @>\tilde \gamma>> F_1^{\times} \\
\endCD
$$
The assertion follows from this. 
\qed
\enddemo

Finally, we shall prove the uniqueness of the field homomorphism 
$\tilde\gamma:\tilde K_2\to\tilde K_1$ that is Galois-compatible 
with respect to $\sigma$ and restricts to a separable homomorphism 
$K_2\to K_1$. In the profinite case, 
this uniqueness follows formally from the uniqueness in the 
assertion of the Isom-form proved in [Uchida2], as in 
the case of rigid homomorphisms in \S 3. (Observe that 
$\tilde\gamma:\tilde K_2\to\tilde K_1$ is then an isomorphism.) 
In general, however, we need some arguments which are not 
entirely formal, as follows. 

So, Let $\tilde\gamma':\tilde K_2\to\tilde K_1$ be 
another such field homomorphism. The 
field homomorphisms $\tilde\gamma$ and $\tilde\gamma'$ 
induce field isomorphisms $\bar k_2\isom \bar k_1$, say, $\bar \gamma$ and $\bar \gamma'$, 
respectively, 
which are Galois-compatible with respect to $\sigma$. We may 
write $\bar\gamma'=\varphi_1^{\alpha}\circ\bar\gamma$ for some $\alpha\in\hat\Bbb Z$, 
where $\varphi_1\in\Gal(\bar k_1/\Bbb F_p)$ stands for the $p$-th power 
Frobenius element. 
Further, the isomorphisms $\bar \gamma$ and $\bar \gamma'$ induce
$\hat\Bbb Z^{p'}$-module isomorphisms $M_2\isom M_1$, say, 
$\tau^{-1}$ and $(\tau')^{-1}$, respectively, which are 
Galois-compatible with respect to $\sigma$. Thus, we have 
$(\tau')^{-1}=[p^\alpha]\circ\tau^{-1}$. By Kummer theory, we have the following commutative diagrams: 
$$
\matrix
K_1^{\times}&\hookrightarrow&H^1(\frak G_1,M_1)\\
&&\\
\gamma\uparrow\phantom{\gamma}&&\phantom{(\sigma,\tau^{-1})}\uparrow(\sigma,\tau^{-1})\\
&&\\
K_2^{\times}&\hookrightarrow&\phantom{,}H^1(\frak G_2,M_2), 
\endmatrix
$$
and: 
$$
\matrix
K_1^{\times}&\hookrightarrow&H^1(\frak G_1,M_1)\\
&&\\
\gamma'\uparrow\phantom{\gamma'}&&\phantom{(\sigma,(\tau')^{-1})}\uparrow(\sigma,(\tau')^{-1})\\
&&\\
K_2^{\times}&\hookrightarrow&\phantom{.}H^1(\frak G_2,M_2).
\endmatrix
$$
Thus, for each $g\in K_2^{\times}$, we have $\gamma'(g)=\gamma(g)^{p^\alpha}$ in 
$(K_1^{\times})^{\wedge p'}$. 
Since both $\gamma$ and $\gamma'$ are field homomorphisms, we deduce $p^{\alpha}\in\Bbb Q_{>0}$, 
by taking a non-constant function $g$ and considering valuations at suitable points. Thus, 
$\alpha\in\Bbb Z$. Exchanging $\gamma$ and $\gamma'$ if necessary, we may assume that 
$\alpha\geq 0$. Thus, $\gamma'=F^\alpha\circ \gamma$, where $F$ stands for the $p$-th power 
Frobenius map. Since $\gamma'$ is separable, we conclude $\alpha=0$, hence $\gamma'=\gamma$. 
Passing to the open subgroups $\frak H_1\subset\frak G_1$, $\frak H_2\subset\frak G_2$ 
with $\sigma(\frak H_1)\subset \frak H_2$, we conclude that $\tilde\gamma:\tilde K_2\to\tilde K_1$ 
is unique. 

Thus, the proof of Theorem 4.8 is completed. \qed

\subhead 
\S 5. Recovering the additive structure
\endsubhead

This section is devoted to the proof of Proposition 5.3, used in 
the proof of Theorem 4.8 in \S 4. We shall 
first axiomatize the set-up. We will use the following notations.
For $i\in \{1,2\}$, let $X_i$ be a proper, smooth, geometrically connected 
curve over a field $k_i$ of characteristic 
$p_i\ge 0$. Let $K_i=K_{X_i}$ be the function field of $X_i$, and 
$\Sigma _{X_i}$ the set of closed points of $X_i$. 
Let
$$\iota: K_2^{\times} \hookrightarrow K_1 ^{\times}$$ 
be an embedding between multiplicative groups, which we extend to an embedding 
$\iota: K_2 \hookrightarrow K_1$ between multiplicative monoids 
by setting $\iota (0)=0$. 
We assume that we are given a map
$$\phi :\Sigma _{X_1}\to \Sigma _{X_2}$$
which has finite fibers, i.e., for any $x_2\in X_2$, 
the inverse image $\phi ^{-1}(x_2)\subset \Sigma _{X_1}$ is a finite set.

\definition
{Definition 5.1} (Order-Preserving Maps) The map 
$\iota : K_2\to K_1$ is called order-preserving with respect to the 
map $\phi$, if for any $x_2\in \Sigma _{X_2}$ and any $x_1\in 
\phi ^{-1}(x_2)$, there exists positive integers $e_{x_1x_2}>0$ such
that the following diagram commutes:
$$
\CD
K_1 @>\ord _{x_1}>> \Bbb Z\cup \{\infty\} \\ 
@A{\iota}AA       @A{[e_{x_1x_2}]}AA \\
K_2 @>\ord _{x_2}>> \Bbb Z\cup \{\infty\} \\
\endCD
$$
Here, $[e_{x_1x_2}]$ denotes the map of multiplication by $e_{x_1x_2}$
in $\Bbb Z$, which we extend naturally to $\Bbb Z\cup \{\infty\}$ by
mapping $\infty$ to $\infty$.
\enddefinition

Next, we assume that the map $\iota : K_2\to K_1$ is order-preserving with 
respect to the map $\phi: \Sigma_{X_1}\to\Sigma_{X_2}$. Further, we assume that 
we are given 
an embedding
$$\iota _{x_1x_2}:k (x_2)^{\times}\hookrightarrow k (x_1)^{\times}$$
between multiplicative groups for 
any $x_2\in \Sigma _{X_2}$ and any $x_1\in \phi ^{-1}(x_2)$. 

\definition
{Definition 5.2} (Value-Preserving Maps) The map 
$\iota : K_2\hookrightarrow K_1$ is called value-preserving with respect to the 
map $\phi$ and the embeddings $\{\iota _{x_1x_2}\}_{(x_1,x_2)}$, 
where $(x_1,x_2)$ runs over all pairs of points
$x_2\in \Sigma _{X_2}$ and $x_1\in \phi ^{-1}(x_2)$,
if for any $f_2\in K_2^{\times}$ and any point $x_2\in \Sigma _{X_2}$ such that
$x_2\cap \Supp \div (f_2)=\emptyset $, the following equality holds:
$$\iota _{x_1,x_2}(f_2(x_2))=\iota (f_2)(x_1)$$
where $f_2(x_2)$ (resp. $\iota (f_2)(x_1)$) denotes the value of $f_2$ at 
$x_2$ (resp. the value of $\iota (f_2)$ at $x_1$).
\enddefinition

Note that if $\iota : K_2\hookrightarrow K_1$ is value-preserving, 
it particularly fits into the following commutative diagram:
$$
\CD
k (x_2)^{\times} @>\iota_{x_1x_2}>> k (x_1)^{\times} \\
@AAA                   @AAA \\
k_2^{\times} @>\iota >> k_1^{\times} \\
\endCD
$$
where the vertical maps are the natural embeddings. (Observe that $\iota$
maps $k_2$ into $k_1$, by the order-preserving assumption.) 

\proclaim
{Proposition 5.3} (Recovering the Additive Structure) Assume that
the embedding $\iota : K_2\hookrightarrow K_1$ is order-preserving with respect to the 
map $\phi$, and value-preserving with respect to the 
map $\phi$ and the embeddings $\{\iota _{x_1x_2}\}_{(x_1,x_2)}$,
where the pair $(x_1,x_2)$ runs over all points $x_2\in \Sigma _{X_2}$ and
$x_1\in \phi ^{-1}(x_2)$. 
Assume further that 
$X_2(k_2)$ is an infinite set. 
Then the map $\iota $ is additive (hence, a homomorphism of fields). 
\endproclaim

\demo
{Proof} 
First, we shall prove that $\iota^{-1}(k_1)=k_2$. 
(Namely, $f\in K_2$ is constant if and only if $\iota(f)\in K_1$ is constant.) 
Indeed, 
set $F_2\defeq\iota^{-1}(k_1)$. Note that $k_i^{\times}$ 
coincides with the set of functions in $K_i^{\times}$ with neither 
zeroes nor poles (or, equivalently, with no poles) 
anywhere in $\Sigma_{X_i}$. Now, by the 
order-preserving property of $\iota$, 
$
F_2\smallsetminus\{0\}$ 
coincides with the set of functions in $K_2^{\times}$ with 
neither zeroes nor poles (or, equivalently, with no poles) 
in $\phi(\Sigma_{X_1})\subset \Sigma_{X_2}$. It follows 
easily from this characterization that $F_2$ is a subfield of $K_2$ 
containing $k_2$. Since $K_2$ is a function field of one variable
over $k_2$ and since $k_2$ is algebraically closed in $K_2$, we have 
either $F_2=k_2$ or that $F_2$ is also a function field of one variable 
over $k_2$. Suppose 
the latter, and let $W_2$ be the (proper, smooth, geometrically 
connected) curve over $k_2$ with function field $F_2$. Take any 
point $x_1\in\Sigma_{X_1}$ and let $w\in\Sigma_{W_2}$ be the image 
of $x_1$ under the composite map 
$\Sigma_{X_1}\overset{\phi}\to{\to}\Sigma_{X_2}\to\Sigma_{W_2}$, 
where the second map arises from the cover $X_2\to W_2$ corresponding 
to the extension $L_2/F_2$. Now, by the Riemann-Roch theorem, 
there exists a function $f\in F_2$ having a pole at $w$. By the 
order-preserving property of $\iota$, the function 
$\iota(f)\in K_1$ must have a pole at $x_1$. 
This contradicts the definition of $F_2$. Therefore, 
we must have $F_2=k_2$, as desired. 

Next, we shall prove that $\phi:\Sigma_{X_1}\to\Sigma_{X_2}$ is 
surjective. Indeed, suppose otherwise and take 
$x_2\in\Sigma_{X_2}\smallsetminus\phi(\Sigma_{X_1})\neq\emptyset$. 
By the Riemann-Roch theorem, there exists a non-constant function 
$f\in K_2$ such that the pole divisor of $f$
is supported on $x_2\in \Sigma _{X_2}$. Then, by the order-preserving 
property of $\iota$, the function $\iota(f)\in K_1$ admits no poles, hence 
$\iota(f)\in k_1$. As $\iota^{-1}(k_1)=k_2$, we thus have 
$f\in k_2$, which is absurd. 

The rest of the proof is similar to the proof of 
Proposition 4.4 
in [ST], 
where $\phi$ is a bijection. 
We shall first prove that $\iota$ restricted to $k_2$
is additive. 
Again by the Riemann-Roch theorem, there exists a non-constant function 
$f\in K_2$ such that the pole divisor $\div (f)_{\infty}$ of $f$
is supported on a unique point $x_2\in \Sigma _{X_2}$: 
$\div (f)_{\infty}=nx_2$, $n>0$.  
For a non-zero constant $\alpha \in k_2$ we shall 
analyze the divisor of the function
$\iota (f+\alpha)-\iota (f)$. We claim that 
$\Supp \div(\iota (f+\alpha)-\iota (f))\subset \phi ^{-1}(x_2)$. 
Indeed, if $y_1\in \Sigma _{X_1}$ is such that $y_2\defeq \phi (y_1)\neq x_2$, then
$\ord _{y_1}(\iota (f+\alpha))\ge 0$, and $\ord _{y_1}(\iota (f))\ge 0$.
Moreover, $\iota (f+\alpha)(y_1)\neq \iota (f)(y_1)$, as 
follows from the value-preserving assumption, 
since $(f+\alpha)(y_1)\neq f(y_1)$. Thus,
$y_1\notin \Supp \div(\iota (f+\alpha)-\iota (f))$ and our claim follows.
Further, 
if $x_1\in \phi ^{-1}(x_2)$ is a pole of $\iota (f+\alpha)-\iota (f)$, 
we have $\vert \ord _{x_1}(\iota (f+\alpha)-\iota (f)) \vert \le ne_{x_1 x_2}$. 
We deduce easily from this that there 
are only finitely many possibilities for the divisor
$\div(\iota (f+\alpha)-\iota (f))$. Since $k_2$ is infinite 
(as $X_2(k_2)$ is infinite), there exists an 
infinite subset $A\subset k_2^{\times}$ such that 
$\div(\iota (f+\alpha)-\iota (f))
$ is constant, 
for all $\alpha \in A$. 

Let $\alpha \neq \beta$ be elements of $A$. Thus, $\div(\iota(f+\alpha)
-\iota (f))=\div(\iota (f+\beta)-\iota (f))$, which implies that
$\frac {\iota (f+\beta)-\iota (f)} {\iota (f+\alpha)-\iota (f)}=c \in k_1^{\times}$.
Observe that $\iota (f+\alpha)-\iota (f)\neq 0$, 
by 
the injectivity of $\iota$. 
Further, $c=\frac 
{\iota (\beta)}{\iota (\alpha)}$, as is easily seen by evaluating the function
$\frac {\iota (f+\beta)-\iota (f)} {\iota (f+\alpha)-\iota (f)}$ at a zero 
of the non-constant function $\iota(f)$. Thus, we have the equality
$\iota (\beta) (\iota (f+\alpha)-\iota (f))=\iota (\alpha) (\iota (f+\beta)-
\iota (f))$, which is equivalent to
$$\iota (f)(\iota (\alpha)-\iota (\beta))=\iota (\alpha)\iota (f+\beta)-
\iota (\beta)\iota (f+\alpha)\tag {$*$}$$
Let $g\defeq \frac {\beta(f+\alpha)}{(\alpha-\beta)f}=\frac 
{\beta (1+\alpha f^{-1})}{(\alpha-\beta)}$. Note that $g$ is a 
non-constant function, since $f$ is non-constant. 
We have $g+1=\frac {\beta(f+\alpha)}{(\alpha-\beta)f}
+\frac {(\alpha-\beta)f}{(\alpha-\beta)f}=\frac {\beta \alpha +\alpha f}
{\alpha f-\beta f}=\frac {\alpha(\beta+f)}{(\alpha-\beta)f}$. 
By dividing the above equality ($*$) by 
$\iota (\alpha-\beta) \iota (f)\neq 0$, we obtain 
$\frac {\iota (\alpha)-\iota (\beta)}{\iota (\alpha-\beta)}
=\frac {\iota (\alpha) \iota (f+\beta)-\iota (\beta)\iota (\alpha+f)}{\iota (\alpha-\beta) 
\iota (f)}$. 
Thus, $\frac {\iota (\alpha)-\iota (\beta)}
{\iota (\alpha-\beta)}=\frac {\iota (\alpha) \iota (f+\beta)}
{\iota (\alpha-\beta) \iota (f)}-\frac {\iota (\beta)\iota (\alpha+f)}
{\iota (\alpha-\beta) \iota (f)}$, which equals $\iota (g+1)-\iota (g)$.
Further, $\frac {\iota (\alpha)-\iota (\beta)}{\iota (\alpha-\beta)}=1$, 
as follows by evaluating the function $\iota (g+1)-\iota (g)$ at a zero 
of the non-constant function $\iota (g)$. Thus, 
$$\iota (g+1)=\iota (g)+1.$$ 

Take any $\zeta\in k_2$. Then, evaluating this equation at a zero of 
$\iota(g-\zeta)$, we obtain 
$\iota (\zeta+1)=\iota(\zeta)+1.$ 
Now, for any $\xi,\eta\in k_2$, we have 
$\iota(\xi+\eta)=\iota(\xi)+\iota(\eta).$ 
Indeed, if $\eta=0$, this follows from $\iota(0)=0$. 
If $\eta\neq 0$, we have 
$$\iota(\xi+\eta)
=\iota\left(\frac{\xi}{\eta}+1\right)\iota(\eta)
=\left(\iota\left(\frac{\xi}{\eta}\right)+1\right)\iota(\eta)=
\iota(\xi)+\iota(\eta).$$
Thus, $\iota|_{k_2}$ is additive. 

{}From this it follows easily that $\iota $ itself is additive. 
Indeed, let $h$ and $l$ be any elements of $K_2$, and let us 
prove $\iota(h+l)=\iota(h)+\iota(l)$.  
Take any $x_2\in X_2(k_2)$ which is neither a pole of $h$ nor 
a pole of $l$. Then, evaluating at any $x_1\in \phi^{-1}(x_2)$, 
we obtain
$$
\align
(\iota (h+l))(x_1)
&=\iota_{x_1x_2}((h+l)(x_2)) \\
&=\iota_{x_1x_2}(h(x_2)+l(x_2)) \\
&=\iota(h(x_2)+l(x_2)) \\
&=\iota(h(x_2))+\iota(l(x_2)) \\
&=\iota_{x_1x_2}(h(x_2))+\iota_{x_1x_2}(l(x_2)) \\
&=(\iota (h))(x_1)+(\iota (l))(x_1) 
\\&
=(\iota (h)+\iota (l))(x_1) 
\endalign
$$
where the first and the sixth equalities follow from the value-preserving 
property, the second and the last equalities follow from the additivity of the 
evaluation maps, the third and the fifth equalities follow from 
the value-preserving property and the fact that 
$h(x_2),l(x_2)\in k_2$ (as $x_2\in X_2(k_2)$), and the fourth equality follows from 
the additivity of $\iota|_{k_2}$. Now, since 
there are infinitely many such $x_1$ by assumption, 
the equality $\iota (h+l)=\iota (h)+\iota (l)$ must hold. 
Thus, the assertion is proved. 
\qed
\enddemo

\bigskip
$$\text{References.}$$
\noindent
[Chevalley] Chevalley, C., 
Deux th\'eor\`emes d'arithm\'etique. J. Math. Soc. Japan, 3 (1951), 
36--44. 

\noindent
[EK] Engler, A.~J. and Koenigsmann, J., Abelian subgroups of pro-$p$ Galois 
groups.  Trans. Amer. Math. Soc.  350  (1998),  no. 6, 2473--2485.

\noindent
[EN] Engler, A.~J. and  Nogueira, J.~B., Maximal abelian normal subgroups of 
Galois pro-$2$-groups.  J. Algebra  166  (1994),  no. 3, 481--505.

\noindent
[FJ] Fried, M.~D. and  Jarden, M., Field arithmetic. Ergebnisse der 
Mathematik und ihrer Grenzgebiete (3), 
11. Springer-Verlag, Berlin, 1986.

\noindent
[Harbater] Harbater, D., Fundamental groups and embedding problems in 
characteristic $p$.  Recent developments in the inverse Galois problem 
(Seattle, 
1993),  353--369, Contemp. Math., 186, Amer. 
Math. Soc., Providence, 
1995.

\noindent
[Koenigsmann] Koenigsmann, J., Encoding valuations in absolute Galois 
groups.  Valuation theory and its applications, Vol. II 
(Saskatoon, 
1999),  107--132, Fields Inst. Commun., 33, 
Amer. Math. Soc., Providence, 
2003.

\noindent
[Mochizuki1] Mochizuki, S., The local pro-$p$ anabelian geometry of curves.  
Invent. Math.  138  (1999),  no. 2, 319--423. 

\noindent
[Mochizuki2] Mochizuki, S., Topics surrounding the anabelian geometry of 
hyperbolic curves. Galois groups and fundamental groups, 119--165, 
Math. Sci. Res. Inst. Publ., 41, Cambridge Univ. Press, Cambridge, 2003. 


\noindent
[Mochizuki3] Mochizuki, S., Galois sections in absolute anabelian geometry.
Nagoya Math. J. Vol. 179 (2005), 17-45.

\noindent
[Mochizuki4] Mochizuki, S., Absolute anabelian cuspidalizations of proper hyperbolic curves. 
J. Math. Kyoto Univ. 47 (2007), no. 3, 451--539.


\noindent
[Neukirch1]
Neukirch, J. Kennzeichnung der $p$-adischen und der endlichen algebraischen Zahlk\"orper. 
Invent. Math. 6 (1969), 296--314.

\noindent
[Neukirch2] Neukirch, J., Kennzeichnung der endlich-algebraischen Zahlk\"orper 
durch die Galoisgruppe der maximal aufl\"osbaren Erweiterungen. 
J. Reine Angew. Math.  238  (1969), 135--147. 

\noindent
[Pop1] Pop, F., On Grothendieck's conjecture of birational anabelian 
geometry.  Ann. of Math. (2)  139  (1994),  no. 1, 145--182.

\noindent
[Pop2] Pop, F., \'Etale Galois covers of affine smooth curves. 
The geometric case of a conjecture of Shafarevich. On Abhyankar's 
conjecture.  Invent. Math.  120  (1995),  no. 3, 555--578.

\noindent
[Pop3] Pop, F., The birational anabelian conjecture -- revisited --, preprint, 
available in 
http://www.math.upenn.edu/$\sim$pop/Research/Papers.html

\noindent
[Rosen] Rosen, M., Number Theory in Function Fields, 
Graduate Texts in Mathematics, 210, Springer-Verlag, New York, 2002.

\noindent
[ST] Sa\"\i di, M. and Tamagawa, A., A prime-to-$p$ version
of the Grothendieck anabelian conjecture for hyperbolic curves
over finite fields of characteristic $p>0$, 
Publ. RIMS 45 (2009), no.1, 135--186. 

\noindent
[SGA-1] Grothendieck, A. and Mme. Raynaud, M., S\'eminaire de G\'eom\'etrie
Alg\'ebrique du Bois Marie 1960/1961, Rev\^etement Etales et Groupe 
Fondamantal (SGA-1), Lecture Notes in Mathematics, 224, Springer-Verlag,
Berlin/Heidelberg/New York, 1971.

\noindent
[Serre1] Serre, J.-P., 
Cohomologie galoisienne, 
Seconde edition. 
Lecture Notes in Mathematics 5, Springer-Verlag, Berlin-Heidelberg-New York, 1962/1963. 

\noindent
[Serre2] Serre, J.-P., Local class field theory. Algebraic Number 
Theory (Proc. Instructional Conf., Brighton, 1965),  128--161, Thompson, 
Washington, D.C., 1967.

\noindent
[Serre3] Serre, J.-P., Corps locaux. Hermann, Paris, 1968.

\noindent
[Szamuely] Szamuely, T., Groupes de Galois de corps de type fini (d'apr\`es 
Pop). Ast\'erisque  No. 294 (2004), ix, 403--431. 

\noindent
[Tamagawa] Tamagawa, A., The Grothendieck conjecture for affine curves.
Compositio Mathematica 109 (1997), 135--194.


\noindent
[Uchida1] Uchida, K., Isomorphisms of Galois groups. 
J. Math. Soc. Japan 28 (1976), no. 4, 617--620. 

\noindent
[Uchida2] Uchida, K., Isomorphisms of Galois groups of algebraic function 
fields.  Ann. Math. (2)  106  (1977), no. 3, 589--598. 

\noindent
[Uchida3] Uchida, K., Homomorphisms of Galois groups of solvably closed 
Galois extensions. J. Math. Soc. Japan  33  (1981), no. 4, 595--604.

\bigskip

\noindent
Mohamed Sa\"\i di

\noindent
School of Engineering, Computer Science, and Mathematics

\noindent
University of Exeter

\noindent
Harrison Building

\noindent
North Park Road

\noindent
EXETER EX4 4QF 

\noindent
United Kingdom

\noindent
M.Saidi\@exeter.ac.uk

\bigskip

\noindent
Akio Tamagawa

\noindent
Research Institute for Mathematical Sciences

\noindent
Kyoto University

\noindent
KYOTO 606-8502

\noindent
Japan

\noindent
tamagawa\@kurims.kyoto-u.ac.jp

\enddocument